\documentclass[12pt]{amsart}
\usepackage{epsfig}
\usepackage{mathtools}
\usepackage[subnum]{cases}
\usepackage{enumerate}
\usepackage{bm}
\usepackage{amsmath,amssymb,mathrsfs,amsthm}
\usepackage{color}
\usepackage{verbatim}
\usepackage{pgf,tikz,pgfplots,tikz-3dplot,tkz-euclide}
\usetikzlibrary{intersections}
\pgfplotsset{compat=1.14}
\usetikzlibrary{arrows, shapes.geometric}
\usepackage{graphicx}
\graphicspath{ {./images/} }
\usepackage{hyperref}

\usepackage[hang]{footmisc}

\newtheorem{theorem}{Theorem}[section]
\newtheorem{definition}[theorem]{Definition}
\newtheorem{lemma}[theorem]{Lemma}
\newtheorem{corollary}[theorem]{Corollary}
\newtheorem{prop}[theorem]{Proposition}
\newtheorem{remark}[theorem]{Remark}

\newtheorem*{property*}{Property}

\newcommand{\R}{\mathbb{R}}
\newcommand{\Sp}{\mathbb{S}^{2}}

\DeclareMathOperator{\vol}{vol}
\DeclareMathOperator{\area}{area}

\DeclareMathOperator{\conv}{conv}
\DeclareMathOperator{\diag}{diag}
\DeclareMathOperator{\sgn}{sgn}

\definecolor{electricviolet}{rgb}{0.56, 0.0, 1.0}
\begin{document}
\setcounter{footnote}{0}

\title[On maximum volume polytopes inscribed in the sphere]{On maximum volume polytopes inscribed in the sphere}

    \author[S. Hoehner]{Steven Hoehner}
    \date{\today}


\begin{abstract}
We combine degree sequence bounds with a supporting tangent plane method to study maximum volume polytopes inscribed in the sphere $\mathbb S^2$. Building on the vertex-star estimates and combinatorial reductions in Wimmer's dissertation, we determine the exact maximum for ten vertices and the unique maximizing combinatorial type for eleven vertices. Moreover, using explicit computational certificates we identify a unique possible maximizing type for polytopes with thirteen and fifteen vertices, and two possible types for sixteen vertices; this approach also recovers Wimmer's fourteen-vertex combinatorial theorem. Finally, we derive a facial inequality which holds for every degree pattern, and prove that the maximum volume for polytopes with thirty-two vertices is $\frac{2}{3}\sqrt{10(5-\sqrt{5})}$, attained precisely by an inscribed realization of the pentakis dodecahedron, which is the convex hull of the regular icosahedron together with its normalized facet centers. Wimmer had already determined this configuration and its optimality within its own combinatorial type. We also obtain a global degree sequence envelope and a universal affine volume bound.
\end{abstract}

\maketitle

\section{Introduction and main results}\label{sec:introduction}

Let $\mathcal{P}_N$ denote the class of convex polytopes with at most $N$ vertices inscribed in the Euclidean unit sphere $\Sp$. A classical problem
in discrete geometry asks for
\begin{equation}\label{eq:max-volume-problem}
    \max\{\vol(P):P\in\mathcal{P}_N\},
\end{equation}
and for a characterization of the maximizing polytopes. Equivalently, we may view this as the problem of finding the best approximation of the three-dimensional Euclidean ball by an inscribed polytope with a prescribed number of vertices, where the accuracy of the approximation is measured by the
volume difference.

For small values of $N$, this problem has been solved by exploiting the combinatorial structure of the polytope together with local geometric inequalities; see, e.g., \cite{BermanHanes1970,Toth,Hoehner-Ledford-2026,Wimmer1998}. As $N$ increases, however, the number of possible combinatorial types grows rapidly. In general, we must first rule out the competing combinatorial classes, and then determine the best realization within each remaining class.

In this direction, Wimmer's dissertation \cite{Wimmer1998} developed substantial earlier progress on this problem. In particular, he proves a vertex-star volume bound, derives local exclusions, and reduces the ten- and eleven-vertex problems to two possible combinatorial types in each case. He also identifies the unique possible maximizing type for fourteen vertices. In \cite[Chapter IV]{Wimmer1998}, Wimmer  discusses the thirteen-vertex candidates, obtains bounds in terms of degree-5 and degree-6 vertices, and, for $N=32$, he determines the volume maximizer within a special combinatorial type. In \cite[Chapter VI]{Wimmer1998}, he supplies numerical configurations for the larger values of $N$ considered in the present paper. We will use these results explicitly in our proofs below. 

In this paper, we combine degree sequence bounds with a supporting ``tangent plane method'' for facial cone volumes. We will employ Wimmer's reductions for small cases to shorten the combinatorial part of the ten- and eleven-vertex arguments. The remaining ten-vertex type is optimized without any symmetry assumptions. For thirteen and fifteen vertices, the methods in the present paper identify a unique possible maximizing combinatorial type; for sixteen vertices, our methods leave two possible types. Finally, a universal choice of supporting planes proves a sharp global theorem for polytopes with thirty-two vertices.

\vspace{2mm}

Our first main result concerns polytopes with ten vertices.
\begin{theorem}\label{mainThm}
Let $P\in\mathcal{P}_{10}$. Then
\begin{equation}
\begin{aligned}
\vol(P)
\leq V_*&:=\frac{4}{3}(1-h_*^2)\left(1+(1+\sqrt{2})h_*\right)\\
&=\frac{16}{81}\left(3+\sqrt{2}+(3-\sqrt{2})\sqrt{10+6\sqrt{2}}\right)= 2.2187111315\ldots
\end{aligned}
\end{equation}
where $h_*=\frac{\sqrt{10+6\sqrt{2}}-1}{3(1+\sqrt{2})}$. Equality holds if and only if, up to an orthogonal transformation, $P$ is the gyroelongated square bipyramid with vertices
\begin{align*}
a_+=(0,0,1), \quad a_-=(0,0,-1),
\end{align*}
and
\begin{equation}
    v_j=\left(r_*\cos\left(\tfrac{\pi}{8}+\tfrac{j\pi}{4}\right),r_*\sin\left(\tfrac{\pi}{8}+\tfrac{j\pi}{4}\right),(-1)^j h_*\right),\quad j\in\{0,1,\ldots,7\},
\end{equation}
where $r_*=\sqrt{1-h_*^2}$.
\end{theorem}

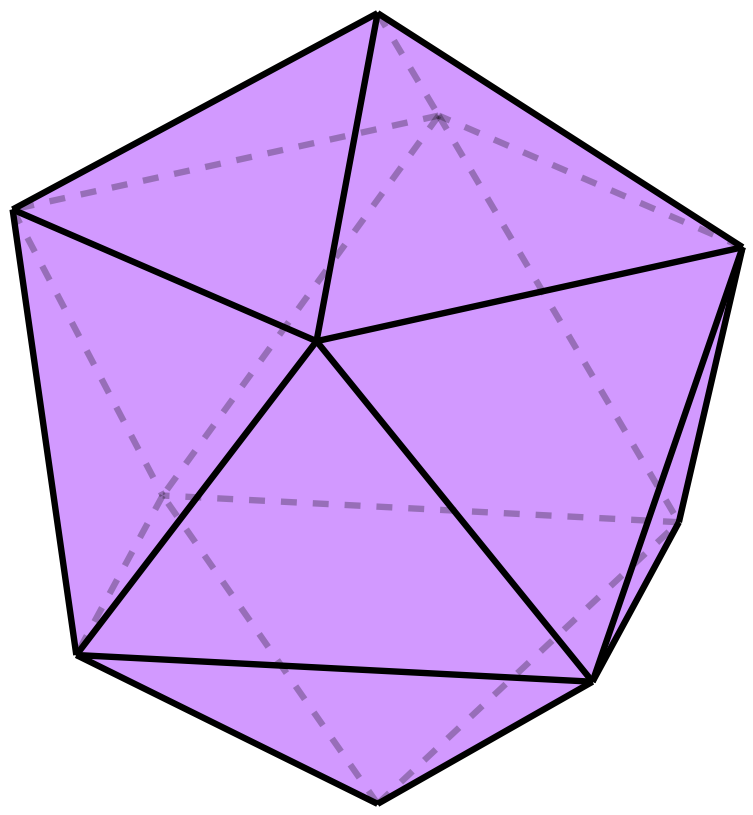
\begin{figure}
\begin{center}
\tdplotsetmaincoords{72}{122}

\begin{tikzpicture}[
    scale=2.75,
    line join=bevel,
    tdplot_main_coords
]

\def\h{0.455559157234444}
\def\a{0.822442652610432}
\def\b{0.340666900985345}

\coordinate (T) at (0,0,1);
\coordinate (D) at (0,0,-1);

\coordinate (A1) at ( \a, \b, \h);
\coordinate (A2) at (-\b, \a, \h);
\coordinate (A3) at (-\a,-\b, \h);
\coordinate (A4) at ( \b,-\a, \h);

\coordinate (B1) at ( \a,-\b,-\h);
\coordinate (B2) at ( \b, \a,-\h);
\coordinate (B3) at (-\a, \b,-\h);
\coordinate (B4) at (-\b,-\a,-\h);

\draw[fill opacity=0.4,fill=electricviolet]
    (T)--(A1)--(A2)--cycle;

\draw[fill opacity=0.4,fill=electricviolet]
    (T)--(A4)--(A1)--cycle;

\draw[fill opacity=0.4,fill=electricviolet]
    (A1)--(A2)--(B2)--cycle;

\draw[fill opacity=0.4,fill=electricviolet]
    (A1)--(B1)--(B2)--cycle;

\draw[fill opacity=0.4,fill=electricviolet]
    (A2)--(B2)--(B3)--cycle;

\draw[fill opacity=0.4,fill=electricviolet]
    (A4)--(A1)--(B1)--cycle;

\draw[fill opacity=0.4,fill=electricviolet]
    (D)--(B1)--(B2)--cycle;

\draw[very thick] (T)--(A1);
\draw[very thick] (T)--(A2);
\draw[very thick] (T)--(A4);

\draw[very thick] (A1)--(A2);
\draw[very thick] (A4)--(A1);

\draw[very thick] (A1)--(B1);
\draw[very thick] (A1)--(B2);

\draw[very thick] (A2)--(B2);
\draw[very thick] (A2)--(B3);

\draw[very thick] (A4)--(B1);

\draw[very thick] (B1)--(B2);
\draw[very thick] (B2)--(B3);

\draw[very thick] (D)--(B1);
\draw[very thick] (D)--(B2);

\draw[very thick,dashed,opacity=0.28] (T)--(A3);

\draw[very thick,dashed,opacity=0.28] (A2)--(A3);
\draw[very thick,dashed,opacity=0.28] (A3)--(A4);

\draw[very thick,dashed,opacity=0.28] (A3)--(B3);
\draw[very thick,dashed,opacity=0.28] (A3)--(B4);
\draw[very thick,dashed,opacity=0.28] (A4)--(B4);

\draw[very thick,dashed,opacity=0.28] (B3)--(B4);
\draw[very thick,dashed,opacity=0.28] (B4)--(B1);

\draw[very thick,dashed,opacity=0.28] (D)--(B3);
\draw[very thick,dashed,opacity=0.28] (D)--(B4);

\end{tikzpicture}
\end{center}
\caption{The global volume maximizer in $\mathcal{P}_{10}$ is a gyroelongated square bipyramid.}
    \label{fig:gyroelongated-square-bipyramid}
\end{figure}

The four even-indexed $v_j$'s form the upper square, the four odd-indexed $v_j$'s form the lower square, and the two squares are parallel and twisted by $45^\circ$. The other two vertices, $a_+$ and $a_-$, are the apices which cap the squares. See Figure \ref{fig:gyroelongated-square-bipyramid} for an illustration.

For polytopes with eleven vertices, we determine the maximizing combinatorial type.

\begin{theorem}\label{mainThm-N11}
Let $P_{11}^*\in\mathcal{P}_{11}$ have maximum volume. Then $P_{11}^*$ is combinatorially equivalent to the edge-contracted icosahedron. In particular, its degree sequence is $(6,5^8,4^2)$.
\end{theorem}

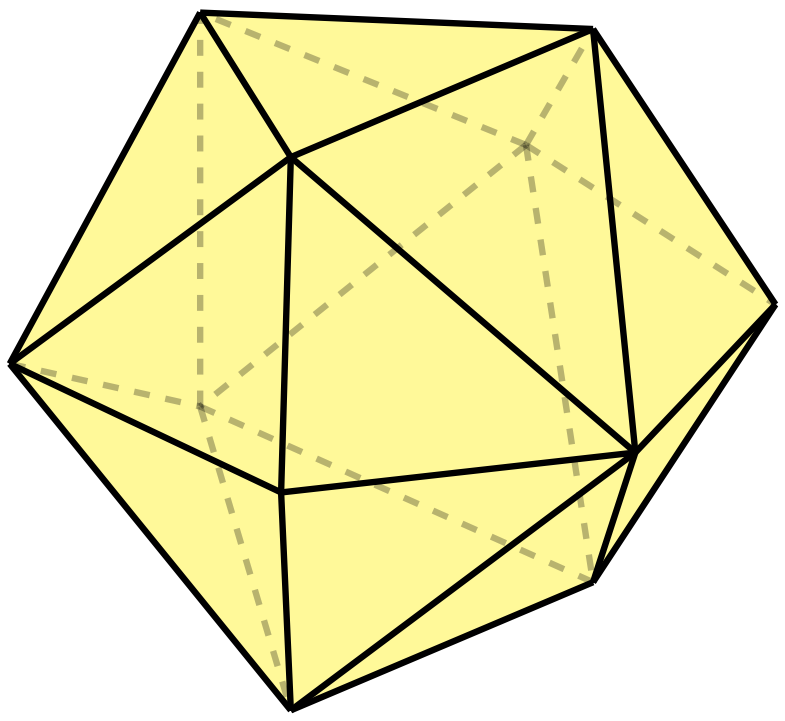
\begin{figure}[ht]
\begin{center}
\begin{tikzpicture}[
    scale=0.8,
    line join=bevel
]

\coordinate (p)   at ( 1.9990,-0.9500);
\coordinate (ap)  at ( 3.1578, 0.2745);
\coordinate (am)  at (-0.9300,-1.2800); 
\coordinate (bp)  at ( 1.0949, 1.5931);
\coordinate (bm)  at (-3.1739,-0.2162);
\coordinate (cp)  at (-1.5992, 2.6875);
\coordinate (cm)  at (-1.5992,-0.5693);
\coordinate (dpp) at ( 1.6492, 2.5516);
\coordinate (dpm) at ( 1.6492,-2.0221);
\coordinate (dmp) at (-0.8496, 1.4925);
\coordinate (dmm) at (-0.8496,-3.0812);


\draw[fill opacity=0.4,fill=yellow]
    (p)--(ap)--(dpp)--cycle;

\draw[fill opacity=0.4,fill=yellow]
    (p)--(ap)--(dpm)--cycle;

\draw[fill opacity=0.4,fill=yellow]
    (p)--(am)--(dmp)--cycle;

\draw[fill opacity=0.4,fill=yellow]
    (p)--(am)--(dmm)--cycle;

\draw[fill opacity=0.4,fill=yellow]
    (p)--(dpp)--(dmp)--cycle;

\draw[fill opacity=0.4,fill=yellow]
    (p)--(dpm)--(dmm)--cycle;

\draw[fill opacity=0.4,fill=yellow]
    (am)--(bm)--(dmp)--cycle;

\draw[fill opacity=0.4,fill=yellow]
    (am)--(bm)--(dmm)--cycle;

\draw[fill opacity=0.4,fill=yellow]
    (bm)--(cp)--(dmp)--cycle;

\draw[fill opacity=0.4,fill=yellow]
    (cp)--(dpp)--(dmp)--cycle;


\begin{scope}[very thick,dashed,opacity=0.28]
    \draw (bm)--(cm);
    \draw (cm)--(cp);
    \draw (cm)--(dmm);
    \draw (cm)--(dpm);
    \draw (bp)--(cm);
    \draw (bp)--(cp);
    \draw (bp)--(dpp);
    \draw (bp)--(ap);
    \draw (bp)--(dpm);
\end{scope}


\begin{scope}[very thick]
    \draw (am)--(dmp);
    \draw (dmp)--(p);
    \draw (p)--(am);
    \draw (ap)--(dpm);
    \draw (ap)--(p);
    \draw (p)--(dpm);
    \draw (bm)--(cp);
    \draw (bm)--(dmp);
    \draw (cp)--(dmp);
    \draw (am)--(bm);
    \draw (dpm)--(dmm);
    \draw (am)--(dmm);
    \draw (p)--(dmm);
    \draw (bm)--(dmm);
    \draw (dpp)--(dmp);
    \draw (p)--(dpp);
    \draw (ap)--(dpp);
    \draw (cp)--(dpp);
\end{scope}



\end{tikzpicture}
\end{center}
\caption{The combinatorial type of every global volume maximizer in $\mathcal{P}_{11}$ is that of an edge-contracted icosahedron.}
\label{fig:edge-contracted-icosahedron}
\end{figure}

Here an edge-contracted icosahedron means the triangulation obtained by contracting one edge of the icosahedral graph and suppressing the resulting pairs of parallel edges. Note that Theorem \ref{mainThm-N11} does not determine the best geometric realization of this graph; we leave this as an open problem for further research. See Figure~\ref{fig:edge-contracted-icosahedron} for an illustration of the maximizing combinatorial type for $N=11$.

Our next result provides combinatorial reductions for the cases $N=13$, 15, and 16. Here and in what follows, $(ijk)$ denotes a triangular facet whose vertices have degrees $i,j,k$, and $m(ijk)$ means there are $m$ facets of that type.

\begin{theorem}\label{thm:larger-combinatorial-reductions}
Let $P_N^*\in\mathcal{P}_N$ have maximum volume.
\begin{itemize}
    \item[(i)] If $N=13$, then $P_N^*$ has the unique combinatorial type with degree sequence $(6^2,5^{10},4)$. Its facial degree multiset is $4(456)+10(555)+8(556)$.
    \item[(ii)] If $N=15$, then $P_N^*$ has the unique combinatorial type with degree sequence $(6^3,5^{12})$. Its facial degree multiset is $8(555)+18(556)$.
    \item[(iii)] If $N=16$, then $P_N^*$ has one of the two combinatorial types with degree sequence $(6^4,5^{12})$. Their facial degree multisets are
    \[
        8(555)+16(556)+4(566)\quad\text{and}\quad 4(555)+24(556).
    \]
\end{itemize}
\end{theorem}
 The combinatorial reductions use exhaustive graph lists and explicit interval certificates, supplied in the accompanying Zenodo archive; see Appendix \ref{sec:certificates}. The geometric optimization of the types in Theorem \ref{thm:larger-combinatorial-reductions} is not completed here; this is also left as an open problem for further research. Let us note here that our computations for this result also recover Wimmer's fourteen-vertex  theorem from \cite{Wimmer1998} in which he determined the maximizer's combinatorial type; see Remark \ref{rem:N14} below.
 
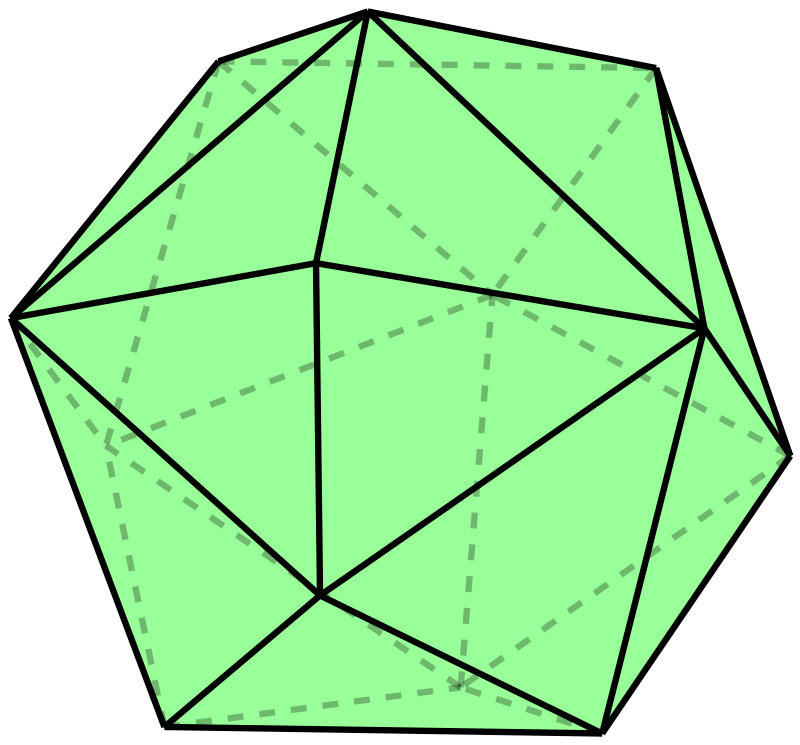
\begin{figure}[ht]
\begin{center}
\tdplotsetmaincoords{72}{122}

\begin{tikzpicture}[
    scale=2.75,
    line join=bevel,
    tdplot_main_coords
]

\coordinate (V1) at ( 0.7844132594, 0.2457500141, 0.5694758723);
\coordinate (V2) at ( 0.9150283136, 0.3387736150,-0.2189877238);
\coordinate (V3) at ( 0.1207139198, 0.9326800220, 0.3399060549);
\coordinate (V4) at ( 0.1781027249, 0.6785691290,-0.7126172581);
\coordinate (V5) at ( 0.5411976825,-0.3361653288,-0.7707774907);
\coordinate (V6) at (-0.7060395651, 0.6597483311,-0.2573718557);
\coordinate (V7) at (-0.3265337647,-0.7248352052, 0.6066214848);
\coordinate (V8) at ( 0.1507412569,-0.0033237372, 0.9885676640);
\coordinate (V9) at (-0.6887741547, 0.2897382441, 0.6645614447);
\coordinate (V10) at ( 0.6942252021,-0.6760338142, 0.2470418000);
\coordinate (V11) at (-0.4288319131,-0.1021342874,-0.8975922112);
\coordinate (V12) at (-0.9397406787,-0.3311439317,-0.0850361886);
\coordinate (V13) at (-0.1385840045,-0.9269212541,-0.3487283504);

\draw[fill opacity=0.4,fill=green]
    (V3)--(V6)--(V9)--cycle;

\draw[fill opacity=0.4,fill=green]
    (V7)--(V10)--(V8)--cycle;

\draw[fill opacity=0.4,fill=green]
    (V3)--(V4)--(V6)--cycle;

\draw[fill opacity=0.4,fill=green]
    (V3)--(V9)--(V8)--cycle;

\draw[fill opacity=0.4,fill=green]
    (V2)--(V5)--(V4)--cycle;

\draw[fill opacity=0.4,fill=green]
    (V2)--(V10)--(V5)--cycle;

\draw[fill opacity=0.4,fill=green]
    (V1)--(V8)--(V10)--cycle;

\draw[fill opacity=0.4,fill=green]
    (V2)--(V4)--(V3)--cycle;

\draw[fill opacity=0.4,fill=green]
    (V1)--(V3)--(V8)--cycle;

\draw[fill opacity=0.4,fill=green]
    (V1)--(V10)--(V2)--cycle;

\draw[fill opacity=0.4,fill=green]
    (V1)--(V2)--(V3)--cycle;

\draw[very thick,dashed,opacity=0.28] (V4)--(V11);
\draw[very thick,dashed,opacity=0.28] (V5)--(V11);
\draw[very thick,dashed,opacity=0.28] (V5)--(V13);
\draw[very thick,dashed,opacity=0.28] (V6)--(V11);
\draw[very thick,dashed,opacity=0.28] (V6)--(V12);
\draw[very thick,dashed,opacity=0.28] (V7)--(V9);
\draw[very thick,dashed,opacity=0.28] (V7)--(V12);
\draw[very thick,dashed,opacity=0.28] (V7)--(V13);
\draw[very thick,dashed,opacity=0.28] (V9)--(V12);
\draw[very thick,dashed,opacity=0.28] (V10)--(V13);
\draw[very thick,dashed,opacity=0.28] (V11)--(V12);
\draw[very thick,dashed,opacity=0.28] (V11)--(V13);
\draw[very thick,dashed,opacity=0.28] (V12)--(V13);

\draw[very thick] (V1)--(V2);
\draw[very thick] (V1)--(V3);
\draw[very thick] (V1)--(V8);
\draw[very thick] (V1)--(V10);
\draw[very thick] (V2)--(V3);
\draw[very thick] (V2)--(V4);
\draw[very thick] (V2)--(V5);
\draw[very thick] (V2)--(V10);
\draw[very thick] (V3)--(V4);
\draw[very thick] (V3)--(V6);
\draw[very thick] (V3)--(V8);
\draw[very thick] (V3)--(V9);
\draw[very thick] (V4)--(V5);
\draw[very thick] (V4)--(V6);
\draw[very thick] (V5)--(V10);
\draw[very thick] (V6)--(V9);
\draw[very thick] (V7)--(V8);
\draw[very thick] (V7)--(V10);
\draw[very thick] (V8)--(V9);
\draw[very thick] (V8)--(V10);

\end{tikzpicture}
\end{center}
\caption{Every global volume maximizer in $\mathcal{P}_{13}$ is combinatorially equivalent to the polytope shown, whose degree sequence is $(6^2,5^{10},4)$.}
\label{fig:N13-maximizing-type}
\end{figure}

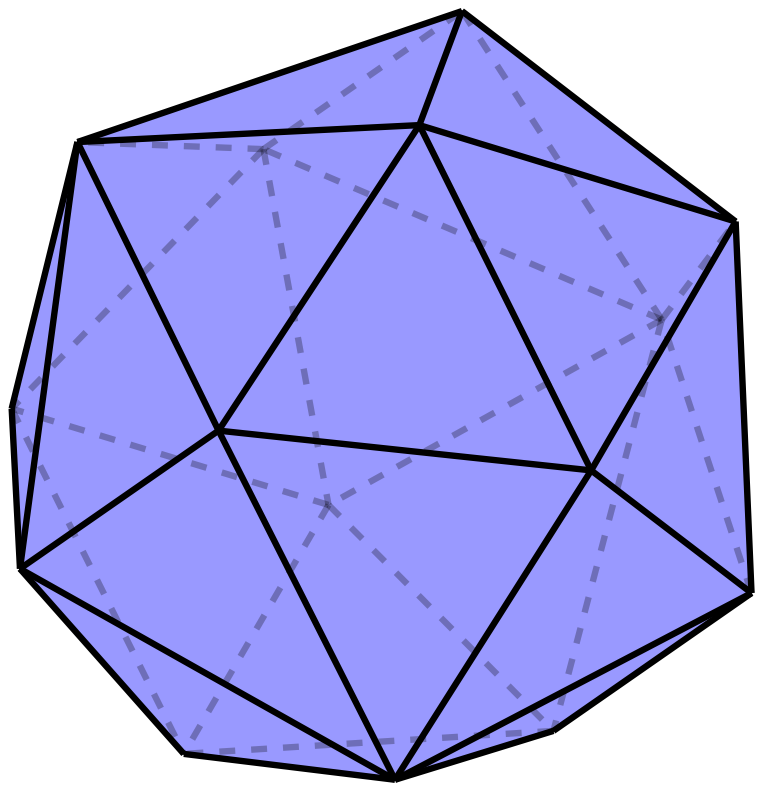
\begin{figure}[ht]
\begin{center}
\tdplotsetmaincoords{72}{122}

\begin{tikzpicture}[
    scale=2.75,
    line join=bevel,
    tdplot_main_coords
]

\coordinate (V1) at (-0.2941446969, 0.8755062771,-0.3833636082);
\coordinate (V2) at (-0.9815230391, 0.1913422526, 0.0008161094);
\coordinate (V3) at (-0.6014951000,-0.7009363160, 0.3832648243);
\coordinate (V4) at (-0.4757729424,-0.0593412048, 0.8775640881);
\coordinate (V5) at (-0.6605595067,-0.5533123111,-0.5074511056);
\coordinate (V6) at (-0.4039333072, 0.7618270086, 0.5064163231);
\coordinate (V7) at ( 0.3399671402, 0.3294699909, 0.8808358920);
\coordinate (V8) at (-0.4194106322, 0.2353169081,-0.8767671723);
\coordinate (V9) at ( 0.8326080552,-0.4954942001,-0.2474860078);
\coordinate (V10) at ( 0.0827894281,-0.9875507812,-0.1337511313);
\coordinate (V11) at ( 0.1910500644,-0.4317190080,-0.8815432893);
\coordinate (V12) at ( 0.3910274019,-0.6083923113, 0.6906202767);
\coordinate (V13) at ( 0.4477680852, 0.8842386103, 0.1327622761);
\coordinate (V14) at ( 0.5906791009, 0.4170690233,-0.6907616300);
\coordinate (V15) at ( 0.9579731705, 0.1461465335, 0.2468371837);

\draw[fill opacity=0.4,fill=blue]
    (V8)--(V1)--(V14)--cycle;

\draw[fill opacity=0.4,fill=blue]
    (V9)--(V12)--(V10)--cycle;

\draw[fill opacity=0.4,fill=blue]
    (V9)--(V11)--(V14)--cycle;

\draw[fill opacity=0.4,fill=blue]
    (V6)--(V4)--(V7)--cycle;

\draw[fill opacity=0.4,fill=blue]
    (V4)--(V12)--(V7)--cycle;

\draw[fill opacity=0.4,fill=blue]
    (V13)--(V1)--(V6)--cycle;

\draw[fill opacity=0.4,fill=blue]
    (V13)--(V14)--(V1)--cycle;

\draw[fill opacity=0.4,fill=blue]
    (V15)--(V12)--(V9)--cycle;

\draw[fill opacity=0.4,fill=blue]
    (V15)--(V9)--(V14)--cycle;

\draw[fill opacity=0.4,fill=blue]
    (V13)--(V6)--(V7)--cycle;

\draw[fill opacity=0.4,fill=blue]
    (V15)--(V7)--(V12)--cycle;

\draw[fill opacity=0.4,fill=blue]
    (V15)--(V14)--(V13)--cycle;

\draw[fill opacity=0.4,fill=blue]
    (V15)--(V13)--(V7)--cycle;

\draw[very thick,dashed,opacity=0.28] (V1)--(V2);
\draw[very thick,dashed,opacity=0.28] (V2)--(V3);
\draw[very thick,dashed,opacity=0.28] (V2)--(V4);
\draw[very thick,dashed,opacity=0.28] (V2)--(V5);
\draw[very thick,dashed,opacity=0.28] (V2)--(V6);
\draw[very thick,dashed,opacity=0.28] (V2)--(V8);
\draw[very thick,dashed,opacity=0.28] (V3)--(V4);
\draw[very thick,dashed,opacity=0.28] (V3)--(V5);
\draw[very thick,dashed,opacity=0.28] (V3)--(V10);
\draw[very thick,dashed,opacity=0.28] (V3)--(V12);
\draw[very thick,dashed,opacity=0.28] (V5)--(V8);
\draw[very thick,dashed,opacity=0.28] (V5)--(V10);
\draw[very thick,dashed,opacity=0.28] (V5)--(V11);
\draw[very thick,dashed,opacity=0.28] (V8)--(V11);
\draw[very thick,dashed,opacity=0.28] (V10)--(V11);

\draw[very thick] (V1)--(V6);
\draw[very thick] (V1)--(V8);
\draw[very thick] (V1)--(V13);
\draw[very thick] (V1)--(V14);
\draw[very thick] (V4)--(V6);
\draw[very thick] (V4)--(V7);
\draw[very thick] (V4)--(V12);
\draw[very thick] (V6)--(V7);
\draw[very thick] (V6)--(V13);
\draw[very thick] (V7)--(V12);
\draw[very thick] (V7)--(V13);
\draw[very thick] (V7)--(V15);
\draw[very thick] (V8)--(V14);
\draw[very thick] (V9)--(V10);
\draw[very thick] (V9)--(V11);
\draw[very thick] (V9)--(V12);
\draw[very thick] (V9)--(V14);
\draw[very thick] (V9)--(V15);
\draw[very thick] (V10)--(V12);
\draw[very thick] (V11)--(V14);
\draw[very thick] (V12)--(V15);
\draw[very thick] (V13)--(V14);
\draw[very thick] (V13)--(V15);
\draw[very thick] (V14)--(V15);

\end{tikzpicture}
\end{center}
\caption{Every global volume maximizer in $\mathcal{P}_{15}$ is combinatorially equivalent to the polytope shown, whose degree sequence is $(6^3,5^{12})$.}
\label{fig:N15-maximizing-type}
\end{figure}

\begin{figure}[ht]
\begin{center}
\begin{minipage}[b]{0.48\textwidth}
\centering
\tdplotsetmaincoords{72}{122}

\begin{tikzpicture}[
    scale=2.75,
    line join=bevel,
    tdplot_main_coords
]

\coordinate (A1) at (-0.1686733064, 0.0667776398, 0.9834073736);
\coordinate (A2) at ( 0.1891680674, 0.8083410552, 0.5574945566);
\coordinate (A3) at ( 0.7526702292, 0.0835560572, 0.6530742005);
\coordinate (A4) at (-0.0403325064,-0.7173078318, 0.6955880702);
\coordinate (A5) at (-0.7964244326,-0.4306440701, 0.4245630790);
\coordinate (A6) at (-0.7311207049, 0.4928698251, 0.4717434158);
\coordinate (A7) at (-0.3190171105, 0.9408135448,-0.1144463064);
\coordinate (A8) at ( 0.6156373418, 0.7759343054,-0.1375376933);
\coordinate (A9) at ( 0.9913819056,-0.0667889760,-0.1126993786);
\coordinate (A10) at ( 0.6168684552,-0.7537208431, 0.2266675973);
\coordinate (A11) at (-0.4380081402,-0.8864148499,-0.1497250242);
\coordinate (A12) at (-0.9302994308, 0.0269244874,-0.3658114830);
\coordinate (A13) at (-0.2575188382, 0.5302151300,-0.8078093611);
\coordinate (A14) at ( 0.4860019287, 0.1325495425,-0.8638476394);
\coordinate (A15) at ( 0.4137920825,-0.6921970075,-0.5913031501);
\coordinate (A16) at (-0.3841255404,-0.3109080092,-0.8693582570);

\draw[fill opacity=0.4,fill=red]
    (A2)--(A7)--(A6)--cycle;

\draw[fill opacity=0.4,fill=red]
    (A9)--(A10)--(A15)--cycle;

\draw[fill opacity=0.4,fill=red]
    (A9)--(A15)--(A14)--cycle;

\draw[fill opacity=0.4,fill=red]
    (A1)--(A2)--(A6)--cycle;

\draw[fill opacity=0.4,fill=red]
    (A7)--(A8)--(A13)--cycle;

\draw[fill opacity=0.4,fill=red]
    (A8)--(A14)--(A13)--cycle;

\draw[fill opacity=0.4,fill=red]
    (A3)--(A4)--(A10)--cycle;

\draw[fill opacity=0.4,fill=red]
    (A1)--(A4)--(A3)--cycle;

\draw[fill opacity=0.4,fill=red]
    (A2)--(A8)--(A7)--cycle;

\draw[fill opacity=0.4,fill=red]
    (A8)--(A9)--(A14)--cycle;

\draw[fill opacity=0.4,fill=red]
    (A3)--(A10)--(A9)--cycle;

\draw[fill opacity=0.4,fill=red]
    (A1)--(A3)--(A2)--cycle;

\draw[fill opacity=0.4,fill=red]
    (A3)--(A9)--(A8)--cycle;

\draw[fill opacity=0.4,fill=red]
    (A2)--(A3)--(A8)--cycle;

\draw[very thick,dashed,opacity=0.28] (A1)--(A5);
\draw[very thick,dashed,opacity=0.28] (A4)--(A5);
\draw[very thick,dashed,opacity=0.28] (A4)--(A11);
\draw[very thick,dashed,opacity=0.28] (A5)--(A6);
\draw[very thick,dashed,opacity=0.28] (A5)--(A11);
\draw[very thick,dashed,opacity=0.28] (A5)--(A12);
\draw[very thick,dashed,opacity=0.28] (A6)--(A12);
\draw[very thick,dashed,opacity=0.28] (A7)--(A12);
\draw[very thick,dashed,opacity=0.28] (A10)--(A11);
\draw[very thick,dashed,opacity=0.28] (A11)--(A12);
\draw[very thick,dashed,opacity=0.28] (A11)--(A15);
\draw[very thick,dashed,opacity=0.28] (A11)--(A16);
\draw[very thick,dashed,opacity=0.28] (A12)--(A13);
\draw[very thick,dashed,opacity=0.28] (A12)--(A16);
\draw[very thick,dashed,opacity=0.28] (A13)--(A16);
\draw[very thick,dashed,opacity=0.28] (A14)--(A16);
\draw[very thick,dashed,opacity=0.28] (A15)--(A16);

\draw[very thick] (A1)--(A2);
\draw[very thick] (A1)--(A3);
\draw[very thick] (A1)--(A4);
\draw[very thick] (A1)--(A6);
\draw[very thick] (A2)--(A3);
\draw[very thick] (A2)--(A6);
\draw[very thick] (A2)--(A7);
\draw[very thick] (A2)--(A8);
\draw[very thick] (A3)--(A4);
\draw[very thick] (A3)--(A8);
\draw[very thick] (A3)--(A9);
\draw[very thick] (A3)--(A10);
\draw[very thick] (A4)--(A10);
\draw[very thick] (A6)--(A7);
\draw[very thick] (A7)--(A8);
\draw[very thick] (A7)--(A13);
\draw[very thick] (A8)--(A9);
\draw[very thick] (A8)--(A13);
\draw[very thick] (A8)--(A14);
\draw[very thick] (A9)--(A10);
\draw[very thick] (A9)--(A14);
\draw[very thick] (A9)--(A15);
\draw[very thick] (A10)--(A15);
\draw[very thick] (A13)--(A14);
\draw[very thick] (A14)--(A15);

\end{tikzpicture}
\par\smallskip
{\small (a) $8(555)+16(556)+4(566)$\par}
\end{minipage}\hfill
\begin{minipage}[b]{0.48\textwidth}
\centering
\tdplotsetmaincoords{72}{122}

\begin{tikzpicture}[
    scale=2.75,
    line join=bevel,
    tdplot_main_coords
]

\coordinate (B1) at (-0.1055107172,-0.1532281329, 0.9825419217);
\coordinate (B2) at ( 0.6665153675, 0.7385423479, 0.1015503086);
\coordinate (B3) at ( 0.1286199642, 0.6087384551, 0.7828757233);
\coordinate (B4) at (-0.6592557338,-0.5927827019, 0.4625911216);
\coordinate (B5) at ( 0.1821337229,-0.8375151018, 0.5151657609);
\coordinate (B6) at (-0.7753095394, 0.3722905249, 0.5101910262);
\coordinate (B7) at ( 0.7896072778,-0.0831528292, 0.6079522628);
\coordinate (B8) at ( 0.6860983037,-0.7170236946,-0.1230696512);
\coordinate (B9) at (-0.1747007104, 0.9834290897, 0.0484446827);
\coordinate (B10) at (-0.4051552951,-0.1142509851,-0.9070809772);
\coordinate (B11) at (-0.9436124740,-0.2425396592,-0.2253220200);
\coordinate (B12) at (-0.6930192757, 0.5714550752,-0.4395035615);
\coordinate (B13) at (-0.2399250702,-0.9017477333,-0.3595649958);
\coordinate (B14) at ( 0.9365073947, 0.0970206430,-0.3369582977);
\coordinate (B15) at ( 0.2259615375, 0.6121862421,-0.7577396575);
\coordinate (B16) at ( 0.3830217354,-0.3433436806,-0.8575601828);

\draw[fill opacity=0.4,fill=red]
    (B15)--(B12)--(B9)--cycle;

\draw[fill opacity=0.4,fill=red]
    (B6)--(B1)--(B3)--cycle;

\draw[fill opacity=0.4,fill=red]
    (B16)--(B14)--(B8)--cycle;

\draw[fill opacity=0.4,fill=red]
    (B6)--(B3)--(B9)--cycle;

\draw[fill opacity=0.4,fill=red]
    (B1)--(B5)--(B7)--cycle;

\draw[fill opacity=0.4,fill=red]
    (B8)--(B7)--(B5)--cycle;

\draw[fill opacity=0.4,fill=red]
    (B16)--(B15)--(B14)--cycle;

\draw[fill opacity=0.4,fill=red]
    (B15)--(B9)--(B2)--cycle;

\draw[fill opacity=0.4,fill=red]
    (B1)--(B7)--(B3)--cycle;

\draw[fill opacity=0.4,fill=red]
    (B8)--(B14)--(B7)--cycle;

\draw[fill opacity=0.4,fill=red]
    (B15)--(B2)--(B14)--cycle;

\draw[fill opacity=0.4,fill=red]
    (B2)--(B9)--(B3)--cycle;

\draw[fill opacity=0.4,fill=red]
    (B2)--(B3)--(B7)--cycle;

\draw[fill opacity=0.4,fill=red]
    (B2)--(B7)--(B14)--cycle;

\draw[very thick,dashed,opacity=0.28] (B1)--(B4);
\draw[very thick,dashed,opacity=0.28] (B4)--(B5);
\draw[very thick,dashed,opacity=0.28] (B4)--(B6);
\draw[very thick,dashed,opacity=0.28] (B4)--(B11);
\draw[very thick,dashed,opacity=0.28] (B4)--(B13);
\draw[very thick,dashed,opacity=0.28] (B5)--(B13);
\draw[very thick,dashed,opacity=0.28] (B6)--(B11);
\draw[very thick,dashed,opacity=0.28] (B6)--(B12);
\draw[very thick,dashed,opacity=0.28] (B8)--(B13);
\draw[very thick,dashed,opacity=0.28] (B10)--(B11);
\draw[very thick,dashed,opacity=0.28] (B10)--(B12);
\draw[very thick,dashed,opacity=0.28] (B10)--(B13);
\draw[very thick,dashed,opacity=0.28] (B10)--(B15);
\draw[very thick,dashed,opacity=0.28] (B10)--(B16);
\draw[very thick,dashed,opacity=0.28] (B11)--(B12);
\draw[very thick,dashed,opacity=0.28] (B11)--(B13);
\draw[very thick,dashed,opacity=0.28] (B13)--(B16);

\draw[very thick] (B1)--(B3);
\draw[very thick] (B1)--(B5);
\draw[very thick] (B1)--(B6);
\draw[very thick] (B1)--(B7);
\draw[very thick] (B2)--(B3);
\draw[very thick] (B2)--(B7);
\draw[very thick] (B2)--(B9);
\draw[very thick] (B2)--(B14);
\draw[very thick] (B2)--(B15);
\draw[very thick] (B3)--(B6);
\draw[very thick] (B3)--(B7);
\draw[very thick] (B3)--(B9);
\draw[very thick] (B5)--(B7);
\draw[very thick] (B5)--(B8);
\draw[very thick] (B6)--(B9);
\draw[very thick] (B7)--(B8);
\draw[very thick] (B7)--(B14);
\draw[very thick] (B8)--(B14);
\draw[very thick] (B8)--(B16);
\draw[very thick] (B9)--(B12);
\draw[very thick] (B9)--(B15);
\draw[very thick] (B12)--(B15);
\draw[very thick] (B14)--(B15);
\draw[very thick] (B14)--(B16);
\draw[very thick] (B15)--(B16);

\end{tikzpicture}
\par\smallskip
{\small (b) $4(555)+24(556)$\par}
\end{minipage}
\end{center}
\caption{Every global volume maximizer in $\mathcal{P}_{16}$ is combinatorially equivalent to one of the two polytopes shown. Both have degree sequence $(6^4,5^{12})$; their facial degree multisets are given below the corresponding drawings.}
\label{fig:N16-maximizing-types}
\end{figure}
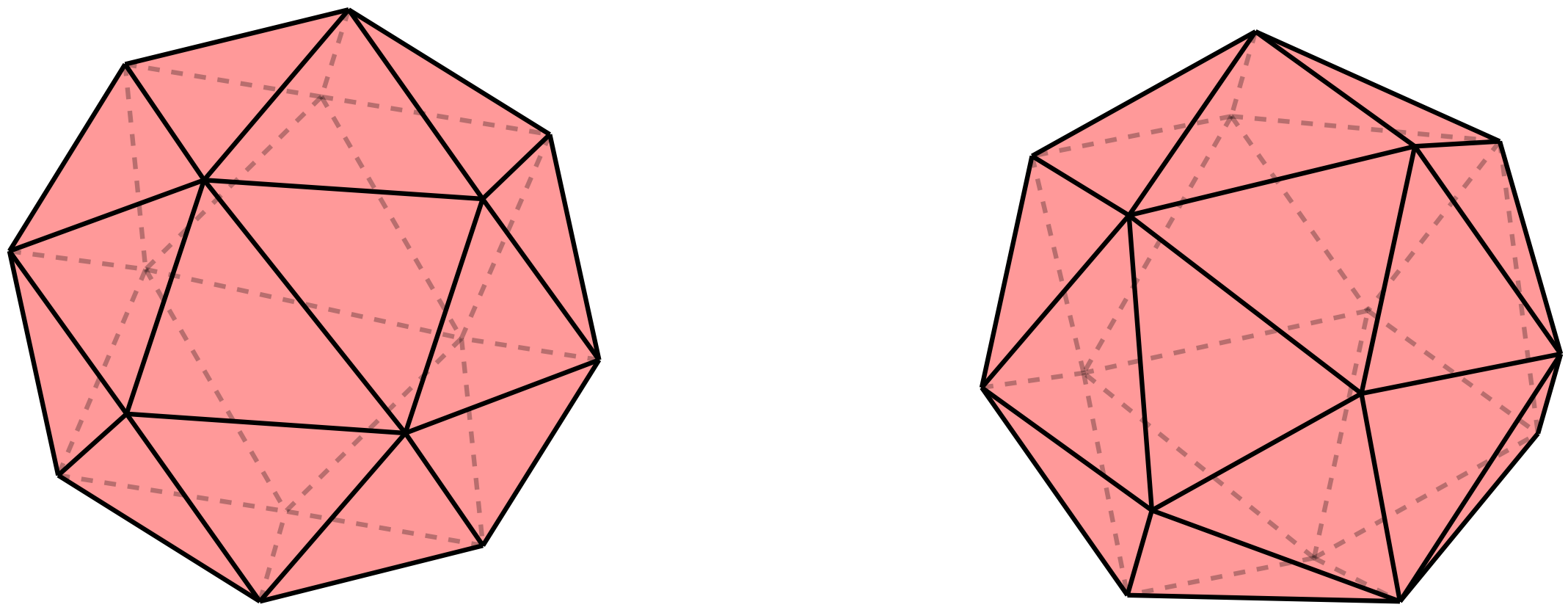

\vspace{2mm}
The maximizing combinatorial types for $N=13$ and $N=15$ are illustrated in Figure~\ref{fig:N13-maximizing-type} and Figure~\ref{fig:N15-maximizing-type}, respectively. The two possible maximizing combinatorial types for $N=16$ are depicted in Figure~\ref{fig:N16-maximizing-types}.

 \vspace{2mm}

In our final main result, we determine the global volume maximizer for polytopes with $N=32$ vertices.

\begin{theorem}\label{mainThm-N32}
Let $P\in\mathcal{P}_{32}$. Then
\begin{equation}\label{eq:sharp-N32-volume}
    \vol(P)\leq V_{32}:=\frac{2}{3}\sqrt{10(5-\sqrt{5})}.
\end{equation}
Equality holds if and only if, up to an orthogonal transformation, $P$ is an inscribed realization of the pentakis dodecahedron, which is the convex hull of the twelve vertices of a regular icosahedron inscribed in $\Sp$ and the radial projections onto $\Sp$ of its twenty facet centers.
\end{theorem}

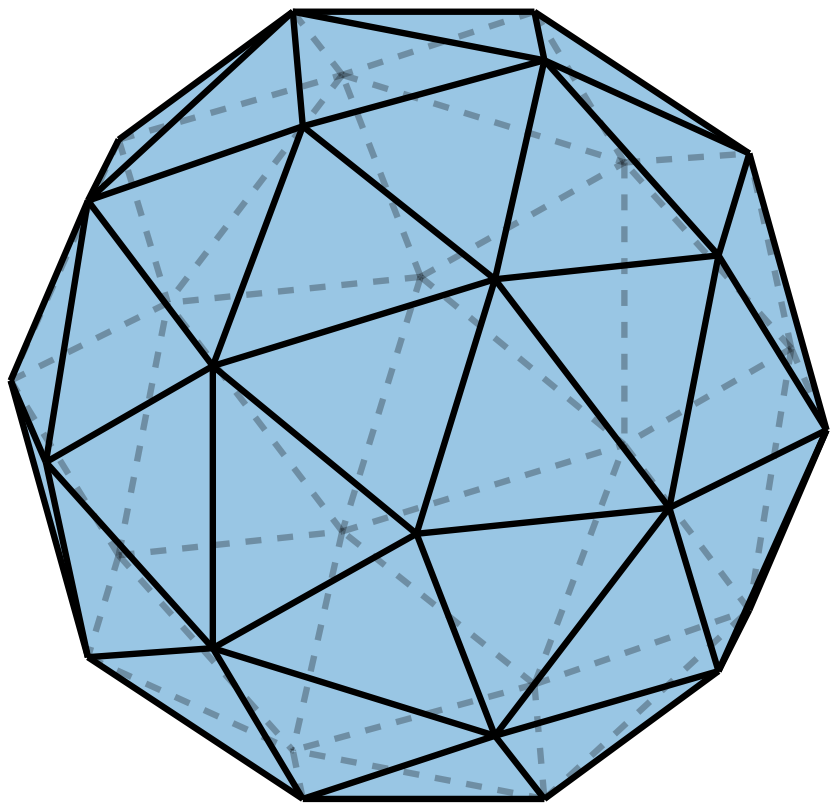
\begin{figure}[ht]
\begin{center}
\tdplotsetmaincoords{72}{122}

\begin{tikzpicture}[
    scale=2.75,
    line join=bevel,
    tdplot_main_coords
]

\pgfmathsetmacro{\a}{sqrt((5-sqrt(5))/10)}
\pgfmathsetmacro{\b}{sqrt((5+sqrt(5))/10)}
\pgfmathsetmacro{\c}{1/sqrt(3)}
\pgfmathsetmacro{\d}{sqrt((3+sqrt(5))/6)}
\pgfmathsetmacro{\e}{sqrt((3-sqrt(5))/6)}

\coordinate (A1) at (0,-\b,-\a);
\coordinate (A2) at (-\a,0,-\b);
\coordinate (A3) at (-\b,-\a,0);
\coordinate (A4) at (0,\b,-\a);
\coordinate (A5) at (-\a,0,\b);
\coordinate (A6) at (\b,-\a,0);
\coordinate (A7) at (0,-\b,\a);
\coordinate (A8) at (\a,0,-\b);
\coordinate (A9) at (-\b,\a,0);
\coordinate (A10) at (0,\b,\a);
\coordinate (A11) at (\a,0,\b);
\coordinate (A12) at (\b,\a,0);

\coordinate (B1) at (-\c,-\c,-\c);
\coordinate (B2) at (-\c,-\c,\c);
\coordinate (B3) at (\e,-\d,0);
\coordinate (B4) at (-\e,-\d,0);
\coordinate (B5) at (0,\e,-\d);
\coordinate (B6) at (\c,-\c,-\c);
\coordinate (B7) at (0,-\e,-\d);
\coordinate (B8) at (-\d,0,\e);
\coordinate (B9) at (-\d,0,-\e);
\coordinate (B10) at (-\c,\c,-\c);
\coordinate (B11) at (-\e,\d,0);
\coordinate (B12) at (-\c,\c,\c);
\coordinate (B13) at (0,\e,\d);
\coordinate (B14) at (0,-\e,\d);
\coordinate (B15) at (\c,-\c,\c);
\coordinate (B16) at (\d,0,-\e);
\coordinate (B17) at (\e,\d,0);
\coordinate (B18) at (\c,\c,-\c);
\coordinate (B19) at (\c,\c,\c);
\coordinate (B20) at (\d,0,\e);

\draw[fill opacity=0.4,fill=cyan!50!blue]
    (A4)--(B10)--(B11)--cycle;

\draw[fill opacity=0.4,fill=cyan!50!blue]
    (A7)--(B15)--(B14)--cycle;

\draw[fill opacity=0.4,fill=cyan!50!blue]
    (A6)--(B3)--(B6)--cycle;

\draw[fill opacity=0.4,fill=cyan!50!blue]
    (B13)--(B12)--(A5)--cycle;

\draw[fill opacity=0.4,fill=cyan!50!blue]
    (B14)--(B13)--(A5)--cycle;

\draw[fill opacity=0.4,fill=cyan!50!blue]
    (A6)--(B15)--(B3)--cycle;

\draw[fill opacity=0.4,fill=cyan!50!blue]
    (A8)--(B5)--(B18)--cycle;

\draw[fill opacity=0.4,fill=cyan!50!blue]
    (A4)--(B18)--(B5)--cycle;

\draw[fill opacity=0.4,fill=cyan!50!blue]
    (B11)--(B12)--(A10)--cycle;

\draw[fill opacity=0.4,fill=cyan!50!blue]
    (A8)--(B16)--(B6)--cycle;

\draw[fill opacity=0.4,fill=cyan!50!blue]
    (A6)--(B6)--(B16)--cycle;

\draw[fill opacity=0.4,fill=cyan!50!blue]
    (B13)--(A10)--(B12)--cycle;

\draw[fill opacity=0.4,fill=cyan!50!blue]
    (B14)--(B15)--(A11)--cycle;

\draw[fill opacity=0.4,fill=cyan!50!blue]
    (A4)--(B11)--(B17)--cycle;

\draw[fill opacity=0.4,fill=cyan!50!blue]
    (B14)--(A11)--(B13)--cycle;

\draw[fill opacity=0.4,fill=cyan!50!blue]
    (A8)--(B18)--(B16)--cycle;

\draw[fill opacity=0.4,fill=cyan!50!blue]
    (B11)--(A10)--(B17)--cycle;

\draw[fill opacity=0.4,fill=cyan!50!blue]
    (A4)--(B17)--(B18)--cycle;

\draw[fill opacity=0.4,fill=cyan!50!blue]
    (A6)--(B20)--(B15)--cycle;

\draw[fill opacity=0.4,fill=cyan!50!blue]
    (B20)--(A11)--(B15)--cycle;

\draw[fill opacity=0.4,fill=cyan!50!blue]
    (A6)--(B16)--(B20)--cycle;

\draw[fill opacity=0.4,fill=cyan!50!blue]
    (B13)--(B19)--(A10)--cycle;

\draw[fill opacity=0.4,fill=cyan!50!blue]
    (B13)--(A11)--(B19)--cycle;

\draw[fill opacity=0.4,fill=cyan!50!blue]
    (A12)--(B16)--(B18)--cycle;

\draw[fill opacity=0.4,fill=cyan!50!blue]
    (B19)--(B17)--(A10)--cycle;

\draw[fill opacity=0.4,fill=cyan!50!blue]
    (A12)--(B18)--(B17)--cycle;

\draw[fill opacity=0.4,fill=cyan!50!blue]
    (B20)--(B16)--(A12)--cycle;

\draw[fill opacity=0.4,fill=cyan!50!blue]
    (B20)--(B19)--(A11)--cycle;

\draw[fill opacity=0.4,fill=cyan!50!blue]
    (A12)--(B17)--(B19)--cycle;

\draw[fill opacity=0.4,fill=cyan!50!blue]
    (B20)--(A12)--(B19)--cycle;

\draw[very thick,dashed,opacity=0.28] (A1)--(B1);
\draw[very thick,dashed,opacity=0.28] (A1)--(B3);
\draw[very thick,dashed,opacity=0.28] (A1)--(B4);
\draw[very thick,dashed,opacity=0.28] (A1)--(B6);
\draw[very thick,dashed,opacity=0.28] (A1)--(B7);
\draw[very thick,dashed,opacity=0.28] (A2)--(B1);
\draw[very thick,dashed,opacity=0.28] (A2)--(B5);
\draw[very thick,dashed,opacity=0.28] (A2)--(B7);
\draw[very thick,dashed,opacity=0.28] (A2)--(B9);
\draw[very thick,dashed,opacity=0.28] (A2)--(B10);
\draw[very thick,dashed,opacity=0.28] (A3)--(B1);
\draw[very thick,dashed,opacity=0.28] (A3)--(B2);
\draw[very thick,dashed,opacity=0.28] (A3)--(B4);
\draw[very thick,dashed,opacity=0.28] (A3)--(B8);
\draw[very thick,dashed,opacity=0.28] (A3)--(B9);
\draw[very thick,dashed,opacity=0.28] (A5)--(B2);
\draw[very thick,dashed,opacity=0.28] (A5)--(B8);
\draw[very thick,dashed,opacity=0.28] (A7)--(B2);
\draw[very thick,dashed,opacity=0.28] (A7)--(B3);
\draw[very thick,dashed,opacity=0.28] (A7)--(B4);
\draw[very thick,dashed,opacity=0.28] (A8)--(B7);
\draw[very thick,dashed,opacity=0.28] (A9)--(B8);
\draw[very thick,dashed,opacity=0.28] (A9)--(B9);
\draw[very thick,dashed,opacity=0.28] (A9)--(B10);
\draw[very thick,dashed,opacity=0.28] (A9)--(B11);
\draw[very thick,dashed,opacity=0.28] (A9)--(B12);
\draw[very thick,dashed,opacity=0.28] (B1)--(B4);
\draw[very thick,dashed,opacity=0.28] (B1)--(B7);
\draw[very thick,dashed,opacity=0.28] (B1)--(B9);
\draw[very thick,dashed,opacity=0.28] (B2)--(B4);
\draw[very thick,dashed,opacity=0.28] (B2)--(B8);
\draw[very thick,dashed,opacity=0.28] (B2)--(B14);
\draw[very thick,dashed,opacity=0.28] (B3)--(B4);
\draw[very thick,dashed,opacity=0.28] (B5)--(B7);
\draw[very thick,dashed,opacity=0.28] (B5)--(B10);
\draw[very thick,dashed,opacity=0.28] (B6)--(B7);
\draw[very thick,dashed,opacity=0.28] (B8)--(B9);
\draw[very thick,dashed,opacity=0.28] (B8)--(B12);
\draw[very thick,dashed,opacity=0.28] (B9)--(B10);

\draw[very thick] (A4)--(B5);
\draw[very thick] (A4)--(B10);
\draw[very thick] (A4)--(B11);
\draw[very thick] (A4)--(B17);
\draw[very thick] (A4)--(B18);
\draw[very thick] (A5)--(B12);
\draw[very thick] (A5)--(B13);
\draw[very thick] (A5)--(B14);
\draw[very thick] (A6)--(B3);
\draw[very thick] (A6)--(B6);
\draw[very thick] (A6)--(B15);
\draw[very thick] (A6)--(B16);
\draw[very thick] (A6)--(B20);
\draw[very thick] (A7)--(B14);
\draw[very thick] (A7)--(B15);
\draw[very thick] (A8)--(B5);
\draw[very thick] (A8)--(B6);
\draw[very thick] (A8)--(B16);
\draw[very thick] (A8)--(B18);
\draw[very thick] (A10)--(B11);
\draw[very thick] (A10)--(B12);
\draw[very thick] (A10)--(B13);
\draw[very thick] (A10)--(B17);
\draw[very thick] (A10)--(B19);
\draw[very thick] (A11)--(B13);
\draw[very thick] (A11)--(B14);
\draw[very thick] (A11)--(B15);
\draw[very thick] (A11)--(B19);
\draw[very thick] (A11)--(B20);
\draw[very thick] (A12)--(B16);
\draw[very thick] (A12)--(B17);
\draw[very thick] (A12)--(B18);
\draw[very thick] (A12)--(B19);
\draw[very thick] (A12)--(B20);
\draw[very thick] (B3)--(B6);
\draw[very thick] (B3)--(B15);
\draw[very thick] (B5)--(B18);
\draw[very thick] (B6)--(B16);
\draw[very thick] (B10)--(B11);
\draw[very thick] (B11)--(B12);
\draw[very thick] (B11)--(B17);
\draw[very thick] (B12)--(B13);
\draw[very thick] (B13)--(B14);
\draw[very thick] (B13)--(B19);
\draw[very thick] (B14)--(B15);
\draw[very thick] (B15)--(B20);
\draw[very thick] (B16)--(B18);
\draw[very thick] (B16)--(B20);
\draw[very thick] (B17)--(B18);
\draw[very thick] (B17)--(B19);
\draw[very thick] (B19)--(B20);

\end{tikzpicture}
\end{center}
\caption{The global volume maximizer in $\mathcal{P}_{32}$ is the convex hull of an icosahedron and its normalized facet centers.}
\label{fig:N32-maximizer}
\end{figure}

The configuration in Theorem \ref{mainThm-N32}, its volume, and its optimality within its own combinatorial type are due to Wimmer \cite[Chapter IV, p.~55]{Wimmer1998}. The additional argument we provide in this paper is an inequality which holds for every facial degree pattern, which allows a comparison with all other combinatorial types. In particular, the proof of Theorem~\ref{mainThm-N32} does not require an enumeration of graphs with 32 vertices. See Figure~\ref{fig:N32-maximizer} for an illustration of the volume maximizer for $N=32$.

\subsection{Overview of the methods}

The degree sequence bound follows by applying Wimmer's vertex-star estimate \cite[Chapter I, Proposition 1, pp.~9--10]{Wimmer1998}, stated as Lemma \ref{lem:vertex-star-estimate} below, and then averaging within each degree class. In \cite{Wimmer1998}, Wimmer used the degree-$5,6$ specialization of this result. We formulate the estimate for arbitrary degree sequences, reduce its optimization to one variable, and then prove that the balanced degree-$5,6$ envelope is in fact an upper bound for every degree sequence when $N\geq 12$.

The ``tangent plane method'' we introduce is, to the best of the author's knowledge, new. Using the three spherical angles of each facet, we show that facial-cone volume is a strictly concave function of the spherical half-angles. A supporting plane (over the feasible set of half-angles) is then used to provide an upper bound whose linear terms can be summed using the fact that the half-angles around each vertex add up to $\pi$. We use coefficients depending on degree, and, in the last fifteen-vertex exclusion, coefficients depending on the individual vertex. For the case $N=32$, a correction involving reciprocal degrees makes the resulting bound independent of the underlying graph.

The common lemmas are proved in Section \ref{sec:lemmas}. Sections \ref{sec:mainThm-proof} and \ref{sec:N11} treat ten and eleven vertices, respectively. Section \ref{sec:larger-cases} proves the combinatorial reductions for larger values of $N$, and in Section \ref{sec:N32} we prove Theorem \ref{mainThm-N32}. The proof of the aforementioned facial concavity and the computational verification details are included at the end of the paper.

     \section{Background and notation}\label{background}

 We will work in the 3-dimensional Euclidean space $\R^3=\{x=(x_1,x_2,x_3): x_1,x_2,x_3\in\R\}$. For $x,y\in\R^3$, the standard inner product is $\langle x,y\rangle=\sum_{i=1}^3 x_i y_i$ and the norm of $x$ is $\|x\|=\sqrt{\langle x,x\rangle}$. The unit sphere in $\R^3$ centered at the origin $o\in\R^3$ is $\Sp=\{x\in\R^3: \|x\|=1\}$.

 The interior of a set $A$ is denoted by $\operatorname{int}(A)$, and the convex hull of $A$ is $\conv(A)$. In particular, if $A=\{y_1,\ldots,y_n\}\subset\R^3$ is a finite point set, then $\conv(A)$ is called a \emph{polytope}; if it is three-dimensional, it is called a \emph{3-polytope}; we also write $[y_1,\ldots,y_n]$ for shorthand. A \emph{face} of a convex polytope $P\subset\R^3$ is the intersection of $P$ with a support plane $H$ of $P$. The 0-dimensional faces of a polytope $P\subset\R^3$ are called \emph{vertices}, the 1-dimensional faces of $P$ are called \emph{edges}, and the 2-dimensional faces of $P$ are called \emph{facets}. We let $\operatorname{Vert}(P)$ and $\mathcal{F}_2(P)$ denote the vertex set and facet set of $P$, respectively. The set of all convex 3-polytopes with at most $N$ vertices and which are inscribed in $\Sp$  is denoted by $\mathcal{P}_N$. Also, a polytope $P\subset\R^3$ is \emph{simplicial} if all of its facets are triangles. For a polytope $P$ and a vertex $v\in\operatorname{Vert}(P)$, the \emph{star of $v$} is the set of all facets of $P$ that contain $v$. Two convex 3-polytopes $P,Q\subset\R^3$ are \emph{combinatorially equivalent} if there exists a bijection $\varphi$ between the set of faces of $P$ and the set of faces of $Q$ such that for any two faces $F_1$ and $F_2$ of $P$, we have $F_1\subset F_2$ if and only if $\varphi(F_1)\subset\varphi(F_2)$.

 For more background on convex polyhedra and discrete geometry, we refer the reader to, e.g., \cite{Brondsted,Grunbaum} and \cite{VDG,FejesToth-FejesToth-Kuperberg-2023}, respectively. For related results on the volume maximization problem, see, e.g.,  \cite{BermanHanes1970,VDG,Toth,Hoehner-Ledford-2022,Hoehner-Ledford-2026,HorvathIcosahedron,HorvathLangi,Wimmer1998}.

\section{Lemmas}\label{sec:lemmas}

We begin with some standard reductions and geometric estimates that will be used throughout the paper.

\subsection{Existence of volume-maximizing polytopes}

\begin{lemma}\label{lem:existence-exact-number}
For every fixed $N\geq 4$, the maximum $\max\{\vol(P): P\in\mathcal{P}_N\}$ is attained. Moreover, every maximizing polytope has exactly $N$ vertices.
\end{lemma}

\begin{proof}
For a multiset $\{p_1,\ldots,p_N\}\subset\Sp$, define the function $V:(\Sp)^N\to[0,\infty)$ by 
\[
V(p_1,\ldots,p_N):=\vol([p_1,\ldots,p_N]).
\]
The function $V$ is continuous on the compact set $(\Sp)^N$, and hence it attains its maximum. This proves the existence.

Next, let $P\in\mathcal{P}_N$ be a volume maximizer, and suppose by way of contradiction that $P$ has
$m<N$ vertices. Note that since $\mathcal{P}_N$ contains full-dimensional polytopes (e.g., the regular inscribed tetrahedron), every global volume maximizer in $\mathcal{P}_N$ must be full-dimensional.  Since $P$ is a proper subset of the unit ball,
there exists $q\in\Sp\setminus P$. Thus
$P\subsetneq \conv(P\cup\{q\})$, and since $P$ is full-dimensional, we have $\vol(\conv(P\cup\{q\}))>\vol(P)$. The new polytope $\conv(P\cup\{q\})$ has at most $m+1\leq N$ vertices lying in $\Sp$, so $\conv(P\cup\{q\})\in\mathcal{P}_N$, which contradicts the maximality of $P$. Thus, every maximizer has exactly $N$ vertices.
\end{proof}

\subsection{Property Z and volume maximizers}

The next definition is a local optimality criterion due to Berman and Hanes \cite{BermanHanes1970}.

\begin{definition}
    Let $P=\conv\{p_1,\ldots,p_N\}\in\mathcal{P}_N$. We say that $P$ satisfies \emph{Property Z} if for each $p_i$ there exists an open set $U_i\subset\Sp$ with $p_i\in U_i$ such that
    \[
\vol\big(\conv\{p_1,\ldots,p_{i-1},q,p_{i+1},\ldots,p_N\}\big)\leq \vol(P)
    \]
    for all $q\in U_i$.
\end{definition}
Since this is a condition that must be satisfied by a local maximizer, it must also be satisfied by any global maximizer.

\subsection{A volume maximizer must contain the origin in its interior}

\begin{lemma}\label{lem:origin-interior}
Let $P\subset\R^3$ be a full-dimensional simplicial polytope whose
vertices $p_1,\ldots,p_m$ lie in $\Sp$. If $P$ satisfies
Property~Z, then $o\in\operatorname{int}(P)$. Consequently, every maximum-volume polytope in $\mathcal{P}_N$
contains the origin in its interior.
\end{lemma}

\begin{proof}
For each $i\in\{1,\ldots,m\}$, set
\[
P_i:=\conv(\{p_1,\ldots,p_{i-1},p_{i+1},\ldots,p_m\}).
\]
Since $P$ is simplicial, every facet is a nondegenerate triangle, and every vertex not in that facet lies strictly on the interior side of its supporting plane. These strict inequalities persist when $p_i$ is replaced by a point $q$ sufficiently close to $p_i$ in $\R^3$. The resulting polytope therefore has the same oriented triangular facets. Indeed, the retained facets contain all vertices, and Euler's formula bounds their number by $2m-4$, which is already the number of facets of $P$.

By the oriented volume formula in Lemma \ref{lem:simp-vol}, its volume depends affinely on $q$: every determinant involving $q$ contains it in exactly one column, while all other determinants are constant. This argument also applies when $m=4$, in which case $P_i$ is only two-dimensional. Thus there exist $c_i\in\R$ and
$g_i\in\R^3$ such that
\[
\vol\big(\conv(P_i\cup\{q\})\big)
=c_i+\langle g_i,q\rangle
\]
for every $q$ sufficiently close to $p_i$. 

We first show that $g_i\neq o$. Since $p_i$ is a vertex of $P$, we have $p_i\notin P_i$. Choose
any point $c\in P_i$, and for $t\in[0,1]$, set $q_t:=(1-t)p_i+tc$. For all sufficiently small $t>0$, the point $q_t$ lies in the neighborhood in which the preceding affine volume formula is
valid. Moreover,
\[
P_i\subset\conv(P_i\cup\{q_t\})
\subsetneq\conv(P_i\cup\{p_i\})=P.
\]
Since $P$ is full-dimensional, proper containment implies strict
inequality of volumes. Hence
\[
\vol\big(\conv(P_i\cup\{q_t\})\big)<\vol(P).
\]
Using the affine volume formula, we obtain
\[
\vol\big(\conv(P_i\cup\{q_t\})\big)
=c_i+\langle g_i,q_t\rangle
=\vol(P)+t\langle g_i,c-p_i\rangle.
\]
It follows that $\langle g_i,c-p_i\rangle<0$, and therefore $g_i\neq o$.

We now apply Property Z. Let $w\in p_i^\perp\cap\Sp$,
and consider the great circle curve $q(s):=(\cos s)p_i+(\sin s)w$. 
For all sufficiently small $s$, Property Z and the affine volume
formula together imply
\[
\langle g_i,q(s)\rangle\leq\langle g_i,p_i\rangle.
\]
Thus, the function $f(s):=\langle g_i,q(s)\rangle$ has a local maximum at $s=0$. Consequently, we get $0=f'(0)=\langle g_i,w\rangle$. Since this holds for every $w\in p_i^\perp$, the vector $g_i$ is
parallel to $p_i$. Therefore, $g_i=\lambda_i p_i$ for some $\lambda_i\in\R$. Furthermore,
\[
0\geq f''(0)=-\langle g_i,p_i\rangle=-\lambda_i,
\]
so $\lambda_i\geq 0$. Since $g_i\neq o$ and $\|p_i\|=1$, we must have
$\lambda_i\neq 0$, so $\lambda_i>0$.

Now we use the translation invariance of volume. Fix an arbitrary
vector $a\in\R^3$ and translate all vertices simultaneously to the new points $p_i(t):=p_i+ta$, $i\in\{1,\ldots,m\}$. Since translation does not change volume, we have
\[
\vol\big(\conv(\{p_1(t),\ldots,p_m(t)\})\big)
=\vol(P)
\]
for all sufficiently small $t$. Differentiating at $t=0$, we get
\[
0=\sum_{i=1}^m\langle g_i,a\rangle
=\left\langle\sum_{i=1}^m g_i,a\right\rangle.
\]
Since $a\in\R^3$ was arbitrary, it follows that $\sum_{i=1}^m g_i=o$. Therefore, $\sum_{i=1}^m\lambda_i p_i=o$, where all $\lambda_i>0$. Dividing by $\Lambda:=\sum_{i=1}^m\lambda_i$, we obtain
\[
o=\sum_{i=1}^m\frac{\lambda_i}{\Lambda}p_i,
\quad
\frac{\lambda_i}{\Lambda}>0,
\quad\text{and}\quad
\sum_{i=1}^m\frac{\lambda_i}{\Lambda}=1.
\]
Thus, the origin is a convex combination of all the vertices with
strictly positive coefficients.

We claim that this implies $o\in\operatorname{int}(P)$. If not, then
since $P$ is full-dimensional and $o\in P$, there would exist a
nonzero vector $u\in\R^3$ defining a supporting hyperplane
of $P$ at the origin such that $\langle u,p_i\rangle\geq 0$ for every $i$. Taking the inner product of $\sum_i\lambda_i p_i=o$ with $u$, we obtain
\[
0=\sum_{i=1}^m\lambda_i\langle u,p_i\rangle.
\]
Every summand is nonnegative and every $\lambda_i$ is positive, so $\langle u,p_i\rangle=0$ for every $i$. Hence all the vertices of $P$ lie in the plane $u^\perp$, a contradiction to the assumption that $P$ is full-dimensional. Therefore, $o\in\operatorname{int}(P)$.

Finally, every global volume maximizer is full-dimensional, satisfies Property~Z, and is simplicial by Lemma~\ref{lem:simplicial}. The result
follows.
\end{proof}

Thus, we may assume that every global volume maximizer considered below contains the origin in its interior, once simpliciality is applied.

\subsection{Volume maximizers must be simplicial}

The following lemma is due to Berman and Hanes \cite{BermanHanes1970}.
\begin{lemma}\label{lem:simplicial}
    Let $P\in\mathcal{P}_N$ be a maximum-volume polytope. Then $P$ is simplicial.
\end{lemma}

\begin{remark}
    Lemma \ref{lem:simplicial} was extended to general dimensions $d\geq 2$ by Horv\'ath and L\'angi \cite{HorvathLangi}.
\end{remark}

\subsection{Volume decomposition of oriented simplicial polytopes}

We will also need the following volume decomposition formula for an oriented simplicial polytope, which is well-known. 

\begin{lemma}\label{lem:simp-vol}
    Let $P\subset\R^3$ be a convex simplicial 3-polytope. If every facet $F=[a,b,c]$ of $P$ is given by its outward orientation, then
    \begin{equation}\label{eq:simp-vol}
        6\vol(P)=\sum_{[a,b,c]\in\mathcal{F}_2(P)}\det(a,b,c).
    \end{equation}
\end{lemma}

\begin{remark}
    This lemma holds even if $o\not\in P$. A proof can be found in, e.g., \cite{Hoehner-Ledford-2026}.
\end{remark}

\subsection{The star volume bound}

For a facet $F$ of a polytope containing $o$ in its interior, write $R(F)$ for its radial projection onto $\Sp$. The following result is Wimmer's vertex-star estimate \cite[Chapter I, Proposition 1, pp.~9--10]{Wimmer1998}. 

\begin{lemma}\label{lem:vertex-star-estimate}
Let $P$ be a convex simplicial $3$-polytope inscribed in $\Sp$, with $o\in\operatorname{int}(P)$. Let $v$ be a vertex of degree $m$, and let $F_1,\ldots,F_m$ be the triangular facets incident with $v$, ordered cyclically around $v$. Set
\[
\alpha_v:=\sum_{i=1}^m \area(R(F_i)).
\]
Then
\[
\sum_{i=1}^m \vol([o,F_i])
\leq\Psi_m(\alpha_v),
\]
where
\begin{equation}\label{eq:Psi-original-form}
\Psi_m(\alpha)
:=\frac{m}{3}(\cot\frac{\pi}{m})
\frac{\sin^2(\pi/m)-\sin^2\left(\frac{2\pi-\alpha}{2m}\right)}{\cos^2\left(\frac{2\pi-\alpha}{2m}\right)}.
\end{equation}
\end{lemma}

We include a proof of this lemma for the reader's convenience.

\begin{proof}
After a rotation, without loss of generality we may assume that $v=e_3$. Let $p_1,\ldots,p_m$ be the neighbors of $v$, ordered cyclically so that $F_i=[v,p_i,p_{i+1}]$, $p_{m+1}=p_1$, and $\det(v,p_i,p_{i+1})>0$. Set $A_i:=\operatorname{area}(R(F_i))$ and $V_i:=\operatorname{vol}([o,F_i])$, and let $\theta_i$ be the angle of the spherical triangle $R(F_i)$ at $v$. Since $o\in\operatorname{int}(P)$, the radial projections of the facets form a geodesic triangulation of $\mathbb{S}^2$. Thus $0<\theta_i<\pi$ and $\sum_{i=1}^m\theta_i=2\pi$. 
    Moreover, $R(F_i)$ is properly contained in the spherical lune of angle $\theta_i$ determined by its two sides through $v$. Since this lune has area $2\theta_i$, we have
 $0<A_i<2\theta_i$. 
    In particular,
    \[
        0<\alpha_v=\sum_{i=1}^m A_i<4\pi.
    \]

    We first estimate the volume of a single cone. Fix an index $i$, and abbreviate
    \[
        p=p_i,\qquad q=p_{i+1},\qquad
        \theta=\theta_i,\qquad A=A_i,\qquad V=V_i.
    \]
    After a rotation about the $e_3$-axis, we may assume that
    \begin{align*}
        p&=(\sin\phi,0,\cos\phi)\quad\text{and}\quad
        q=(\sin\psi\cos\theta,\sin\psi\sin\theta,\cos\psi),
    \end{align*}
    where $0<\phi,\psi<\pi$. Indeed, neither neighbor can equal $-v$, since an edge of $P$ cannot pass through $o$. Set $u:=\tan\frac{\phi}{2}$ and $w:=\tan\frac{\psi}{2}$. 
    Then $u,w>0$, and
    \[
        V=\frac{1}{6}\det(v,p,q)
        =\frac{2uw\sin\theta}{3(1+u^2)(1+w^2)}.
    \]
By the spherical cosine laws and Girard's formula, we derive that
    \[
        \cot\frac{A}{2}
        =\frac{1+\langle v,p\rangle+\langle p,q\rangle
        +\langle q,v\rangle}{\det(v,p,q)}.
    \]
    Substituting the above coordinates, we obtain
    \[
        \cot\frac{A}{2}
        =\frac{1+uw\cos\theta}{uw\sin\theta}.
    \]
    Multiplying by $uw\sin\theta\sin(A/2)$ and rearranging, we then get
    \[
        \sin\frac{A}{2}
        =uw\sin\left(\theta-\frac{A}{2}\right).
    \]
    Since $0<A<2\theta<2\pi$, both sine factors are positive. Thus
    \[
        uw=\frac{\sin(A/2)}
        {\sin(\theta-A/2)}.
    \]

    On the other hand,
    \[
        (1+u^2)(1+w^2)-(1+uw)^2=(u-w)^2\geq 0.
    \]
    Hence
    \[
        V\leq\frac{2uw\sin\theta}{3(1+uw)^2}.
    \]
Let $x:=\theta/2$ and $y:=(\theta-A)/2$. Then $0<x<\pi/2$, $-x<y<x$, and
    \[
        uw=\frac{\sin(x-y)}{\sin(x+y)}.
    \]
    It follows that
    \begin{align*}
        V&\leq
        \frac{2\sin(2x)}{3}
        \frac{\sin(x-y)\sin(x+y)}
        {[\sin(x-y)+\sin(x+y)]^2}\\
        &=\frac{2\sin(2x)}{3}\cdot
        \frac{\sin^2x-\sin^2y}{4\sin^2x\cos^2y}
        =\frac{1}{3}(\cot x)
        \frac{\sin^2x-\sin^2y}{\cos^2y}.
    \end{align*}
Now define
    \[
        g(x,y):=\frac{1}{3}\cot x
        \frac{\sin^2x-\sin^2y}{\cos^2y}
    \]
    on the convex domain
    \[
        \mathcal{D}_g
        :=\left\{(x,y)\in\mathbb{R}^2:
        0<x<\frac{\pi}{2},\ -x<y<x\right\}.
    \]
    We claim that $g$ is strictly concave on $\mathcal{D}_g$. Differentiating, we get
    \[
        g_{xx}(x,y)
        =-\frac{2(2\sin^4x+\sin^2y)\cos x}
        {3\sin^3x\cos^2y}<0
    \]
    and
    \[
        \det(D^2g(x,y))
        =\frac{8(\sin^2x-\sin^2y)^2}
        {9\cos^6y\tan^4x}>0.
    \]
    Here the second inequality follows from $|y|<x<\pi/2$. Thus the Hessian is negative definite throughout $\mathcal{D}_g$, proving the claim.

    Returning to the incident facets, set $x_i:=\theta_i/2$ and $y_i:=(\theta_i-A_i)/2$. The single-cone estimate gives us $V_i\leq g(x_i,y_i)$, and $\sum_{i=1}^m x_i=\pi$ and $\sum_{i=1}^m y_i=\frac{2\pi-\alpha_v}{2}$. 
    Therefore, by Jensen's inequality,
    \begin{align*}
        \sum_{i=1}^m\operatorname{vol}([o,F_i])
        &=\sum_{i=1}^m V_i
        \leq\sum_{i=1}^m g(x_i,y_i)\\
        &\leq mg\left(
        \frac{\pi}{m},\frac{2\pi-\alpha_v}{2m}\right)
        =\frac{m}{3}(\cot\frac{\pi}{m})
        \frac{\sin^2(\pi/m)
        -\sin^2\left(\frac{2\pi-\alpha_v}{2m}\right)}
        {\cos^2\left(\frac{2\pi-\alpha_v}{2m}\right)}
        =\Psi_m(\alpha_v).
    \end{align*}
    This proves the lemma.
\end{proof}

We will also need the following alternative form of $\Psi_m(\alpha)$.

\begin{lemma}\label{lem:Psi-alternate-form}
    For an integer $m\geq 3$, let $\Psi_m(\alpha)$ be the function defined in \eqref{eq:Psi-original-form}. Then $\Psi_m(\alpha)$ can be expressed as
    \[
\Psi_m(\alpha)=A_m-mK_m\tan^2\left(\frac{2\pi-\alpha}{2m}\right)
    \]
    where $A_m:=\frac{m}{6}\sin\frac{2\pi}{m}$ and $K_m:=\frac{\cos^3(\pi/m)}{3\sin(\pi/m)}$.
\end{lemma}

\begin{proof}
    Starting from the fraction in \eqref{eq:Psi-original-form}, we have
    \begin{align*}
        \frac{\sin^2(\pi/m)-\sin^2\left(\frac{2\pi-\alpha}{2m}\right)}{\cos^2\left(\frac{2\pi-\alpha}{2m}\right)}&=\sin^2\left(\frac{\pi}{m}\right)\sec^2\left(\frac{2\pi-\alpha}{2m}\right)-\tan^2\left(\frac{2\pi-\alpha}{2m}\right)\\
        &=\sin^2\left(\frac{\pi}{m}\right)\left[1+\tan^2\left(\frac{2\pi-\alpha}{2m}\right)\right]-\tan^2\left(\frac{2\pi-\alpha}{2m}\right)\\
        &=\sin^2\left(\frac{\pi}{m}\right)+\left[\sin^2\left(\frac{\pi}{m}\right)-1\right]\tan^2\left(\frac{2\pi-\alpha}{2m}\right)\\
        &=\sin^2\left(\frac{\pi}{m}\right)-\cos^2\left(\frac{\pi}{m}\right)\tan^2\left(\frac{2\pi-\alpha}{2m}\right).
    \end{align*}
    Thus
    \begin{align*}
        \Psi_m(\alpha)&=\frac{m}{3}\cot\left(\frac{\pi}{m}\right)\sin^2\left(\frac{\pi}{m}\right)-\frac{m}{3}\cot\left(\frac{\pi}{m}\right)\cos^2\left(\frac{\pi}{m}\right)\tan^2\left(\frac{2\pi-\alpha}{2m}\right).
    \end{align*}
    For the first term, we have
    \begin{align*}
        \frac{m}{3}\cot\left(\frac{\pi}{m}\right)\sin^2\left(\frac{\pi}{m}\right)&=\frac{m}{3}\cos\left(\frac{\pi}{m}\right)\sin\left(\frac{\pi}{m}\right)=\frac{m}{6}\sin\left(\frac{2\pi}{m}\right)=A_m,
    \end{align*}
    and for the coefficient of $\tan^2\left(\frac{2\pi-\alpha}{2m}\right)$ in the second term we have
    \[
\frac{m}{3}\cot\left(\frac{\pi}{m}\right)\cos^2\left(\frac{\pi}{m}\right)=\frac{m\cos^3(\pi/m)}{3\sin(\pi/m)}=mK_m.
    \]
The result follows.
\end{proof}

\subsection{Volume estimates via degree sequences}

We first state one of the main tools of this paper, namely, the degree sequence bound for arbitrary degree classes. It follows from Wimmer's star bound by Jensen's inequality, and includes his degree-$5,6$ estimate as a special case.

\begin{lemma}\label{lem:new-star-bound}
    Let $P$ be a convex simplicial 3-polytope inscribed in $\Sp$, and assume that $o\in\operatorname{int}(P)$. For
$m\geq 3$, let $n_m:=\#\{v\in\operatorname{Vert}(P):\deg(v)=m\}$
denote the number of vertices of $P$ of degree $m$, and let $D(P):=\{m\geq 3:n_m>0\}$. Then
    \begin{equation}\label{eq:degree-sequence-bound}
        3\vol(P) \leq \max\left\{\sum_{m\in D(P)} n_m\Psi_m(a_m): \sum_{m\in D(P)} n_m a_m=12\pi, 0\leq a_m\leq 4\pi\right\}.
    \end{equation}
    In particular, the right-hand side of \eqref{eq:degree-sequence-bound} depends only on the degree sequence
of $P$. Moreover, each function $\Psi_m$ is strictly concave on $(0,4\pi)$. Consequently,
if equality holds in \eqref{eq:degree-sequence-bound}, then all vertices
of a fixed degree have the same total (radial) star area.
\end{lemma}

\begin{proof}
    The radial facets incident with a vertex form a proper subset of the spherical tiling, so $\alpha_v\in(0,4\pi)$ for every $v\in\operatorname{Vert}(P)$.

    Since $o\in\operatorname{int}(P)$, the radial projections of the facets
of $P$ form a tiling of $\Sp$ by spherical triangles with pairwise disjoint interiors. Hence
\[
\sum_{F\in\mathcal{F}_2(P)}\area(R(F))=4\pi.
\]
By hypothesis, every facet of $P$ is a triangle, and is therefore incident with exactly three vertices. Thus, summing $\alpha_v$ over all vertices $v\in\operatorname{Vert}(P)$, every
``radial facet'' is counted three times, so
\begin{equation}\label{eq:total-star-area}
\sum_{v\in\operatorname{Vert}(P)}\alpha_v
=3\sum_{F\in\mathcal{F}_2(P)}\area(R(F))
=12\pi.
\end{equation}
\noindent Similarly, the facial cones $[o,F]:=\conv(\{o\}\cup F)$, where $F\in\mathcal F_2(P)$, form a decomposition of $P$ with pairwise disjoint interiors. Thus,
\[
\vol(P)=\sum_{F\in\mathcal F_2(P)}\vol([o,F]).
\]
Summing the volumes of these facial cones, each
facial cone is counted once for each of the three vertices of its corresponding triangular facet. Hence
\begin{equation}\label{eq:triple-volume-count}
3\vol(P)
=\sum_{v\in\operatorname{Vert}(P)}
\sum_{\substack{F\in\mathcal{F}_2(P)\\v\in F}}\vol([o,F]).
\end{equation}
\noindent By Lemma~\ref{lem:vertex-star-estimate}, for every vertex $v\in\operatorname{Vert}(P)$, we have
\[
\sum_{\substack{F\in\mathcal{F}_2(P)\\v\in F}}\vol([o,F])
\leq\Psi_{\deg(v)}(\alpha_v).
\]
Combining this with \eqref{eq:triple-volume-count}, we obtain
\begin{equation}\label{eq:sum-star-functions}
3\vol(P)\leq
\sum_{v\in\operatorname{Vert}(P)}
\Psi_{\deg(v)}(\alpha_v).
\end{equation}

Now for each $m\in D(P)$, denote the degree-$m$ vertices of $P$ by $v_{m,1},\ldots,v_{m,n_m}$,
and set $\alpha_{m,j}:=\alpha_{v_{m,j}}$ and $a_m:=\frac{1}{n_m}\sum_{j=1}^{n_m}\alpha_{m,j}$.
Then $0<a_m<4\pi$, and \eqref{eq:total-star-area} implies
\begin{equation}\label{eq:weighted-star-constraint}
\sum_{m\in D(P)}n_m a_m
=\sum_{m\in D(P)}\sum_{j=1}^{n_m}\alpha_{m,j}
=12\pi.
\end{equation}

We now verify the concavity. Set $\vartheta:=\pi/m$ and $u:=\frac{2\pi-\alpha}{2m}$. Since
\[
\frac{\sin^2\vartheta-\sin^2 u}{\cos^2 u}
=\sin^2\vartheta-\cos^2\vartheta\tan^2 u,
\]
we can express the function $\Psi_m$ in the form
\[
\Psi_m(\alpha)
=\frac{m}{3}\cot\vartheta
\left(\sin^2\vartheta-\cos^2\vartheta\tan^2 u
\right).
\]
Differentiating and using $u'=-1/(2m)$, we get
\begin{align*}
\Psi_m'(\alpha)
&=\frac{\cos^3\vartheta}{3\sin\vartheta}
\tan u\sec^2 u\\
\Psi_m''(\alpha)
&=-\frac{\cos^3\vartheta}{6m\sin\vartheta}
\sec^2u\left(1+3\tan^2u\right).
\end{align*}
Now for every $m\geq 3$ and $\alpha\in(0,4\pi)$, we have
\[
u\in\left(-\frac{\pi}{m},\frac{\pi}{m}\right)
\subset\left(-\frac{\pi}{2},\frac{\pi}{2}\right),
\]
and hence $\Psi_m''(\alpha)<0$. Therefore, the function $\Psi_m$ is strictly concave on $(0,4\pi)$.

Finally, applying Jensen's inequality separately within each degree class, we obtain
\[
\sum_{j=1}^{n_m}\Psi_m(\alpha_{m,j})
\leq n_m\Psi_m\left(
\frac{1}{n_m}\sum_{j=1}^{n_m}\alpha_{m,j}
\right)
=n_m\Psi_m(a_m).
\]
Summing over $m\in D(P)$ and using
\eqref{eq:sum-star-functions}, we derive that
\[
3\vol(P)
\leq\sum_{m\in D(P)}n_m\Psi_m(a_m).
\]
The $\#D(P)$-tuple $(a_m)_{m\in D(P)}$ satisfies
\eqref{eq:weighted-star-constraint}, and hence
\[
3\vol(P)
\leq\max\left\{
\sum_{m\in D(P)}n_m\Psi_m(a_m): \sum_{m\in D(P)}n_m a_m=12\pi, 0\leq a_m\leq 4\pi\right\}.
\]
Note that the maximum exists because the feasible set is a compact subset of
$\mathbb{R}^{\#D(P)}$ and the objective function is continuous. The
quantities appearing in this optimization problem are only the integers
$n_m$, so the bound depends only on the degree sequence. Finally, equality in Jensen's inequality for a strictly concave
function occurs only when all of its arguments are equal. Thus, if equality
holds in \eqref{eq:degree-sequence-bound}, then for every $m\in D(P)$ we must have
$\alpha_{m,1}=\cdots=\alpha_{m,n_m}=a_m$. This completes the proof of the lemma.
\end{proof}

The next result locates the maximizer of the degree sequence relaxation.

\begin{lemma}\label{lem:degree-sequence-interior-maximizer}
Let $\mathcal{D}\subset\{3,4,\ldots\}$ be a finite set, and let $n_m\in\mathbb{N}$, $m\in\mathcal{D}$, satisfy $\sum_{m\in\mathcal{D}}n_m=N>6$.
Then the optimization problem
\[
\max\left\{
\sum_{m\in\mathcal{D}}n_m\Psi_m(a_m):
\sum_{m\in\mathcal{D}}n_m a_m=12\pi,\quad 0\leq a_m\leq 4\pi
\right\}
\]
has a unique maximizer $(a_m^*)_{m\in\mathcal{D}}$. Moreover, $0<a_m^*<2\pi$ holds
for every $m\in\mathcal{D}$, and there exists a unique $\lambda^*\in(0,1/3)$ such that $\Psi_m'(a_m^*)=\lambda^*$
for every $m\in\mathcal{D}$.
\end{lemma}

\begin{proof}
Recall that
\[
\Psi_m'(\alpha)
=
K_m\tan\left(\frac{2\pi-\alpha}{2m}\right)
\sec^2\left(\frac{2\pi-\alpha}{2m}\right),
\]
where $K_m=\frac{\cos^3(\pi/m)}{3\sin(\pi/m)}$. By Lemma~\ref{lem:new-star-bound}, the function $\Psi_m$ is strictly concave on $(0,4\pi)$. In particular, $\Psi_m'$ is strictly decreasing there. Moreover,
\begin{align*}
\Psi_m'(0)
&=K_m\tan\left(\frac{\pi}{m}\right)
\sec^2\left(\frac{\pi}{m}\right)=
\frac{\cos^3(\pi/m)}{3\sin(\pi/m)}
\frac{\sin(\pi/m)}{\cos(\pi/m)}
\frac{1}{\cos^2(\pi/m)}
=\frac{1}{3},
\end{align*}
while $\Psi_m'(2\pi)=0$,
and, similarly, $\Psi_m'(4\pi)=-1/3$.
Thus $\Psi_m'$ is strictly decreasing from $1/3$ to $-1/3$ on the interval $[0,4\pi]$.

The feasible set is nonempty and compact. Indeed, the choice $a_m=12\pi/N$, where $m\in\mathcal{D}$,
is feasible since
\[
\sum_{m\in\mathcal{D}}n_ma_m
=\frac{12\pi}{N}\sum_{m\in\mathcal{D}}n_m
=12\pi.
\]
Hence a maximizer exists. Since each $\Psi_m$ is strictly concave, the objective function
\[
(a_m)_{m\in\mathcal{D}}
\longmapsto
\sum_{m\in\mathcal{D}}n_m\Psi_m(a_m)
\]
is strictly concave on the feasible set, so the maximizer is unique.

Next, we show that no coordinate of the maximizer can lie on the boundary of $[0,4\pi]$. Suppose first that $a_r^*=0$ for some $r\in\mathcal{D}$. Since
\[
\sum_{m\in\mathcal{D}}n_ma_m^*=12\pi>0,
\]
there exists $s\in\mathcal{D}$ such that $a_s^*>0$. For sufficiently small $\varepsilon>0$, replace $a_r^*$ by  $a_r^*+\frac{\varepsilon}{n_r}$
and $a_s^*$ by $a_s^*-\frac{\varepsilon}{n_s}$. Note that this perturbation preserves the constraint. The derivative at $\varepsilon=0$ of the objective function along this perturbation is $\Psi_r'(0)-\Psi_s'(a_s^*)$. Since $a_s^*>0$ and $\Psi_s'$ is strictly decreasing, from this we derive that  $\Psi_s'(a_s^*)<\Psi_s'(0)=1/3=\Psi_r'(0)$.
Thus the derivative is positive, contradicting the maximality of $(a_m^*)_{m\in\mathcal{D}}$. Hence $a_m^*>0$ for every $m\in\mathcal{D}$.

Suppose next that $a_r^*=4\pi$ for some $r\in\mathcal{D}$. Since
\[
\frac{\sum_{m\in\mathcal{D}}n_ma_m^*}
{\sum_{m\in\mathcal{D}}n_m}
=\frac{12\pi}{N}<4\pi,
\]
there exists $s\in\mathcal{D}$ such that $a_s^*<4\pi$. For sufficiently small $\varepsilon>0$, we now replace
$a_r^*$ by $a_r^*-\frac{\varepsilon}{n_r}$ and $a_s^*$ by $a_s^*+\frac{\varepsilon}{n_s}$.
Again, the constraint is preserved, and the derivative of the objective function at $\varepsilon=0$ is $-\Psi_r'(4\pi)+\Psi_s'(a_s^*)
=\frac{1}{3}+\Psi_s'(a_s^*)$.
Since $a_s^*<4\pi$ and $\Psi_s'$ is strictly decreasing, we have $\Psi_s'(a_s^*)>\Psi_s'(4\pi)=-1/3$.
Hence this derivative is also positive, again contradicting maximality. Therefore, we have $0<a_m^*<4\pi$
for every $m\in\mathcal{D}$.

Consequently, the maximizer lies in the relative interior of the feasible set, so the method of Lagrange multipliers applies. There exists $\lambda^*\in\mathbb{R}$ such that $n_m\Psi_m'(a_m^*)=\lambda^*n_m$ for every $m\in\mathcal{D}$. Since $n_m>0$, this implies that $\Psi_m'(a_m^*)=\lambda^*$ for every $m\in\mathcal{D}$. We claim that $\lambda^*>0$. Otherwise, since $\Psi_m'$ is strictly decreasing and $\Psi_m'(2\pi)=0$, we would have $a_m^*\geq 2\pi$ for every $m\in\mathcal{D}$. Hence
\[
12\pi
=\sum_{m\in\mathcal{D}}n_ma_m^*
\geq
2\pi\sum_{m\in\mathcal{D}}n_m
=2N\pi>12\pi,
\]
which is impossible. Thus $\lambda^*>0$, which implies $a_m^*<2\pi$ for every $m\in\mathcal{D}$. Finally, since $a_m^*>0$ and $\Psi_m'$ is strictly decreasing, we obtain $\lambda^*=\Psi_m'(a_m^*)<\Psi_m'(0)=1/3$.
Therefore, $0<\lambda^*<1/3$.
The uniqueness of $\lambda^*$ follows immediately from the uniqueness of the maximizer. This completes the proof.
\end{proof}

\subsection{One-variable and supporting line forms of the degree sequence bound}

The preceding optimization can be carried out with one scalar parameter.

\begin{lemma}\label{lem:degree-sequence-lambda-reduction}
Under the hypotheses of Lemma \ref{lem:degree-sequence-interior-maximizer}, for $\lambda\in(0,1/3)$ let $t_m(\lambda)$ be the unique solution in $(0,\tan(\pi/m))$ of
\begin{equation}\label{eq:tm-lambda-cubic}
    K_m(t^3+t)=\lambda,
\end{equation}
and set $\alpha_m(\lambda):=2\pi-2m\arctan t_m(\lambda)$. Then
\[
    H(\lambda):=\sum_{m\in\mathcal{D}}n_m\alpha_m(\lambda)-12\pi
\]
is strictly decreasing and has a unique zero $\lambda_*\in(0,1/3)$. The maximizing coordinates are $a_m^*=\alpha_m(\lambda_*)$. If $H(\lambda_0)>0>H(\lambda_1)$, then $\lambda_0<\lambda_*<\lambda_1$. If, in addition, $P$ satisfies the hypotheses of Lemma \ref{lem:new-star-bound} and has $n_m$ vertices of degree $m$ for every $m\in\mathcal D$, and no other vertex degrees, then
\begin{equation}\label{eq:degree-sequence-lambda-upper-bound}
    3\vol(P)<\sum_{m\in\mathcal{D}}n_m\big[A_m-mK_m t_m(\lambda_0)^2\big].
\end{equation}
\end{lemma}

\begin{proof}
By differentiating the alternative expression for $\Psi_m$, we have
\[
    \Psi_m'(\alpha)=K_m\tan\left(\frac{2\pi-\alpha}{2m}\right)
    \sec^2\left(\frac{2\pi-\alpha}{2m}\right).
\]
Thus the equation $\Psi_m'(\alpha)=\lambda$ is equivalent to \eqref{eq:tm-lambda-cubic}. The function $q_m(t):=K_m(t^3+t)$ is strictly increasing for $t\geq0$, with $q_m(0)=0$ and $q_m(\tan(\pi/m))=1/3$. Hence $t_m(\lambda)$ exists and is unique. Moreover,
\[
    t_m'(\lambda)=\frac{1}{K_m(3t_m(\lambda)^2+1)}>0
    \quad\text{and}\quad \alpha_m'(\lambda)=-\frac{2m t_m'(\lambda)}{1+t_m(\lambda)^2}<0.
\]
Consequently, $H$ is strictly decreasing. Its limits at $0$ and $1/3$ are $(2N-12)\pi>0$ and $-12\pi$, respectively. The constraint at the unique maximizer gives $H(\lambda_*)=0$.

Finally, if $\lambda_0<\lambda_*$, then $\alpha_m(\lambda_*)<\alpha_m(\lambda_0)<2\pi$. Since $\Psi_m$ is strictly increasing on $(0,2\pi)$, we get
\[
    3\vol(P)\leq\sum_m n_m\Psi_m(\alpha_m(\lambda_*))
    <\sum_m n_m\Psi_m(\alpha_m(\lambda_0))
    =\sum_m n_m[A_m-mK_m t_m(\lambda_0)^2].
\]
This completes the proof.
\end{proof}

For the finite certificates, it is convenient to use supporting lines directly. If $b_m\in(0,4\pi)$, $g_m:=\Psi_m'(b_m)\leq\lambda$, and $c_m\geq\Psi_m(b_m)-g_mb_m$, then concavity and $a\geq0$ imply $\Psi_m(a)\leq c_m+\lambda a$. Hence
\begin{equation}\label{eq:degree-supporting-line-bound}
    3\vol(P)\leq\sum_m n_m c_m+12\pi\lambda.
\end{equation}
This form requires only finitely many evaluations at specified points, rather than the exact solution of an optimization problem.

\subsection{A global degree sequence envelope}

The degree-$5,6$ expression used by Wimmer in \cite{Wimmer1998} admits the following extension to arbitrary degree sequences.

\begin{prop}\label{prop:global-degree-envelope}
Let $N\geq12$, and let $P$ be a convex simplicial $3$-polytope with $N$ vertices inscribed in $\Sp$, with $o\in\operatorname{int}(P)$. Then
\begin{equation}\label{eq:global-degree-envelope}
\begin{aligned}
    \vol(P)\leq B_N:=\frac{1}{3}\max\big\{&12\Psi_5(x)+(N-12)\Psi_6(y):\\
    &12x+(N-12)y=12\pi,\quad 0\leq x,y\leq4\pi\big\}.
\end{aligned}
\end{equation}
For $N=12$, the term involving $y$ is omitted.
\end{prop}

\begin{proof}
Set
\[
    G(\alpha,\beta):=\frac{\sin(2\alpha)}{6}
    (1-\cot^2\alpha\cot^2\beta).
\]
The function $G$ is strictly concave on the convex domain $0<\alpha<\pi/2$, $|\beta-\pi/2|<\alpha$. Indeed, $G_{\beta\beta}<0$, and, setting $S=6G$, we derive that
\[
    \det(D^2S)=\frac{32\cot^4\alpha}{\sin^6\beta}
    [1-\cos^2\alpha-\cos^2\beta]^2>0.
\]
The last inequality follows since $|\beta-\pi/2|<\alpha$. This is the isosceles cone function used in Wimmer's proof of the star estimate in \cite{Wimmer1998}; the displayed calculation includes both sides of $\beta=\pi/2$.

We extend $\Psi_m(a)$ to real $m>2$ by defining
\[
    \Psi(m,a):=mG\left(\frac{\pi}{m},\frac{\pi}{2}+\frac{a-2\pi}{2m}\right).
\]
The perspective of a concave function is concave, so $\Psi$ is jointly concave in $(m,a)$ for $m>2$, $0<a<4\pi$, and extends continuously to the corresponding closed $a$-interval. For $0<\lambda<1/3$, define
\[
    H_\lambda(m):=\max_{0\leq a\leq4\pi}\{\Psi(m,a)-\lambda a\}.
\]
By joint concavity, we see that $H_\lambda$ is concave in $m$: apply joint concavity to two maximizing pairs and use their convex combination as a feasible pair. For every integer $m\geq 3$, by concavity we obtain
\[
    H_\lambda(m)\leq (6-m)H_\lambda(5)+(m-5)H_\lambda(6).
\]
For $m=5,6$ this is an equality. Every other integer lies outside $(5,6)$, where the extrapolated chord is an upper bound for a concave function.

By Euler's formula, $\sum_m n_m=N$ and $\sum_m mn_m=6N-12$. Combining these identities with the star estimates, we obtain
\begin{align*}
    3\vol(P)&\leq\sum_v\Psi_{\deg(v)}(\alpha_v)
    \leq\sum_m n_m H_\lambda(m)+12\pi\lambda
    \leq12H_\lambda(5)+(N-12)H_\lambda(6)+12\pi\lambda.
\end{align*}
Now we apply Lemma \ref{lem:degree-sequence-interior-maximizer} to the weights $12,N-12$, omitting the zero weight when $N=12$. At the maximizing pair in \eqref{eq:global-degree-envelope}, choose the common derivative $\lambda^*$. Then the two terms $H_{\lambda^*}(5)$ and $H_{\lambda^*}(6)$ are attained at that pair, and the last displayed upper bound is exactly $3B_N$. This proves the proposition.
\end{proof}

\begin{remark}
Wimmer \cite[Chapter IV, pp.~54--57]{Wimmer1998} uses the expression in \eqref{eq:global-degree-envelope} for polytopes whose degrees are $5$ and $6$. Proposition \ref{prop:global-degree-envelope} proves that it holds for arbitrary degree sequences. Numerically, we have $B_{15}\approx 2.822999864$, $B_{16}\approx 2.898413462$, and $B_{32}\approx 3.508374679$. These bounds exceed the candidate volumes, so facial information will be needed in these cases below.
\end{remark}

\subsection{Wimmer's trivalent vertex bound}

We record another consequence of Wimmer's estimates in \cite{Wimmer1998}, which is useful before enumerating graphs.

\begin{lemma}\label{lem:Wimmer-trivalent}
Let $P$ be a convex simplicial $3$-polytope with $N$ vertices inscribed in $\Sp$, with $o\in\operatorname{int}(P)$. If $P$ has a trivalent vertex, then
\begin{equation}\label{eq:Wimmer-trivalent}
\vol(P)\leq\max_{0\leq \theta\leq4\pi}
\left\{\Psi_3(\theta)+(2N-7)U\left(\frac{4\pi-\theta}{2N-7}\right)\right\},
\end{equation}
where
\[
    U(a):=\frac{1}{12}\left(3-\cot^2\frac{a+\pi}{6}\right)\cot\frac{a+\pi}{6}.
\]
\end{lemma}

In order to prove this lemma, we will need the following elementary combinatorial result.
\begin{lemma}\label{facial-tetrahedra-lemma}
    Let $P$ be a convex, simplicial 3-polytope with $N$ vertices that contains the origin in its interior. Then $P$ has $2(N-2)$ facets and $3(N-2)$ edges. Furthermore, if $P$ has a trivalent vertex $v_0$, then $P$ has $2N-7$ facial tetrahedra not incident with $v_0$.
\end{lemma}

\begin{proof}
    Since $P$ is simplicial, by the handshaking lemma, $2f_1(P)=3f_2(P)$. Thus, by Euler's equation,
    \begin{equation*}
        2 = f_0(P)-f_1(P)+f_2(P) = N - \frac{1}{2}f_2(P)=N-\frac{1}{3}f_1(P).
    \end{equation*}
    The first claim follows.

    Next, assume that $P$ has a trivalent vertex $v_0$. By definition, there are 3 facets incident with $v_0$, and since $P$ is simplicial, by the first part we have $f_2(P)=2(N-2)$. Therefore, there are $2(N-2)-3=2N-7$ facets of $P$ that do not contain $v_0$. Since $P$ contains the origin in its interior, this implies that there are $2N-7$ facial tetrahedra not incident with $v_0$.
\end{proof}

\begin{proof}[Proof of Lemma~\ref{lem:Wimmer-trivalent}]
The three cones incident with the trivalent vertex have total volume at most $\Psi_3(\theta)$, where $\theta$ is their total radial area. There are $2N-7$ other facial cones, and their total radial area is $4\pi-\theta$.

For completeness, the equal-area triangular cone estimate follows from the symmetry and concavity of $W$, established below in Lemmas \ref{lem:spherical-cone-formula} and \ref{lem:Types8-12-main-Lemma-concavity}. If a facial cone has radial area $a$, its three spherical half-angles sum to $(a+\pi)/2$. Averaging their permutations and applying concavity, we obtain
\[
    \vol([o,F])\leq W\left(\frac{a+\pi}{6},\frac{a+\pi}{6},\frac{a+\pi}{6}\right)=U(a).
\]
Writing $q=\cot((a+\pi)/6)$ and differentiating, we get
\[
    U'(a)=\frac{q^4-1}{24},\qquad U''(a)=-\frac{q^3(1+q^2)}{36}<0
    \quad\text{for }0<a<2\pi.
\]
Every radial facet has area in this interval. Thus, by Jensen's inequality, the remaining facial cones have total volume at most $(2N-7)U((4\pi-\theta)/(2N-7))$. Adding the estimates proves \eqref{eq:Wimmer-trivalent}. Note that this is the case $k=3$ of Wimmer's $F_k$ bound in \cite[Chapter I, p.~10]{Wimmer1998}.
\end{proof}

\subsection{A spherical cone formula}

\begin{lemma}\label{lem:spherical-cone-formula}
    Let $P$ be a convex simplicial 3-polytope inscribed in $\Sp$ with $o\in\operatorname{int}(P)$. Let $F$ be a facet of $P$, and let $T=R(F)\subset\Sp$ be its radial projection. Let $\alpha,\beta,\gamma$ be the spherical interior angles of $T$, and let $\omega=\area(T)=\alpha+\beta+\gamma-\pi$ be its spherical excess. Set $t_0=\omega/2$, $t_1=\alpha-\omega/2=\frac{\pi+\alpha-\beta-\gamma}{2}$, $t_2=\beta-\omega/2=\frac{\pi-\alpha+\beta-\gamma}{2}$, and $t_3=\gamma-\omega/2=\frac{\pi-\alpha-\beta+\gamma}{2}$. Then
    \begin{equation}
        \vol([o,F]) = \frac{2}{3\sum_{i=0}^3\cot t_i}.
    \end{equation}
\end{lemma}

\begin{proof}
    First note that for a proper spherical triangle, $t_0,t_1,t_2,t_3>0$ and $\sum_{i=0}^3 t_i=\pi$. Let $p_1,p_2,p_3$ be the vertices of $F$, and let $D=|\det(p_1,p_2,p_3)|$. Then $\vol([o,F])=D/6$. Note that since $D$ is unchanged if we interchange two vertices, without loss of generality we may assume that $D=\det(p_1,p_2,p_3)>0$.  We shall prove that
    \begin{equation}\label{eq:D-spherical-geometry-identity}
        D=\frac{1-\cos^2\alpha-\cos^2\beta-\cos^2\gamma-2\cos\alpha\cos\beta\cos\gamma}{\sin\alpha\sin\beta\sin\gamma}.
    \end{equation}

    First, define the spherical side lengths $a=\arccos\langle p_2,p_3\rangle$, $b=\arccos\langle p_3,p_1\rangle$, and $c=\arccos\langle p_1,p_2\rangle$. Thus $\langle p_2,p_3\rangle=\cos a$, $\langle p_3,p_1\rangle=\cos b$, and $\langle p_1,p_2\rangle=\cos c$. Note that since the triangle is nondegenerate, we have $a,b,c\in(0,\pi)$, so the numbers $\sin a$, $\sin b$, and $\sin c$ are all positive.

    Since
    \[
\|p_2-(\cos c)p_1\|^2=\langle p_2-(\cos c)p_1,p_2-(\cos c)p_1\rangle=1-2\cos^2 c+\cos^2 c=\sin^2 c
    \]
    and $\langle p_1,p_2-(\cos c)p_1\rangle=0$, this vector is orthogonal to $p_1$, so at $p_1$ the unit tangent vector pointing along the geodesic toward $p_2$ is $u_{12}:=\frac{p_2-(\cos c)p_1}{\sin c}$. Similarly, at $p_1$ the unit tangent vector pointing along the geodesic toward $p_3$ is $u_{13}:=\frac{p_3-(\cos b)p_1}{\sin b}$. Thus, by the definition of the spherical angle, we have $\langle u_{12},u_{13}\rangle=\cos\alpha$. Moreover, since $u_{12}$ and $u_{13}$ lie in the tangent plane $p_1^\perp$, their cross product is parallel to $p_1$, and a straightforward computation yields $\|u_{12}\times u_{13}\|=\sin\alpha$. Our chosen orientation thus gives $\langle p_1,u_{12}\times u_{13}\rangle=\sin\alpha$. Expanding the inner product, we have
    \[
\sin\alpha=\langle p_1,u_{12}\times u_{13}\rangle=\frac{\langle p_1,(p_2-(\cos c)p_1)\times(p_3-(\cos b)p_1)\rangle}{\sin b\sin c}=\frac{\langle p_1,p_2\times p_3\rangle}{\sin b\sin c}=\frac{D}{\sin b\sin c}.
    \]
    In the same way, we also get $\sin\beta=\frac{D}{\sin c\sin a}$ and $\sin\gamma=\frac{D}{\sin a\sin b}$. Multiplying these three identities together, we obtain
    \begin{equation}\label{eq:sinsinsin}
        \sin\alpha\sin\beta\sin\gamma=\frac{D^3}{\sin^2 a\sin^2 b\sin^2 c}.
    \end{equation}

    Define $n_1=\frac{p_2\times p_3}{\sin a}$, $n_2=\frac{p_3\times p_1}{\sin b}$, and $n_3=\frac{p_1\times p_2}{\sin c}$. Since $\|p_i\times p_j\|=\sqrt{1-\langle p_i,p_j\rangle^2}$, $i\neq j$, each $n_i$ is a unit vector. We will compute their pairwise inner products. On one hand, using elementary properties of the inner and cross products, we have
    \[
\langle p_3\times p_1,p_1\times p_2\rangle=\langle p_3,p_1\rangle\langle p_1,p_2\rangle-\langle p_3,p_2\rangle\langle p_1,p_1\rangle=\cos b\cos c-\cos a.
    \]
    On the other hand, computing $\cos\alpha$ from the tangent vectors above we get
    \[
\cos\alpha=\frac{\langle p_2-(\cos c)p_1,p_3-(\cos b)p_1\rangle}{\sin b\sin c}=\frac{\cos a-\cos b\cos c}{\sin b\sin c}.
    \]
    Therefore, $\langle n_2,n_3\rangle=-\cos\alpha$. Similarly, $\langle n_3,n_1\rangle=-\cos\beta$ and $\langle n_1,n_2\rangle=-\cos\gamma$. Hence the Gram matrix of the vectors $n_1,n_2,n_3$ is
    \[
G=\begin{pmatrix}
    1 & -\cos\gamma & -\cos\beta\\-\cos\gamma & 1 & -\cos\alpha\\ -\cos\beta & -\cos\alpha & 1
\end{pmatrix},
    \]
    and its determinant is
    \begin{equation}\label{eq:Gram-det}
\det(G)=1-\cos^2\alpha-\cos^2\beta-\cos^2\gamma-2\cos\alpha\cos\beta\cos\gamma.
    \end{equation}

    Next we show that $D$ yields the same determinant. Define the matrices $P=\begin{pmatrix}
        p_1 & p_2 & p_3
    \end{pmatrix}$ and $C=\begin{pmatrix}
        p_2\times p_3 & p_3\times p_1 & p_1\times p_2
    \end{pmatrix}$. Then
    \[
P^\top C=\begin{pmatrix}
    \langle p_1,p_2\times p_3\rangle & \langle p_1,p_3\times p_1\rangle & \langle p_1,p_1\times p_2\rangle\\
    \langle p_2,p_2\times p_3\rangle & \langle p_2,p_3\times p_1\rangle & \langle p_2,p_1\times p_2\rangle\\
    \langle p_3,p_2\times p_3\rangle & \langle p_3,p_3\times p_1\rangle & \langle p_3,p_1\times p_2\rangle
\end{pmatrix}
=\begin{pmatrix}
    D & 0 & 0\\0 & D & 0\\ 0 & 0 & D
\end{pmatrix}=DI_3
    \]
    where $I_3$ is the $3\times 3$ identity matrix. Taking determinants, we get $\det(P^\top)\det(C)=D^3$. Since $\det(P^\top)=\det(P)=D$, this implies $\det(C)=D^2$. Now setting $N=\begin{pmatrix}
        n_1 & n_2 & n_3
    \end{pmatrix}$, we get
    \[
N=C\begin{pmatrix}
    1/\sin a & 0 & 0\\ 0 & 1/(\sin b) & 0\\ 0 & 0 & 1/(\sin c)
\end{pmatrix}.
    \]
    Therefore, $\det(N)=\frac{\det(C)}{\sin a\sin b\sin c}=\frac{D^2}{\sin a\sin b\sin c}$. By definition of the Gram matrix, $G=N^\top N$, so
    \[
\det(G)=\det(N^\top N)=\det(N^\top)\det(N)=(\det(N))^2=\frac{D^4}{\sin^2 a\sin^2 b\sin^2 c}.
    \]
    Substituting this back into \eqref{eq:sinsinsin}, we get
    \[
\frac{D^4}{\sin^2 a\sin^2 b\sin^2 c}=D\sin\alpha\sin\beta\sin\gamma.
    \]
    Thus, $\det(G)=D\sin\alpha\sin\beta\sin\gamma$. Comparing this with \eqref{eq:Gram-det}, we find that
    \[
D\sin\alpha\sin\beta\sin\gamma=1-\cos^2\alpha-\cos^2\beta-\cos^2\gamma-2\cos\alpha\cos\beta\cos\gamma.
    \]
    As mentioned earlier, $\alpha,\beta,\gamma\in(0,\pi)$, so their sines are all positive. Therefore,
    \[
D=\frac{1-\cos^2\alpha-\cos^2\beta-\cos^2\gamma-2\cos\alpha\cos\beta\cos\gamma}{\sin\alpha\sin\beta\sin\gamma},
    \]
    which proves \eqref{eq:D-spherical-geometry-identity}.

    In the next step, we show that the numerator of $D$ equals

    \begin{equation}\label{eq:trig-identity-2}
    \begin{aligned}
        1-\cos^2\alpha&-\cos^2\beta-\cos^2\gamma-2\cos\alpha\cos\beta\cos\gamma\\
        &=-4\cos\frac{\alpha+\beta+\gamma}{2}\cos\frac{\beta+\gamma-\alpha}{2}\cos\frac{\gamma+\alpha-\beta}{2}\cos\frac{\alpha+\beta-\gamma}{2}.
            \end{aligned}
        \end{equation}

    To see this, let $\theta_1=(\alpha+\beta+\gamma)/2$, $\theta_2=(\beta+\gamma-\alpha)/2$, $\theta_3=(\gamma+\alpha-\beta)/2$, and $\theta_4=(\alpha+\beta-\gamma)/2$. Then $\theta_1+\theta_2=\beta+\gamma$ and $\theta_1-\theta_2=\alpha$. Therefore, using the identity $2\cos x\cos y=\cos(x+y)+\cos(x-y)$, we obtain $2\cos\theta_1\cos\theta_2=\cos(\beta+\gamma)+\cos\alpha$. Likewise, $\theta_3+\theta_4=\alpha$ and $\theta_3-\theta_4=\gamma-\beta$, so $2\cos\theta_3\cos\theta_4=\cos\alpha+\cos(\gamma-\beta)$. Multiplying the previous two equations, and then multiplying both sides of the resulting equation by $-1$, we get
    \begin{align*}
&-4\prod_{i=1}^4\cos\theta_i=-\bigl(\cos(\beta+\gamma)+\cos\alpha\bigr)\bigl(\cos\alpha+\cos(\gamma-\beta)\bigr)\\
&=-\cos^2\alpha-\cos\alpha\cos(\beta+\gamma)-\cos\alpha\cos(\gamma-\beta)-\cos(\beta+\gamma)\cos(\gamma-\beta)\\
&=-\cos^2\alpha-\cos\alpha\bigl(\cos(\beta+\gamma)+\cos(\beta-\gamma)\bigr)-\cos(\beta+\gamma)\cos(\gamma-\beta)\\
&=-\cos^2\alpha-2\cos\alpha\cos\beta\cos\gamma-\cos(\beta+\gamma)\cos(\gamma-\beta)\\
&=-\cos^2\alpha-2\cos\alpha\cos\beta\cos\gamma-\frac{1}{2}\bigl(\cos((\beta+\gamma)+(\gamma-\beta))+\cos((\beta+\gamma)-(\gamma-\beta))\bigr)\\
&=-\cos^2\alpha-2\cos\alpha\cos\beta\cos\gamma-\frac{\cos(2\beta)+\cos(2\gamma)}{2}\\
&=-\cos^2\alpha-2\cos\alpha\cos\beta\cos\gamma-\frac{2\cos^2\beta-1+2\cos^2\gamma-1}{2}\\
&=1-\cos^2\alpha-\cos^2\beta-\cos^2\gamma-2\cos\alpha\cos\beta\cos\gamma.
    \end{align*}
    This proves \eqref{eq:trig-identity-2}.

    Now note that $\theta_1=\pi/2+t_0$, $\theta_2=\pi/2-t_1$, $\theta_3=\pi/2-t_2$, and $\theta_4=\pi/2-t_3$. Hence the numerator of $D$ becomes
    \[
-4\cos(\pi/2+t_0)\cos(\pi/2-t_1)\cos(\pi/2-t_2)\cos(\pi/2-t_3)=4\prod_{i=0}^3\sin t_i.
    \]
    Moreover, the denominator of $D$ can be expressed as
    \[
\sin\alpha\sin\beta\sin\gamma=\sin(t_0+t_1)\sin(t_0+t_2)\sin(t_0+t_3).
    \]
    To complete the proof, for $t_0,t_1,t_2,t_3\in(0,\pi)$ defined as above and which satisfy $\sum_{i=0}^3 t_i=\pi$, we will prove that
    \begin{equation}\label{eq:cotangent-sine-identity}
        \sum_{i=0}^3\cot t_i=\frac{\sin(t_0+t_1)\sin(t_0+t_2)\sin(t_0+t_3)}{\prod_{i=0}^3\sin t_i}.
    \end{equation}
   Combining denominators, we get
   \begin{equation}\label{eq:sum-of-cotangents}
\sum_{i=0}^3\cot t_i=\sum_{i=0}^3\frac{\cos t_i}{\sin t_i}=\frac{S}{\prod_{i=0}^3\sin t_i}
   \end{equation}
   where
   \begin{align*}
   S=\cos t_0\sin t_1\sin t_2\sin t_3&+\sin t_0\cos t_1\sin t_2\sin t_3\\
   &+\sin t_0\sin t_1\cos t_2\sin t_3+\sin t_0\sin t_1\sin t_2\cos t_3.
   \end{align*}
   The first two summands of $S$ can be expressed as
   \begin{align*}
       \cos t_0\sin t_1\sin t_2\sin t_3+\sin t_0\cos t_1\sin t_2\sin t_3
       &=(\cos t_0\sin t_1+\sin t_0\cos t_1)\sin t_2\sin t_3\\
       &=\sin(t_0+t_1)\sin t_2\sin t_3.
   \end{align*}
   Likewise, the latter two summands of $S$ can be expressed as
   \begin{align*}
       \sin t_0\sin t_1\cos t_2\sin t_3+\sin t_0\sin t_1\sin t_2\cos t_3
       &=\sin t_0\sin t_1(\cos t_2\sin t_3+\sin t_2\cos t_3)\\
       &=\sin t_0\sin t_1\sin(t_2+t_3).
   \end{align*}
   Thus,
   \begin{equation}\label{eq:S-sines}
S=\sin(t_0+t_1)\sin t_2\sin t_3+\sin t_0\sin t_1\sin(t_2+t_3).
   \end{equation}
   Since $t_2+t_3=\pi-(t_0+t_1)$, we have $\sin(t_2+t_3)=\sin(\pi-(t_0+t_1))=\sin(t_0+t_1)$, so
   \[
   S=\sin(t_0+t_1)(\sin t_2\sin t_3+\sin t_0\sin t_1).
   \]
   Thus, it suffices to prove that
   \begin{equation}\label{eq:sines-eqn}
       \sin t_2\sin t_3+\sin t_0\sin t_1=\sin(t_0+t_2)\sin(t_0+t_3).
   \end{equation}
   Using the identity $\sin x\sin y=\frac{1}{2}(\cos(x-y)-\cos(x+y))$, the left-hand side of \eqref{eq:sines-eqn} equals
   \begin{align*}
       \sin t_2\sin t_3+\sin t_0\sin t_1&=\frac{1}{2}(\cos(t_2-t_3)-\cos(t_2+t_3))+\frac{1}{2}(\cos(t_0-t_1)-\cos(t_0+t_1))\\
       &=\frac{1}{2}\bigl(\cos(t_2-t_3)+\cos(t_0-t_1)-\cos(t_2+t_3)-\cos(t_0+t_1)\bigr)
   \end{align*}
   But since $t_2+t_3=\pi-(t_0+t_1)$, we have $\cos(t_2+t_3)=\cos(\pi-(t_0+t_1))=-\cos(t_0+t_1)$, so the final two terms cancel and we get
\begin{equation}\label{eq:trig-identity-3}
 \sin t_2\sin t_3+\sin t_0\sin t_1=\frac{1}{2}\bigl(\cos(t_2-t_3)+\cos(t_0-t_1)\bigr).
 \end{equation}
 Meanwhile, for the right-hand side of \eqref{eq:sines-eqn}, we have
 \begin{align*}
     \sin(t_0+t_2)\sin(t_0+t_3)&=\frac{1}{2}\bigl(\cos((t_0+t_2)-(t_0+t_3))-\cos((t_0+t_2)+(t_0+t_3))\bigr)\\
     &=\frac{1}{2}\bigl(\cos(t_2-t_3)-\cos(2t_0+t_2+t_3)\bigr).
 \end{align*}
 Again using $\sum_{i=0}^3 t_i=\pi$, we obtain
 \[
2t_0+t_2+t_3=t_0+(t_0+t_2+t_3)=\pi+t_0-t_1.
 \]
 Thus $\cos(2t_0+t_2+t_3)=\cos(\pi+t_0-t_1)=-\cos(t_0-t_1)$, which implies
 \begin{align*}
     \sin(t_0+t_2)\sin(t_0+t_3)=\frac{1}{2}\bigl(\cos(t_2-t_3)+\cos(t_0-t_1)\bigr).
     \end{align*}
     Comparing this with \eqref{eq:trig-identity-3}, we obtain \eqref{eq:sines-eqn}.

     Therefore,  \eqref{eq:S-sines} may be expressed as
     \[
S=\sin(t_0+t_1)\sin(t_0+t_2)\sin(t_0+t_3).
     \]
     Substituting this back into \eqref{eq:sum-of-cotangents}, we obtain
     \[
\sum_{i=0}^3\cot t_i=\frac{\sin(t_0+t_1)\sin(t_0+t_2)\sin(t_0+t_3)}{\prod_{i=0}^3\sin t_i}.
     \]
     Combining this with
     \[
D=\frac{4\prod_{i=0}^3\sin t_i}{\sin(t_0+t_1)\sin(t_0+t_2)\sin(t_0+t_3)},
     \]
     we obtain
     \[
D=\frac{4}{\sum_{i=0}^3\cot t_i}.
     \]
     Therefore,
     \[
\vol([o,F])=\frac{D}{6}=\frac{2}{3\sum_{i=0}^3 \cot t_i}.
     \]
     This completes the proof.
\end{proof}

\subsection{Concavity of facial cone volumes}

Define the feasible set
\begin{equation}\label{eq:Types8-12-feasible-set}
\mathcal{D}:=\{(t_0,t_1,t_2,t_3)\in(0,\pi)^4: \sum_{i=0}^3 t_i=\pi\},
\end{equation}
and consider the objective function $f:\mathcal{D}\to(0,\infty)$ defined by
\begin{equation}\label{eq:Types8-12-objective-fn}
f(t_0,t_1,t_2,t_3):=\frac{1}{\sum_{i=0}^3\cot t_i}.
\end{equation}
By Lemma \ref{lem:spherical-cone-formula}, if $F$ is a spherical facet of $P$, then $\vol([o,F])=\frac{2}{3}f(t_0,t_1,t_2,t_3)$.

\begin{lemma}\label{lem:Types8-12-main-Lemma-concavity}
    The function $f:\mathcal{D}\to(0,\infty)$ defined in \eqref{eq:Types8-12-objective-fn} is strictly concave on $\mathcal{D}$.
\end{lemma}

The proof of Lemma \ref{lem:Types8-12-main-Lemma-concavity} is long and quite technical, so we postpone it until the appendix. We first derive the supporting tangent plane inequalities that will be used in the proofs.

Recall that for a spherical facet $R(F)$ of a convex simplicial 3-polytope $P$ with $o\in\operatorname{int}(P)$, we let $\alpha,\beta,\gamma$ denote its interior spherical angles. Set $\theta_1=\alpha/2$, $\theta_2=\beta/2$, and $\theta_3=\gamma/2$. Then $\vol([o,F])=W(\theta_1,\theta_2,\theta_3)=\frac{2}{3}f(t_0,t_1,t_2,t_3)$ where $t_0=\theta_1+\theta_2+\theta_3-\pi/2$, $t_1=\pi/2+\theta_1-\theta_2-\theta_3$, $t_2=\pi/2-\theta_1+\theta_2-\theta_3$, and $t_3=\pi/2-\theta_1-\theta_2+\theta_3$.

Since $f$ is strictly concave on $\mathcal{D}$, the function $W$ is also strictly concave on its domain since it is obtained from $f$ by an affine transformation of its coordinates. Note that the affine map from the half-angles to the $t_i$ is injective since $2\theta_j=t_0+t_j$. Thus, for any fixed point $s=(a,b,c)$ in the domain of $W$, we have
    \begin{equation*}
        W(\theta_1,\theta_2,\theta_3) \leq W(a,b,c)+W_{\theta_1}(a,b,c)(\theta_1-a)+W_{\theta_2}(a,b,c)(\theta_2-b)+W_{\theta_3}(a,b,c)(\theta_3-c),
    \end{equation*}
    or, equivalently,
    \begin{equation}
        W(\theta_1,\theta_2,\theta_3)\leq C+g_1\theta_1+g_2\theta_2+g_3\theta_3
    \end{equation}
where $g_1=W_{\theta_1}(a,b,c)$, $g_2=W_{\theta_2}(a,b,c)$, $g_3=W_{\theta_3}(a,b,c)$, and $C=W(a,b,c)-g_1a-g_2b-g_3c$. We have
\[
\frac{\partial(t_0,t_1,t_2,t_3)}{\partial\theta_1}=(1,1,-1,-1), \frac{\partial(t_0,t_1,t_2,t_3)}{\partial\theta_2}=(1,-1,1,-1), \frac{\partial(t_0,t_1,t_2,t_3)}{\partial\theta_3}=(1,-1,-1,1),
\]
and, writing $S:=\sum_{i=0}^3\cot t_i$, we have $\frac{\partial f}{\partial \theta_1}=-\frac{1}{S^2}\frac{\partial S}{\partial\theta_1}=f^2(\csc^2 t_0+\csc^2 t_1-\csc^2 t_2-\csc^2 t_3)$. Continuing in this fashion, we get:
\begin{equation}\label{eq:partials-of-W}
\begin{aligned}
    \frac{\partial W}{\partial\theta_1}&=\frac{2}{3}f^2(\csc^2 t_0+\csc^2 t_1-\csc^2 t_2-\csc^2 t_3)\\
    \frac{\partial W}{\partial\theta_2}&=\frac{2}{3}f^2(\csc^2 t_0-\csc^2 t_1+\csc^2 t_2-\csc^2 t_3)\\
      \frac{\partial W}{\partial\theta_3}&=\frac{2}{3}f^2(\csc^2 t_0-\csc^2 t_1-\csc^2 t_2+\csc^2 t_3).
\end{aligned}
\end{equation}
For the subsequent numerical comparisons, we will also use the rational estimate
\begin{equation}\label{eq:pi-rational-estimates}
    3.1415926<\pi<3.1415927.
\end{equation}

\subsection{Summing supporting tangent planes}

For brevity, we write $(ijk)$ for a triangular facet whose three vertices have degrees $i,j,k$. The following lemma allows coefficients either by degree or by vertex.

\begin{lemma}\label{lem:tangent-plane-bound}
Let $P$ be a convex simplicial $3$-polytope inscribed in $\Sp$ with $o\in\operatorname{int}(P)$. For each facet $F$, choose a supporting half-angle triple $s_F$ in the domain of $W$. For a corner $v\in F$, let $g_{F,v}$ be the corresponding partial derivative at $s_F$, and set
\[
    C_F:=W(s_F)-\sum_{v\in F}g_{F,v}s_{F,v}.
\]
If $\lambda_v\geq g_{F,v}$ whenever $v\in F$, and $\overline C_F\geq C_F$, then
\begin{equation}\label{eq:vertex-tangent-plane-bound}
    \vol(P)\leq\sum_{F\in\mathcal F_2(P)}\overline C_F+\pi\sum_{v\in\operatorname{Vert}(P)}\lambda_v.
\end{equation}
In particular, if the points and intercept bounds depend only on the facial degree pattern and $\lambda_v=\lambda_{\deg(v)}$, then
\begin{equation}\label{eq:degree-tangent-plane-bound}
    \vol(P)\leq\sum_{(ijk)}f_{ijk}\overline C_{(ijk)}+\pi\sum_m n_m\lambda_m,
\end{equation}
where $f_{ijk}$ is the number of facets of pattern $(ijk)$.
\end{lemma}

\begin{proof}
Let $\theta_{F,v}$ denote the half-angle at the corner $v$ of $R(F)$. By concavity,
\[
    \vol([o,F])=W(\theta_F)\leq C_F+\sum_{v\in F}g_{F,v}\theta_{F,v}
    \leq\overline C_F+\sum_{v\in F}\lambda_v\theta_{F,v}
\]
where in the last inequality we used $\theta_{F,v}>0$. The radial facets tile $\Sp$, so the spherical angles around every vertex add up to $2\pi$. Thus $\sum_{F\ni v}\theta_{F,v}=\pi$. Summing the displayed inequalities proves \eqref{eq:vertex-tangent-plane-bound}, and grouping equal degrees and equal facial patterns, we obtain \eqref{eq:degree-tangent-plane-bound}.
\end{proof}

\begin{remark}\label{rem:fixed-type-sharpness}
Suppose a realization of a fixed triangulation has the property that all partial derivatives of $W$ at corners incident with a given vertex agree. Its supporting planes then give equality in \eqref{eq:vertex-tangent-plane-bound}, so it is optimal among all inscribed realizations of that triangulation containing the origin. Strict concavity determines every facial angle in the equality case. The angles determine the spherical side lengths, and propagation across shared edges determines the realization up to an orthogonal transformation.

In particular, this applies when every face has one degree pattern $(pqr)$ and the corresponding half-angles are $(\pi/p,\pi/q,\pi/r)$. Examples are regular bipyramids, the thirty-two-vertex configuration, and the spherical flag subdivisions of a regular tetrahedron, cube, and icosahedron, with, respectively, $14,26,62$ vertices and facial patterns $(466)$, $(468)$, $(4,6,10)$. The flag vertices are the radial projections of the original vertices, edge midpoints, and face centers. These are statements about a prescribed combinatorial type.
\end{remark}

\subsection{Wimmer's ten- and eleven-vertex reductions}

We will use the following consequences of Wimmer  \cite{Wimmer1998}.

\begin{lemma}\label{lem:Wimmer-small-reductions}
Every volume maximizer in $\mathcal P_{10}$ has one of two combinatorial types. Their degree sequences are $(6,5^6,4^3)$ and $(5^8,4^2)$, and the latter type is the gyroelongated square bipyramid.

Every volume maximizer in $\mathcal P_{11}$ has one of two combinatorial types, with degree sequences $(6^2,5^6,4^3)$ and $(6,5^8,4^2)$, respectively. The second degree sequence has a unique triangulation in which all vertices have minimum degree 4, namely, the edge-contracted icosahedron.
\end{lemma}

The ten-vertex reduction can be found in \cite[Chapter I, pp.~12--14; Chapter III, p.~29; Chapter IV, p.~36]{Wimmer1998}, while the eleven-vertex reduction is in \cite[Chapter IV, p.~43]{Wimmer1998}. 

For the eleven-vertex argument, we will use a slightly weaker consequence of the second reduction: only the two displayed degree sequences need to be considered. Among all triangulations with no trivalent vertices, there are two graphs with degree sequence $(6^2,5^6,4^3)$. We will use the labels Types 32 and 33 for these graphs, and Type 34 for the graph with degree sequence $(6,5^8,4^2)$. For ten vertices, we similarly use Types 11 and 12 for the two graphs in the first part of the lemma. The complete degree sequence tables, with the Type numbering used here, are given in Appendix \ref{sec:degree-sequence-tables}; see Tables \ref{tab:N10-degree-sequences} and \ref{tab:N11-degree-sequences}.

\begin{table}[ht]
\centering
\begin{tabular}{c|c|l}
$N$ & Type & Facial degree patterns\\
\hline
10 & 11 & $6(455)+6(456)+4(555)$\\
10 & 12 & $8(455)+8(555)$\\
11 & 32 & $8(456)+4(455)+2(555)+4(556)$\\
11 & 33 & $2(446)+4(455)+4(456)+2(555)+6(556)$\\
11 & 34 & $4(455)+4(456)+8(555)+2(556)$
\end{tabular}
\caption{Facial degree patterns used in the ten- and eleven-vertex proofs.}
\label{tab:facial-degree-patterns}\label{tab:facial-degree-patterns-N11}
\end{table}

The facial data can be obtained by enumerating planar triangulations with \texttt{plantri} \cite{BrinkmannMcKay,plantri58}. The source certificates include all twelve ten-vertex and all thirty-four eleven-vertex triangulations with no trivalent vertices and check the entries of Table \ref{tab:facial-degree-patterns}. Note that the proofs below do not repeat the types already eliminated by Wimmer.

\section{Proof of Theorem \ref{mainThm}}\label{sec:mainThm-proof}

By Lemma \ref{lem:Wimmer-small-reductions}, only Types 11 and 12 need be considered. We first eliminate Type 11 by the tangent plane method, and then determine the best realization of Type 12. Throughout the argument, it suffices to consider a global volume maximizer, which, by Lemma~\ref{lem:origin-interior}, must contain the origin in its interior.

\subsection{Type 11}

Let $P\in\mathcal{P}_{10}$ have Type 11 combinatorial type. To eliminate Type 11, we will use the following supporting half-angle triples in our estimation of $W$:

\begin{table}[ht]
\centering
\begin{tabular}{c|c}
Facial degree pattern $(ijk)$ & Supporting half-angle triple $s_{(ijk)}=(s_i,s_j,s_k)$\\
\hline
$(455)$ & $(0.754932,\;0.613161,\;0.613161)$\\
$(456)$ & $(0.815873,\;0.601358,\;0.523605)$\\
$(555)$ & $(0.656954,\;0.656954,\;0.656954)$\\
\end{tabular}
\caption{Supporting half-angle triples for Type $11$.}
\label{tab:type11-supporting-triples}
\end{table}

\noindent Start with pattern $(456)$ and its corresponding triple $s_{(456)}$. We have:
\begin{align*}
    t_0&=0.815873+0.601358+0.523605-\pi/2=0.3700396732\ldots\\
    t_1&=\pi/2+0.815873-0.601358-0.523605=1.2617063268\ldots\\
    t_2&=\pi/2-0.815873+0.601358-0.523605=0.8326763268\ldots\\
    t_3&=\pi/2-0.815873-0.601358+0.523605=0.6771703268\ldots
\end{align*}
as well as $S=5.0506914864\ldots$ and $f=1/S=0.1979926912\ldots$. Hence
$W(s_{(456)})=\frac{2}{3}f=0.1319951275\ldots$. Moreover, substituting the preceding values of $t_0$, $t_1$, $t_2$, and $t_3$ into \eqref{eq:partials-of-W}, we get:
\begin{equation}
    \begin{aligned}
        W_{\theta_1}(s_{(456)})&=0.1142891545\ldots<0.11431\\
        W_{\theta_2}(s_{(456)})&=0.1522092431\ldots<0.15223\\
        W_{\theta_3}(s_{(456)})&=0.1898194553\ldots<0.18984.
    \end{aligned}
\end{equation}
Therefore, $C_{(456)}=-0.1521729698\ldots<-0.15217$.

Continuing in this fashion, we obtain the following partial derivative and tangent plane intercept data for the Type 11 patterns (rounded to ten decimal places):

\begin{table}[ht]
\centering
\begin{tabular}{c|cccc}
Pattern $(ijk)$ & $W_{\theta_1}$ & $W_{\theta_2}$ & $W_{\theta_3}$ & $C_{(ijk)}$\\
\hline
$(455)$ &
$0.1142905213$ &
$0.1522103700$ &
$0.1522103700$ &
$-0.1307741524$\\
$(456)$ &
$0.1142891545$ &
$0.1522092431$ &
$0.1898194553$ &
$-0.1521729698$\\
$(555)$ &
$0.1522097145$ &
$0.1522097145$ &
$0.1522097145$ &
$-0.1574884207$\\
\end{tabular}
\caption{Values of the partial derivatives and tangent plane intercepts for the supporting half-angle triples used for Type $11$.}
\label{tab:type11-tangent-data}
\end{table}

For Type 11, now we will take $\lambda_4=0.11431$, $\lambda_5=0.15223$, and $\lambda_6=0.18984$. Using these values and increasing the preceding slopes and intercepts slightly, we get the following bounds. Here the $\theta_m$'s denote arbitrary facial half-angles at degree-$m$ vertices, with primes distinguishing corners of the same degree:

\begin{equation}
\begin{aligned}
    W_{(455)}&\leq -0.13076+\lambda_4\theta_4+\lambda_5\theta_5+\lambda_5\theta_5'\\
    W_{(456)}&\leq -0.15217+\lambda_4\theta_4+\lambda_5\theta_5+\lambda_6\theta_6\\
    W_{(555)}&\leq -0.15748+\lambda_5\theta_5+\lambda_5\theta_5'+\lambda_5\theta_5''.
\end{aligned}
\end{equation}

Using these data and Lemma \ref{lem:tangent-plane-bound}, we obtain
\begin{align*}
    \vol(P) &\leq 4(-0.15748)+6(-0.13076)+6(-0.15217)+\pi(3\lambda_4+6\lambda_5+\lambda_6)\\
    &<2.215715<V_*.
\end{align*}
Therefore, Type 11 is eliminated.
\subsection{Type 12}

We have thus shown that the global volume maximizer in $\mathcal{P}_{10}$ must have Type 12 combinatorial type, which corresponds to the combinatorial class of the gyroelongated square bipyramid. To complete the proof of Theorem \ref{mainThm}, it  remains to determine the volume maximizer in this combinatorial class. The method we use is similar to the one employed in \cite{Hoehner-Ledford-2026} to determine the global maximizer in the final combinatorial class of the $N=9$ case.

To be clear about the goal we are marching toward, we formulate the statement in a theorem.

\begin{theorem}\label{thm:Class-12-maximizer}
    Let $P\in\mathcal{P}_{10}$ have Type 12 combinatorial type $(5^8, 4^2)$. Then
 \[
\vol(P)
\leq \frac{4}{3}(1-h_*^2)\left(1+(1+\sqrt{2})h_*\right)=:V_*
\]
where $h_*=\frac{\sqrt{10+6\sqrt{2}}-1}{3(1+\sqrt{2})}$. Equality holds if and only if, up to an orthogonal transformation, $P$ is the gyroelongated square bipyramid with vertices
\begin{align*}
a_+=(0,0,1), \quad a_-=(0,0,-1),
\end{align*}
and
\begin{equation}
    v_j=\left(r_*\cos\left(\tfrac{\pi}{8}+\tfrac{j\pi}{4}\right),r_*\sin\left(\tfrac{\pi}{8}+\tfrac{j\pi}{4}\right),(-1)^j h_*\right),\quad j\in\{0,1,\ldots,7\},
\end{equation}
where $r_*=\sqrt{1-h_*^2}$.
\end{theorem}
Therefore, once we prove Theorem \ref{thm:Class-12-maximizer}, the proof of Theorem \ref{mainThm} will be complete.

\subsection{Proof of Theorem \ref{thm:Class-12-maximizer}}

We break the proof into several steps.

\subsubsection{Step 0: Set up}

Let $P\in\mathcal{P}_{10}$ have Type 12 degree sequence $(5^8,4^2)$. Label the two degree-4 vertices by $a_+$ and $a_-$, and label the remaining eight vertices cyclically by $v_0,v_1,\ldots,v_7$. The polytope $P$ has 16 facets:
\begin{align*}
    &[a_+,v_0,v_2], \, [a_+,v_2,v_4], \, [a_+,v_4,v_6],\,[a_+,v_6,v_0],\\
    &[a_-,v_1,v_3], \, [a_-,v_3,v_5],\, [a_-,v_5,v_7],\,[a_-,v_7,v_1],
\end{align*}
together with the eight ``belt facets'' labeled cyclically by
\[
[v_i,v_{i+1},v_{i+2}],\quad i\in\mathbb{Z}/8\mathbb{Z}.
\]

\subsubsection{Step 1: The oriented volume decomposition}

We choose the orientation so that the four upper cap facets are
\[
[a_+,v_0,v_2],\,[a_+,v_2,v_4],\,[a_+,v_4,v_6],\,[a_+,v_6,v_0],
\]
and the lower cap facets have the opposite orientation. The contribution to $6\vol(P)$ of the belt facets is denoted by
\begin{equation}\label{eq:def-Vbelt}
V_{\mathrm{belt}}:=\sum_{i=0}^7(-1)^i\det(v_i,v_{i+1},v_{i+2}).
\end{equation}
Here and throughout the remainder of the proof, indices are understood to be taken modulo 8. By Lemma \ref{lem:simp-vol}, the volume of $P$ can be expressed as
\begin{equation}\label{eq:6volP-Type12}
\begin{aligned}
      6\vol(P)&=\langle a_+,v_0\times v_2+v_2\times v_4+v_4\times v_6+v_6\times v_0\rangle\\
      &-\langle a_-,v_1\times v_3+v_3\times v_5+v_5\times v_7+v_7\times v_1\rangle+V_{\mathrm{belt}}.
\end{aligned}
\end{equation}
We will pass to new coordinates more amenable to this problem by first pairing opposite belt vertices. For $i\in\{0,1,2,3\}$, we set
\[
m_i:=\frac{v_i+v_{i+4}}{2}\quad\text{and}\quad d_i:=\frac{v_i-v_{i+4}}{2}.
\]
Then $v_i=m_i+d_i$, $v_{i+4}=m_i-d_i$, and since $v_i,v_{i+4}\in\Sp$, we have $\|m_i+d_i\|^2=\|m_i-d_i\|^2=1$. Subtracting and adding these equations together, we get
\begin{equation}\label{eq:midi}
    \langle m_i,d_i\rangle=0\quad\text{and}\quad \|m_i\|^2+\|d_i\|^2=1.
\end{equation}
In particular, $\sum_{i=0}^3(\|m_i\|^2+\|d_i\|^2)=4$, and combining this with the cross product identity
\[
x_0\times x_1+x_1\times x_2+x_2\times x_3+x_3\times x_0=(x_0-x_2)\times(x_1-x_3),
\]
we get
\[
\sum_{i\in\{0,2,4,6\}}(v_i\times v_{i+2})=(v_0-v_4)\times(v_2-v_6)=4d_0\times d_2\quad\text{and}\quad \sum_{i\in\{1,3,5,7\}}(v_i\times v_{i+2})=4d_1\times d_3.
\]
Hence
\[
6\vol(P)=\langle 4a_+,d_0\times d_2\rangle-\langle 4a_-,d_1\times d_3\rangle+V_{\mathrm{belt}}.
\]
Since $\|a_+\|=\|a_-\|=1$, we get
\begin{equation}\label{eq:cap-estimates}
\begin{aligned}
    \langle 4a_+,d_0\times d_2\rangle-\langle 4a_-,d_1\times d_3\rangle &\leq 4\|d_0\times d_2\|+4\|d_1\times d_3\|\\
    &\leq 4\|d_0\|\,\|d_2\|+4\|d_1\|\,\|d_3\|
    \leq 2\sum_{i=0}^3\|d_i\|^2.
\end{aligned}
\end{equation}
The first inequality follows from two applications of the Cauchy--Schwarz inequality. The second inequality follows since $\|u\times v\|=\|u\|\,\|v\|\sin(\measuredangle(u,v))\leq \|u\|\,\|v\|$, and the final inequality follows from the elementary $p=2$ case of Young's inequality, namely $2ab\leq a^2+b^2$ where $a,b\geq 0$, applied to the two pairs $4\|d_0\|\,\|d_2\|\leq 2\|d_0\|^2+2\|d_2\|^2$ and $4\|d_1\|\,\|d_3\|\leq 2\|d_1\|^2+2\|d_3\|^2$.

Therefore,
\begin{equation}\label{eq:6volP-Step1}
    6\vol(P) \leq 2\sum_{i=0}^3\|d_i\|^2+V_{\mathrm{belt}}.
\end{equation}

\subsubsection{Step 2: Changing to orthogonal coordinates}

Borrowing from the ideas in \cite{Hoehner-Ledford-2026}, we change to orthogonal coordinates, adapted to the combinatorial class under consideration. Define
\begin{equation}\label{eq:def-Ais}
A_0:=\frac{m_0+m_1+m_2+m_3}{\sqrt{2}},\, A_1:=\frac{m_0-m_1+m_2-m_3}{\sqrt{2}},\, A_2:=m_0-m_2,\, A_3:=m_1-m_3.
\end{equation}
For $i\in\{0,1,2,3\}$, we also set $\theta_i=\frac{\pi}{8}+\frac{i\pi}{4}$, and then we define the coordinates
\begin{equation}\label{eq:def-Bi-from-di}
B_0:=\sum_{i=0}^3\cos\theta_i d_i,\, B_1:=\sum_{i=0}^3\sin\theta_i d_i,\,B_2:=\sum_{i=0}^3\cos(3\theta_i)d_i,\,B_3:=\sum_{i=0}^3\sin(3\theta_i)d_i.
\end{equation}
The two coefficient matrices, after division by $\sqrt{2}$, are orthogonal. Thus
\begin{equation}\label{eq:sum-squares}
    \sum_{j=0}^3\|A_j\|^2=2\sum_{i=0}^3\|m_i\|^2\quad\text{and}\quad\sum_{j=0}^3\|B_j\|^2=2\sum_{i=0}^3\|d_i\|^2.
\end{equation}
Hence
\begin{equation}\label{eq:sum-squares-equals-8}
    \sum_{j=0}^3\|A_j\|^2+\sum_{j=0}^3\|B_j\|^2=2\left(\sum_{i=0}^3\|m_i\|^2+\sum_{i=0}^3\|d_i\|^2\right)=2(4)=8.
\end{equation}
Expanding \eqref{eq:def-Vbelt} in these coordinates, we obtain
\begin{equation}\label{eq:Vbelt-main}
    V_{\mathrm{belt}}=\sqrt{2}\det(A_1,A_2,A_3)+\mathcal{M}(\mathbf{A},\mathbf{B}),
\end{equation}
where for $C=1+\sqrt{2}$, $D=\sqrt{2}-1$, and $\kappa=\frac{C}{2\sqrt{2}}=\frac{2+\sqrt{2}}{4}$, the quantity $\mathcal{M}(\mathbf{A},\mathbf{B})$ is defined by
\begin{equation}\label{eq:MAB-def}
    \begin{aligned}
        \mathcal{M}(\mathbf{A},\mathbf{B}):=&-\frac{1}{\sqrt{2}}\det(A_0,B_0,B_2)+\frac{1}{\sqrt{2}}\det(A_0,B_1,B_3)
        \\
        &+\frac{C}{\sqrt{2}}\det(A_1,B_0,B_1)+\frac{D}{\sqrt{2}}\det(A_1,B_2,B_3)\\
        &+\frac{1}{2}\det(A_2,B_0,B_2)-\frac{1}{2}\det(A_2,B_0,B_3)\\
        &+\frac{1}{2}\det(A_2,B_1,B_2)+\frac{1}{2}\det(A_2,B_1,B_3)\\
        &-\frac{1}{2}\det(A_3,B_0,B_2)-\frac{1}{2}\det(A_3,B_0,B_3)\\
        &+\frac{1}{2}\det(A_3,B_1,B_2)-\frac{1}{2}\det(A_3,B_1,B_3).
    \end{aligned}
\end{equation}

\subsubsection{Step 3: Estimating $\mathcal{M}(\mathbf{A},\mathbf{B})$}

The next lemma is used to estimate $\mathcal{M}(\mathbf{A},\mathbf{B})$ from above.

\begin{lemma}\label{lem:main-est-MAB}
    For any arbitrary vectors $A_0,\ldots,A_3,B_0,\ldots,B_3\in\R^3$, we have
    \begin{equation}\label{eq:main-estimate-MAB}
        |\mathcal{M}(\mathbf{A},\mathbf{B})|
        \leq\kappa\left(\sum_{j=0}^3\|B_j\|^2\right)\sqrt{\|A_1\|^2+\frac{9}{16}(\|A_0\|^2+\|A_2\|^2+\|A_3\|^2)}.
    \end{equation}
\end{lemma}

\begin{proof}
    Using the scalar triple product, we write $\mathcal{M}(\mathbf{A},\mathbf{B})=\sum_{j=0}^3\langle A_j,C_j\rangle$, where
    \begin{align*}
        C_0&=\frac{1}{\sqrt{2}}(-B_0\times B_2+B_1\times B_3),\quad C_1=\frac{C}{\sqrt{2}}(B_0\times B_1)+\frac{D}{\sqrt{2}}(B_2\times B_3),\\
        C_2&=\frac{1}{2}(B_0\times B_2-B_0\times B_3+B_1\times B_2+B_1\times B_3),\\
        C_3&=\frac{1}{2}(-B_0\times B_2-B_0\times B_3+B_1\times B_2-B_1\times B_3).
    \end{align*}
    Let $\mathbf{X}=(\tfrac{3}{4}A_0,A_1,\tfrac{3}{4}A_2,\tfrac{3}{4}A_3)\in(\R^3)^4$ and $\mathbf{Y}=(\tfrac{4}{3}C_0,C_1,\tfrac{4}{3}C_2,\tfrac{4}{3}C_3)\in(\R^3)^4$. Then
    \[
\langle \mathbf{X},\mathbf{Y}\rangle=\langle A_1,C_1\rangle+\sum_{j\neq 1}\left\langle\tfrac{3}{4}A_j,\tfrac{4}{3}C_j\right\rangle=\langle A_1,C_1\rangle+\sum_{j\neq 1}\langle A_j,C_j\rangle=\mathcal{M}(\mathbf{A},\mathbf{B}).
    \]
    Thus, by the Cauchy--Schwarz inequality,
    \begin{align*}
        |\mathcal{M}(\mathbf{A},\mathbf{B})| \leq \|\mathbf{X}\|\cdot\|\mathbf{Y}\|&=\sqrt{\|A_1\|^2+\sum_{j\neq 1}\|\tfrac{3}{4}A_j\|^2}\sqrt{\|C_1\|^2+\sum_{j\neq 1}\|\tfrac{4}{3}C_j\|^2}\\
        &=\sqrt{\|A_1\|^2+\frac{9}{16}\sum_{j\neq 1}\|A_j\|^2}\sqrt{\|C_1\|^2+\frac{16}{9}\sum_{j\neq 1}\|C_j\|^2}.
    \end{align*}

   \noindent We claim that
    \begin{equation}\label{eq:Cj-Bj-kappa-relation}
        \|C_1\|^2+\frac{16}{9}\sum_{j\neq 1}\|C_j\|^2\leq\kappa^2\left(\sum_{j=0}^3 \|B_j\|^2\right)^2.
    \end{equation}
    Combining this with the previous inequality for $|\mathcal{M}(\mathbf{A},\mathbf{B})|$ will prove the lemma.

    Let $B$ be the $4\times 3$ matrix
    \[
B=\begin{pmatrix}
    B_0^\top\\B_1^\top\\B_2^\top\\B_3^\top
\end{pmatrix},
    \]
    and consider the singular value decomposition $B=U\Sigma V^\top$ where $\Sigma=\diag(\sigma_1,\sigma_2,\sigma_3)$ and $\sigma_i\geq 0$, and the three columns of $U=\begin{pmatrix}
        u_1 & u_2 & u_3
    \end{pmatrix}$ are orthonormal vectors in $\R^4$. Without loss of generality, we may assume that $V\in\mathrm{SO}(3)$ since a sign can be absorbed into one of the $u_i$'s. In particular, we have $\sum_{j=0}^3\|B_j\|^2=\sigma_1^2+\sigma_2^2+\sigma_3^2$.

    For two vectors $u,v\in\R^4$, define their six $2\times 2$ minors by $p_{ij}(u,v)=u_i v_j-u_j v_i$, $0\leq i<j\leq 3$. Letting $u$ and $v$ be two columns of $U$, then expanding the cross products $B_i\times B_j$ in the orthonormal basis given by the columns of $V$, we get
    \begin{equation}\label{eq:singular-values}
\|C_1\|^2+\frac{16}{9}\sum_{j\neq 1}\|C_j\|^2=\tau_1\sigma_2^2\sigma_3^2+\tau_2\sigma_3^2\sigma_1^2+\tau_3\sigma_1^2\sigma_2^2
    \end{equation}
    where $\tau_1=\tau(u_2,u_3)$, $\tau_2=\tau(u_3,u_1)$, and $\tau_3=\tau(u_1,u_2)$ where
    \begin{align*}
        \tau(u,v)&:=\left(\frac{1+\sqrt{2}}{\sqrt{2}}p_{01}(u,v)+\frac{\sqrt{2}-1}{\sqrt{2}}p_{23}(u,v)\right)^2\\
        &+\frac{16}{9}\bigg[\frac{1}{2}\bigl(-p_{02}(u,v)+p_{13}(u,v)\bigr)^2+\frac{1}{4}\bigl(p_{02}(u,v)-p_{03}(u,v)+p_{12}(u,v)+p_{13}(u,v)\bigr)^2\\
        &+\frac{1}{4}\bigl(-p_{02}(u,v)-p_{03}(u,v)+p_{12}(u,v)-p_{13}(u,v)\bigr)^2\bigg].
    \end{align*}
    \noindent It remains to estimate the three numbers $\tau_i$.

    For an orthonormal pair $u,v\in\R^4$, write $p_{ij}=p_{ij}(u,v)$, and define
    \begin{alignat*}{2}
        x_1&=p_{01}+p_{23},  \quad y_1&&=p_{01}-p_{23},\\
        x_2&=p_{02}-p_{13},\quad y_2&&= p_{02}+p_{13},\\
        x_3&=p_{03}+p_{12}, \quad y_3&&=p_{03}-p_{12}.
    \end{alignat*}
    By Lagrange's identity, for the minors we get
    \[
\sum_{i<j}p_{ij}^2=\|u\|^2\|v\|^2-\langle u,v\rangle^2=1,
    \]
    and by definition of the $p_{ij}$, we also have $p_{01}p_{23}-p_{02}p_{13}+p_{03}p_{12}=0$. Thus,
    \[
x_1^2+x_2^2+x_3^2=\sum_{i<j}p_{ij}^2+2(p_{01}p_{23}-p_{02}p_{13}+p_{03}p_{12})=1,
    \]
    and similarly $y_1^2+y_2^2+y_3^2=1$, so $x=(x_1,x_2,x_3)$ and $y=(y_1,y_2,y_3)$ are unit vectors in $\R^3$. Now $p_{01}=\frac{x_1+y_1}{2}$ and $p_{23}=\frac{x_1-y_1}{2}$, so
    \[
\frac{1+\sqrt{2}}{\sqrt{2}}p_{01}+\frac{\sqrt{2}-1}{\sqrt{2}}p_{23}=x_1+\frac{y_1}{\sqrt{2}}.
    \]
    Furthermore, $\frac{-p_{02}+p_{13}}{\sqrt{2}}=-\frac{x_2}{\sqrt{2}}$, and
    \begin{align*}
        \frac{1}{2}(p_{02}-p_{03}+p_{12}+p_{13})&=\frac{y_2-y_3}{2}\\
        \frac{1}{2}(-p_{02}-p_{03}+p_{12}-p_{13})&=-\frac{y_2+y_3}{2}.
    \end{align*}
    Therefore,
    \begin{equation}\label{eq:tau-est-1}
        \tau(u,v)=\left(x_1+\frac{y_1}{\sqrt{2}}\right)^2+\frac{8}{9}(x_2^2+y_2^2+y_3^2).
    \end{equation}

    Now that we have $\tau$ in a convenient form, we will estimate it from above. Since $x,y\in\Sp$ we have $x_1^2\leq 1-x_2^2$ and $y_2^2+y_3^2=1-y_1^2$. Hence
    \[
\tau(u,v)\leq \left(x_1+\frac{y_1}{\sqrt{2}}\right)^2+\frac{8}{9}(2-x_1^2-y_1^2)=\frac{16}{9}+\frac{1}{9}x_1^2-\frac{7}{18}y_1^2+\sqrt{2}x_1 y_1.
    \]
    For fixed $y_1$, the right-hand side is a quadratic function of $x_1$. Hence, for $x_1\in[-1,1]$, its maximum occurs at an endpoint of this interval. Thus, without loss of generality we may choose the sign of $x_1$ to agree with that of $y_1$, and letting $y=|y_1|\in[0,1]$ we get
    \[
\tau(u,v)\leq\frac{17}{9}+\sqrt{2}y-\frac{7}{18}y^2.
    \]
    The derivative of the right-hand side is $\sqrt{2}-\frac{7}{9}y$, which is positive because $y\in[0,1]$. Thus, its maximum occurs at $y=1$, which implies that
    \[
0\leq \tau(u,v)\leq \frac{3}{2}+\sqrt{2}=:L.
    \]
    Therefore, $\tau_i\in[0,L]$ for $i\in\{1,2,3\}$.

    Next, we show that
    \begin{equation}\label{eq:sum-tau-2}
        \tau_1+\tau_2+\tau_3\leq 2L.
    \end{equation}
    For $u,v\in\R^4$, define
    \begin{align*}
        X(u,v)&:=(p_{01}+p_{23},p_{02}-p_{13},p_{03}+p_{12})\\
        Y(u,v)&:=(p_{01}-p_{23},p_{02}+p_{13},p_{03}-p_{12}).
    \end{align*}
    If $w$ is also a vector in $\R^4$, then expanding the inner products we obtain
    \begin{align*}
    \langle
    X(u,v),X(u,w)\rangle&=\|u\|^2\langle v,w\rangle-\langle u,v\rangle\langle u,w\rangle    \\
       \langle Y(u,v),Y(u,w)\rangle&=\|u\|^2\langle v,w\rangle-\langle u,v\rangle\langle u,w\rangle .
    \end{align*}
    Since $u_1$, $u_2$, and $u_3$ are orthonormal, it follows that $X(u_2,u_3)$, $X(u_3,u_1)$, and $X(u_1,u_2)$ form an orthonormal basis of $\R^3$. For example, $\|X(u_2,u_3)\|=1$ and the preceding identities imply
    \[
\langle X(u_2,u_3),X(u_3,u_1)\rangle = -\langle X(u_3,u_2),X(u_3,u_1)\rangle=0,
    \]
and the other pairs are handled in the same way. Similarly, $Y(u_2,u_3)$, $Y(u_3,u_1)$, and $Y(u_1,u_2)$ form another orthonormal basis of $\R^3$.   For shorthand, let $x^{(1)}=X(u_2,u_3)$, $x^{(2)}=X(u_3,u_1)$, $x^{(3)}=X(u_1,u_2)$, $y^{(1)}=Y(u_2,u_3)$, $y^{(2)}=Y(u_3,u_1)$, and $y^{(3)}=Y(u_1,u_2)$. By \eqref{eq:tau-est-1}, for $i\in\{1,2,3\}$ we have
\[
\tau_i=\left(x_1^{(i)}+\frac{y_1^{(i)}}{\sqrt{2}}\right)^2+\frac{8}{9}\left((x_2^{(i)})^2+(y_2^{(i)})^2+(y_3^{(i)})^2\right).
\]
Summing over $i$ and using the orthonormality of the $x^{(i)}$ and the $y^{(i)}$, we obtain
\[
\sum_{i=1}^3 (x_1^{(i)})^2=1,\quad\sum_{i=1}^3 (y_1^{(i)})^2=1,\quad \sum_{i=1}^3 (x_2^{(i)})^2=1,\quad \text{and}\quad \sum_{i=1}^3 \left((y_2^{(i)})^2+(y_3^{(i)})^2\right)=2.
\]
Hence
\begin{align*}
\tau_1+\tau_2+\tau_3&=\sum_{i=1}^3 (x_1^{(i)})^2+\sqrt{2}\sum_{i=1}^3 x_1^{(i)}y_1^{(i)}+\frac{1}{2}\sum_{i=1}^3 (y_1^{(i)})^2\\
&+\frac{8}{9}\left(\sum_{i=1}^3 (x_2^{(i)})^2+\sum_{i=1}^3 \left((y_2^{(i)})^2+(y_3^{(i)})^2\right)\right)\\
&=\frac{3}{2}+\sqrt{2}\sum_{i=1}^3 x_1^{(i)}y_1^{(i)}+\frac{8}{3}=\frac{25}{6}+\sqrt{2}\sum_{i=1}^3 x_1^{(i)}y_1^{(i)}.
\end{align*}
By the Cauchy--Schwarz inequality,
\[
\sum_{i=1}^3 x_1^{(i)}y_1^{(i)}\leq\left|\sum_{i=1}^3 x_1^{(i)}y_1^{(i)}\right|\leq\sqrt{\sum_{i=1}^3 (x_1^{(i)})^2}\sqrt{\sum_{i=1}^3 (y_1^{(i)})^2}=1.
\]
Since $\sqrt{2}>7/6$, we have
\[
\tau_1+\tau_2+\tau_3\leq \frac{25}{6}+\sqrt{2}<3+2\sqrt{2}=2L.
\]
This proves \eqref{eq:sum-tau-2}.

Now we can return to \eqref{eq:singular-values}. Set $z_i=\sigma_i^2$, so that \eqref{eq:singular-values} becomes
\begin{equation}
    \|C_1\|^2+\frac{16}{9}\sum_{j\neq 1}\|C_j\|^2 = \tau_1 z_2 z_3+\tau_2 z_3 z_1+\tau_3 z_1 z_2.
\end{equation}
Without loss of generality, suppose that $z_3=\max\{z_1,z_2,z_3\}$. Then the two largest products are $z_3z_1$ and $z_3z_2$. Since $0\leq\tau_i\leq L$ and $\tau_1+\tau_2+\tau_3\leq 2L$, we have
\[
 \|C_1\|^2+\frac{16}{9}\sum_{j\neq 1}\|C_j\|^2\leq L(z_2z_3+z_1z_3).
\]
Moreover, we have $(z_3-z_1-z_2)^2\geq 0$, which is equivalent to $4z_3(z_1+z_2)\leq(z_1+z_2+z_3)^2$. Therefore,
\begin{align*}
     \|C_1\|^2+\frac{16}{9}\sum_{j\neq 1}\|C_j\|^2 &\leq Lz_3(z_1+z_2)\leq\frac{L}{4}(z_1+z_2+z_3)^2.
\end{align*}
Finally, since
\[
\frac{L}{4}=\frac{3+2\sqrt{2}}{8}=\left(\frac{1+\sqrt{2}}{2\sqrt{2}}\right)^2=\kappa^2\quad\text{and}\quad z_1+z_2+z_3=\sum_{j=0}^3\|B_j\|^2,
\]
we  obtain
\[
 \|C_1\|^2+\frac{16}{9}\sum_{j\neq 1}\|C_j\|^2\leq \kappa^2\left(\sum_{j=0}^3\|B_j\|^2\right)^2.
\]
This completes the proof of the lemma.
\end{proof}

\subsubsection{Step 4: Reduction to a two-variable optimization problem}

Set $U:=\|A_1\|$, $R^2=\|A_0\|^2+\|A_2\|^2+\|A_3\|^2$, and $E=\sum_{j=0}^3\|B_j\|^2$. Then the equation \eqref{eq:sum-squares-equals-8} becomes $U^2+R^2+E=8$. Returning to equation \eqref{eq:Vbelt-main} for $V_{\mathrm{belt}}$, the first term there can be estimated by
\begin{equation}\label{eq:A-terms-Vbelt}
\sqrt{2}\det(A_1,A_2,A_3)\leq \sqrt{2}U\|A_2\|\,\|A_3\|\leq\frac{U}{\sqrt{2}}(\|A_2\|^2+\|A_3\|^2)\leq\frac{UR^2}{\sqrt{2}}
\end{equation}
where we used Young's inequality. Combining \eqref{eq:6volP-Step1}, \eqref{eq:sum-squares}, \eqref{eq:Vbelt-main}, \eqref{eq:main-estimate-MAB}, and \eqref{eq:A-terms-Vbelt}, we derive that
\begin{equation}\label{eq:6volP-2vars}
    6\vol(P) \leq E+\frac{UR^2}{\sqrt{2}}+\kappa E\sqrt{U^2+\frac{9}{16}R^2}.
\end{equation}
\subsubsection{Step 5: Sharp optimization of \eqref{eq:6volP-2vars}}

Let $T=\sqrt{U^2+\frac{9}{16}R^2}$. Then $T\geq U$, $R^2=\frac{16}{9}(T^2-U^2)$, and
\[
E=8-U^2-R^2=\frac{72+7U^2-16T^2}{9}.
\]
Thus,  \eqref{eq:6volP-2vars} becomes
\[
6\vol(P)\leq \Phi(U,T):=\frac{72+7U^2-16T^2}{9}(1+\kappa T)+\frac{16}{9\sqrt{2}}U(T^2-U^2).
\]
Therefore, we want to solve the following maximization problem:
\begin{equation}\label{eq:2-var-max-problem}
    \max_{(U,T)\in\mathcal{S}}\Phi(U,T)
\end{equation}
where $\mathcal{S}=\{(U,T): U\geq 0, T\geq U, 72+7U^2-16T^2\geq 0\}$. Since $\Phi(U,T)$ is continuous on $\mathcal{S}$ and $\mathcal{S}$ is a compact set, the function $\Phi(U,T)$ attains a maximum on $\mathcal{S}$.

Suppose by way of contradiction that $\Phi(U,T)$ attains its maximum on the interior of $\mathcal{S}$, which is
\[
\operatorname{int}(\mathcal{S})=\{(U,T): U>0, T>U, 72+7U^2-16T^2>0\}.
\]
Setting the partial derivatives of $\Phi$ equal to zero, we get
\begin{equation}\label{eq:partial-derivatives-Phi}
\begin{aligned}
    \frac{\partial\Phi}{\partial U}&=\frac{14}{9}U+\frac{14}{9}\kappa TU+\frac{16}{9\sqrt{2}}T^2-\frac{48}{9\sqrt{2}}U^2=0\\
    \frac{\partial\Phi}{\partial T}&=8\kappa+\frac{7}{9}\kappa U^2-\frac{32}{9}T-\frac{48}{9}\kappa T^2+\frac{32}{9\sqrt{2}}UT=0.
\end{aligned}
\end{equation}
The equation $\frac{\partial\Phi}{\partial U}=0$ is equivalent to $28U=48\sqrt{2}U^2-16\sqrt{2}T^2-28\kappa TU$. For $(U,T)\in\operatorname{int}(\mathcal{S})$, we have $U>0$, so this equation is equivalent to
\[
28=U\left[48\sqrt{2}-16\sqrt{2}\left(\frac{T}{U}\right)^2-28\kappa\left(\frac{T}{U}\right)\right].
\]
Setting $q=T/U$, for $(U,T)\in\operatorname{int}(\mathcal{S})$  we have $T/U>1$, so the previous equation becomes $28=UD(q)$ where $D(q):=48\sqrt{2}-16\sqrt{2}q^2-28\kappa q$. Note that
\[
D'(q)=-32\sqrt{2}q-28\kappa=-32\sqrt{2}q-7(2+\sqrt{2})<0
\]
since $q\geq 1$. Moreover,
\[
D(13/10)=48\sqrt{2}-16\sqrt{2}\left(\frac{13}{10}\right)^2-28\kappa\left(\frac{13}{10}\right)=\frac{-910+593\sqrt{2}}{50}<0.
\]
Combining the previous two inequalities and using $D(q)>0$, we must have $1<q<13/10$.

Now we turn our attention to the second equation in \eqref{eq:partial-derivatives-Phi}, which we write as
\begin{equation}\label{eq:36PhiT}
    36\frac{\partial\Phi}{\partial T}=-(96+48\sqrt{2})T^2+64\sqrt{2}TU-128T+(14+7\sqrt{2})U^2+144+72\sqrt{2}=0.
\end{equation}
Plugging $U=28/D(q)$ and $T=qU$ into \eqref{eq:36PhiT}, we get
\[
36\frac{\partial\Phi}{\partial T}=\frac{16P(q)}{D(q)^2}
\]
where
\begin{align*}
    P(q)&:=(4608+2304\sqrt{2})q^4+(8064+9632\sqrt{2})q^3-(24806+11521\sqrt{2})q^2\\
    &-(24192+25760\sqrt{2})q+42158+21079\sqrt{2}.
\end{align*}
Thus, if $(U,T)\in\operatorname{int}(\mathcal{S})$ were a critical point, then we would need $P(q)=0$. Note that $P(1)=5832-4266\sqrt{2}<0$ and
\[
P(13/10)=-\frac{3}{2500}(280186+3447943\sqrt{2})<0.
\]
It suffices to show that $P$ cannot become nonnegative between these endpoints. Differentiating, we obtain
\begin{align*}
    P''(1)&=14(3862+4457\sqrt{2})>0\\
    P'''(q)&=192\left[(576+288\sqrt{2})q+252+301\sqrt{2}\right].
\end{align*}
Every term inside the brackets is positive since $q\geq 1$, so $P'''(q)>0$ for all $q\geq 1$. Hence $P''(q)$ is strictly increasing on $[1,\infty)$. Thus, since $P''(1)>0$, we have $P''(q)>0$ for all $q\geq 1$. Therefore, $P(q)$ is strictly convex on $[1,13/10]$. For each $q\in[1,13/10]$, we have $q=\lambda\cdot 1+(1-\lambda)\cdot\frac{13}{10}$ for some $\lambda\in[0,1]$. By the strict convexity of $P$ on this interval, we have
\[
\forall q\in[1,13/10],\quad P(q)\leq \lambda P(1)+(1-\lambda)P(13/10)<0.
\]
Hence $\frac{\partial\Phi}{\partial T}=\frac{16P(q)}{36D(q)^2}<0$, a contradiction to the necessary condition $\frac{\partial\Phi}{\partial T}=0$ required for a critical point. This proves that $\Phi$ has no critical points in $\operatorname{int}(\mathcal{S})$. Therefore, the maximum of $\Phi$ must be attained on the boundary of $\mathcal{S}$.

The boundary of $\mathcal{S}$ can be decomposed into three pieces: (i) the part on which $T=U$; (ii) the part on which $U=0$; and (iii) the part on which $E=0$. We begin with case (i) when $T=U$, which is equivalent to $R=0$. Since $72+7U^2-16T^2\geq 0$, we have $9U^2\leq 72$, so $0\leq U\leq 2\sqrt{2}$. Thus, we want to maximize $\Phi(U,U)=(8-U^2)(1+\kappa U)=:G(U)$
on the interval $0\leq U\leq 2\sqrt{2}$. Note that $G'(U)=-3\kappa U^2-2U+8\kappa$ and $G''(U)=-6\kappa U-2<0$, so $G$ is strictly concave on $[0,2\sqrt{2}]$ and has exactly one positive critical point on this interval. Setting $G'(U)=0$, we obtain the positive root $U_*=\frac{-1+\sqrt{24\kappa^2+1}}{3\kappa}\in[0,2\sqrt{2}]$, while the other root is negative. Therefore,
\begin{align*}
\max_{0\leq U\leq 2\sqrt{2}}\Phi(U,U)=G(U_*)=(8-U_*^2)(1+\kappa U_*)&=\frac{16\left(2+\sqrt{10+6\sqrt{2}}\right)^2}{9\left(1+\sqrt{10+6\sqrt{2}}\right)}\\
&=13.3122667893\ldots>12.
\end{align*}
Note that if we set $U=2\sqrt{2}h$, since $2\sqrt{2}\kappa =C$, the critical point equation can be expressed as $3Ch^2+2h-C=0$, whose unique positive solution is $h=h_*=\frac{\sqrt{10+6\sqrt{2}}-1}{3C}$, and thus
\[
\max_{\substack{(U,T)\in\mathcal{S}\\U=T}}\Phi(U,T)=8(1-h_*^2)(1+Ch_*).
\]

Next, we consider the boundary case (ii). Assume that $U=0$. In this case, the objective function becomes
\[
\Phi(0,T)=\left(8-\frac{16}{9}T^2\right)(1+\kappa T).
\]
Based on the computations in case (i), it suffices to prove that $\Phi(0,T)<12$. Indeed, note that $\Phi(0,T)=12-\frac{2}{9}p_0(T)$, where
\[
p_0(T):=(4+2\sqrt{2})T^3+8T^2-(18+9\sqrt{2})T+18.
\]
Since $T\geq 0$ on the feasible set, we have $p_0''(T)=6(4+2\sqrt{2})T+16\geq 16>0$, so $p_0$ is strongly convex with parameter 16. By the strong convexity of $p_0$ and Taylor's theorem at $T=9/10$, we get
\[
p_0(T) \geq p_0(9/10)+p_0'(9/10)(T-9/10)+8(T-9/10)^2.
\]
The minimum of the quadratic on the right-hand side is
\[
p_0(9/10)-\frac{1}{32}[p_0'(9/10)]^2=\frac{9(39797-22482\sqrt{2})}{40000}>0.
\]
Therefore, $\Phi(0,T)=12-\frac{2}{9}p_0(T)<12$.

Finally, we consider the boundary case (iii). Assume that $E=0$. Then $U^2+R^2=8$ and $\Phi(U,T)=UR^2/\sqrt{2}=U(8-U^2)/\sqrt{2}$. Hence, if $E=0$, then
\[
\Phi(U,T)\leq \max_{0\leq U\leq\sqrt{8}}\frac{U(8-U^2)}{\sqrt{2}}=\frac{32}{3\sqrt{3}}<12.
\]

Altogether, we have thus shown that
\[
6\vol(P) \leq \max_{(U,T)\in\mathcal{S}}\Phi(U,T) \leq \frac{16\left(2+\sqrt{10+6\sqrt{2}}\right)^2}{9\left(1+\sqrt{10+6\sqrt{2}}\right)},
\]
and equality requires $R=0$.
\subsubsection{Step 6: Equality yields the gyroelongated square bipyramid}

It remains to prove that equality forces the maximizer in $\mathcal{P}_{10}$ to be the gyroelongated square bipyramid. As shown in the last step, equality requires $R=0$, which implies  $A_0=A_2=A_3=o$. By \eqref{eq:def-Ais}, the condition $A_2=o$ implies that $m_0=m_2$, while the condition $A_3=o$ implies that $m_1=m_3$. Moreover, the condition $A_0=o$ implies that $\sum_{i=0}^3 m_i=o$, and combining this with $m_0=m_2$ and $m_1=m_3$ we deduce that $m_1=-m_0$. This also implies that $A_1=2\sqrt{2}m_0$, so $m_0=\frac{A_1}{2\sqrt{2}}$. Since $R=0$, we have $U=\|A_1\|$, and the maximum occurs when $U_*=2\sqrt{2}h_*>0$. Define $n=A_1/\|A_1\|\in\Sp$. Then $A_1=\|A_1\|n=2\sqrt{2}h_*n$. Combining these identities, we get $m_i=(-1)^i h_* n$ for $i\in\{0,1,2,3\}$.

Since $A_0=A_2=A_3=o$ and $A_1\neq o$, from \eqref{eq:MAB-def} we derive that
\begin{align*}
    \mathcal{M}(\mathbf{A},\mathbf{B}) &=\frac{C}{\sqrt{2}}\det(A_1,B_0,B_1)+\frac{D}{\sqrt{2}}\det(A_1,B_2,B_3)\\
    &=\left\langle A_1,\frac{C}{\sqrt{2}}(B_0\times B_1)+\frac{D}{\sqrt{2}}(B_2\times B_3)\right\rangle\\
    &\leq\frac{\|A_1\|}{\sqrt{2}}(C\|B_0\times B_1\|+D\|B_2\times B_3\|)\\
    &\leq\frac{\|A_1\|}{\sqrt{2}}(C\|B_0\|\,\|B_1\|+D\|B_2\|\,\|B_3\|)\\
    &\leq\frac{\|A_1\|}{2\sqrt{2}}\left(C\left(\|B_0\|^2+\|B_1\|^2\right)+D\left(\|B_2\|^2+\|B_3\|^2\right)\right)\\
    &\leq\frac{C\|A_1\|}{2\sqrt{2}}\sum_{i=0}^3\|B_i\|^2.
\end{align*}
The first inequality is Cauchy--Schwarz, the second is from $\|u\times v\|=\|u\|\,\|v\|\sin(\measuredangle(u,v))\leq \|u\|\,\|v\|$, and the third is a special case of Young's inequality. The last inequality follows since $D<C$, and thus equality can hold there if and only if $\|B_2\|=\|B_3\|=0$, i.e., if and only if $B_2=B_3=o$. Equality in $\|B_0\times B_1\|=\|B_0\|\,\|B_1\|$ requires $B_0$ to be orthogonal to $B_1$, equality in $2ab\leq a^2+b^2$ gives $\|B_0\|=\|B_1\|$, and equality in the first inequality forces $A_1$ to be parallel to $B_0\times B_1$. Hence we may choose an orthonormal basis, say $\{e_1,e_2,n\}$, with $e_1\times e_2=n$, such that $B_0=2r_* e_1$, $B_1=2r_* e_2$, and $A_1=2\sqrt{2}h_*n$ where $r_*=\sqrt{1-h_*^2}$.

Next, recall the orthogonal transformation in \eqref{eq:def-Bi-from-di}, which in matrix form says $\mathbf{B}=M\mathbf{d}$, where
\begin{align*}
    \mathbf{B}=\begin{pmatrix}
        B_0\\B_1\\B_2\\B_3
    \end{pmatrix}, \quad \mathbf{d}=\begin{pmatrix}
        d_0\\d_1\\d_2\\d_3
    \end{pmatrix},\quad\text{and}\quad M=\begin{pmatrix}
        \cos\theta_0 & \cos\theta_1 & \cos\theta_2 & \cos\theta_3\\
        \sin\theta_0 & \sin\theta_1 & \sin\theta_2 & \sin\theta_3\\
        \cos 3\theta_0 & \cos3\theta_1 & \cos 3\theta_2 & \cos 3\theta_3\\
        \sin 3\theta_0 & \sin 3\theta_1 & \sin 3\theta_2 & \sin 3\theta_3
    \end{pmatrix}.
\end{align*}
Since $MM^\top =2I_4$, we have $M^{-1}=\frac{1}{2}M^\top$, so $\mathbf{d}=M^{-1}\mathbf{B}=\frac{1}{2}M^\top\mathbf{B}$. Thus, for $i\in\{0,1,2,3\}$ we obtain
\begin{align*}
    d_i &=\frac{1}{2}(\cos\theta_i B_0+\sin\theta_i B_1+\cos 3\theta_i B_2+\sin 3\theta_i B_3)\\
    &=\frac{1}{2}(\cos\theta_i B_0+\sin\theta_i B_1)\quad\text{since }B_2=B_3=o\\
    &=r_*(\cos\theta_i e_1+\sin\theta_i e_2)\quad\text{since }B_0=2r_*e_1,\, B_1=2r_*e_2.
\end{align*}
Now since $v_i=m_i+d_i$, $v_{i+4}=m_i-d_i$, and $m_i=(-1)^i h_*n$, from this we get
\begin{align*}
    v_i&=(-1)^i h_*n+r_*(\cos\theta_i e_1+\sin\theta_i e_2)\\
    v_{i+4}&=(-1)^i h_*n-r_*(\cos\theta_i e_1+\sin\theta_i e_2)
\end{align*}
for $i\in\{0,1,2,3\}$. Extending the definition $\theta_j=\frac{\pi}{8}+\frac{j\pi}{4}$ to $j\in\{0,1,\ldots,7\}$, we have $\theta_{i+4}=\theta_i+\pi$, so $\cos\theta_{i+4}=-\cos\theta_i$ and $\sin\theta_{i+4}=-\sin\theta_i$. Also, $(-1)^{i+4}=(-1)^i$. Hence
\[
v_{i+4}=(-1)^{i+4} h_*n+r_*(\cos\theta_{i+4}e_1+\sin\theta_{i+4}e_2).
\]
Therefore, a single formula that holds for all of the vertices $v_j$ is
\begin{equation}\label{eq:vertices-final-form}
v_j = (-1)^j h_*n+r_*\left(\cos\left(\frac{\pi}{8}+\frac{j\pi}{4}\right)e_1+\sin\left(\frac{\pi}{8}+\frac{j\pi}{4}\right)e_2\right),\quad j\in\{0,1,\ldots,7\}.
\end{equation}

Finally, equality in the cap estimates in \eqref{eq:cap-estimates} requires
\[
a_+=\frac{d_0\times d_2}{\|d_0\times d_2\|}=n\quad\text{and}\quad a_-=-\frac{d_1\times d_3}{\|d_1\times d_3\|}=-n.
\]
This implies that $a_+=n$ and $a_-=-n$.
\subsubsection{Step 7: Computation of the volume of the maximizer}

By the rotation invariance of the sphere and the volume functional, without loss of generality we may assume that $n=e_3$. Every cap facet has facial cone volume $r_*^2/6$, and there are 8 such facets. For example, $\vol([o,a_+,v_0,v_2]) =\frac{1}{6}|\det(a_+,v_0,v_2)|$ and
\begin{align*}
    \det(a_+,v_0,v_2)&=\det\begin{pmatrix}
        0 & r_*\cos\tfrac{\pi}{8} & r_*\cos\tfrac{5\pi}{8} \\ 0 & r_*\sin\tfrac{\pi}{8} & r_*\sin\tfrac{5\pi}{8}\\ 1 & h_* & h_*
    \end{pmatrix}\\
    &=r_*^2\left(\cos\frac{\pi}{8}\sin\frac{5\pi}{8}-\cos\frac{5\pi}{8}\sin\frac{\pi}{8}\right)=r_*^2\sin\left(\frac{5\pi}{8}-\frac{\pi}{8}\right)=r_*^2.
\end{align*}
The remaining cap facets have the same volume $r_*^2/6$. There are 8 belt facets, and a similar determinant computation shows that each one has volume $\frac{1}{6}(1+\sqrt{2})h_* r_*^2$. Therefore,
\[
6V_*=8r_*^2+8(1+\sqrt{2})h_* r_*^2=8r_*^2\left(1+(1+\sqrt{2})h_*\right)=8(1-h_*^2)\left(1+(1+\sqrt{2})h_*\right).
\]
\subsubsection{Step 8: Verifying that the equality configuration lies in class 12}

To make the computations simpler, rotate the $e_1 e_2$-plane by the angle $\theta_1=3\pi/8$. This does not change any geometric property in the equality configuration. In this coordinate system, the belt vertices become
\begin{alignat*}{3}
v_0&=\left(\frac{r_*}{\sqrt{2}},-\frac{r_*}{\sqrt{2}},h_*\right),
&\qquad v_1&=(r_*,0,-h_*),
&\qquad v_2&=\left(\frac{r_*}{\sqrt{2}},\frac{r_*}{\sqrt{2}},h_*\right),\\
v_3&=(0,r_*,-h_*),
&\qquad v_4&=\left(-\frac{r_*}{\sqrt{2}},\frac{r_*}{\sqrt{2}},h_*\right),
&\qquad v_5&=(-r_*,0,-h_*),\\
v_6&=\left(-\frac{r_*}{\sqrt{2}},-\frac{r_*}{\sqrt{2}},h_*\right),
&\qquad v_7&=(0,-r_*,-h_*).
\end{alignat*}
and the apices are $a_+=(0,0,1)$ and $a_-=(0,0,-1)$.

A triangle $[p,q,s]$ is a facet of the maximizer $P_*$ if there is an affine function $\ell$ such that $\ell(p)=\ell(q)=\ell(s)=0$ and $\ell(w)<0$ for every other vertex $w$ of $P_*$. We will demonstrate the computations for two representative facets: the cap facet $F_c:=[a_+,v_0,v_2]$ and the belt facet $F_b:=[v_0,v_1,v_2]$.

First consider the cap facet $F_c$. The vertices $v_0$ and $v_2$ have the same $x$- and $z$-coordinates and opposite $y$-coordinates. Consequently, the equation of the support plane containing $F_c$ is independent of $y$. For a point $(x,y,z)\in\R^3$, define the affine support function $\ell_c(x,y,z):=\sqrt{2}(1-h_*)x+r_*z-r_*$. Note that $\ell_c(a_+)=\ell_c(v_0)=\ell_c(v_2)=0$, which shows that $\{(x,y,z)\in\R^3: \ell_c(x,y,z)=0\}$ is the  plane containing $F_c$. Next, we have $\ell_c(v_1)=r_*[\sqrt{2}-1-(1+\sqrt{2})h_*]$, which is negative if and only if $h_*>\frac{\sqrt{2}-1}{\sqrt{2}+1}$. This does, in fact, hold true, as we will verify at the end of the argument. Moreover, $\ell_c(v_5)=-r_*[\sqrt{2}+1-(\sqrt{2}-1)h_*]<0$ and $\ell_c(a_-)=-2r_*<0$. We also have $\ell_c(v_3)=\ell_c(v_7)=-r_*(1+h_*)<0$. Similarly, $\ell_c(v_4)=\ell_c(v_6)=2r_*(h_*-1)<0$ since $h_*<1$.

Next, consider the facet $F_b$, and define $\ell_b(x,y,z):=2\sqrt{2}h_*x+r_*(\sqrt{2}-1)z-(\sqrt{2}+1)h_* r_*$. Simple computations yield $\ell_b(v_0)=\ell_b(v_1)=\ell_b(v_2)=0$, so $\{(x,y,z)\in\R^3: \ell_b(x,y,z)=0\}$ is the  plane containing $F_b$. We have $\ell_b(a_+)=r_*[\sqrt{2}-1-(1+\sqrt{2})h_*]$, which is negative if and only if $h_*>\frac{\sqrt{2}-1}{\sqrt{2}+1}$. We also have $\ell_b(a_-)=r_*[1-\sqrt{2}-(1+\sqrt{2})h_*]<0$ since $h_*>0$. Moreover, $\ell_b(v_3)=\ell_b(v_7)=-2\sqrt{2}h_* r_*<0$, $\ell_b(v_4)=\ell_b(v_6)=-4h_* r_*<0$, and $\ell_b(v_5)=-4\sqrt{2}h_* r_*<0$.

By rotational and reflectional symmetry of the configuration, similar computations hold for the remaining 14 facets. To complete the proof, it remains to verify that $h_*>\frac{\sqrt{2}-1}{\sqrt{2}+1}=3-2\sqrt{2}$. On one hand, note that $3-2\sqrt{2}<1/5$. On the other hand, recall that $h_*=\frac{\sqrt{10+6\sqrt{2}}-1}{3(\sqrt{2}+1)}$ is the positive root of $p(h)=3(\sqrt{2}+1)h^2+2h-(\sqrt{2}+1)=0$. Since $p'(h)=6(\sqrt{2}+1)h+2>0$, the function $p$ is increasing on $(0,\infty)$. Since $p(1/5)=-\frac{12+22\sqrt{2}}{25}<0$, $p(h_*)=0$, and $p$ is increasing on $(0,\infty)$, we have $h_*>1/5$. All ten stated points are distinct points of $\Sp$, so each is a vertex of their convex hull. A convex 3-polytope with ten vertices has at most $16$ facets by Euler's formula. Hence the sixteen supporting triangles exhibited above are all its facets, and the configuration has Type 12 combinatorial type. This completes the proof of Theorem \ref{thm:Class-12-maximizer}. \qed

\section{Eleven vertices}\label{sec:N11}

We prove Theorem \ref{mainThm-N11}. By Lemma \ref{lem:Wimmer-small-reductions}, it is enough to exclude the degree sequence $(6^2,5^6,4^3)$. We first give an explicit lower bound in the remaining class, and then apply the tangent plane method to the two graphs with this competing degree sequence. In the exclusions below, $P$ is assumed to be a global volume maximizer, so $o\in\operatorname{int}(P)$ by Lemma \ref{lem:origin-interior}.

\subsection{Identification of Type 34}\label{subsec:Type34-identification}

It remains to identify Type 34, whose degree sequence is $(6,5^8,4^2)$. We show in the following lemma that this combinatorial type has a natural description as an
``edge-contracted icosahedron''.

\begin{lemma}\label{lem:type34-edge-contracted-icosahedron}
A Type 34 polytope is combinatorially equivalent to the polytope obtained by contracting an edge of an icosahedron.
\end{lemma}

\begin{proof}
Let $G$ denote the graph of the icosahedron, and let $[u,v]$ be any
edge of $G$. Since $G$ is a triangulation of the sphere, the edge
$[u,v]$ belongs to precisely two triangular faces, say $[u,v,x]$ and $[u,v,y]$. Moreover, in the icosahedral graph every vertex has degree 5, and the only common neighbors of the adjacent vertices $u$ and
$v$ are $x$ and $y$.

Contract the edge $[u,v]$ to a single vertex $w$, suppressing the two
pairs of parallel edges that arise from the common neighbors $x$
and $y$. Thus, the resulting planar graph has $12-1=11$ vertices. The new vertex $w$ is adjacent precisely to the union
of the neighbors of $u$ and $v$, with $u$ and $v$ themselves
removed. Each of $u$ and $v$ has four neighbors other than the
other endpoint of $[u,v]$, and they have exactly the two common
neighbors $x$ and $y$. Hence $\deg(w)=4+4-2=6$.

The vertices $x$ and $y$ each lose one edge under the contraction:
the two edges $[x,u]$ and $[x,v]$ are identified to a single edge $[x,w]$, and
similarly the two edges $[y,u]$ and $[y,v]$ are identified to a single edge
$[y,w]$. Consequently, $\deg(x)=\deg(y)=4$.
None of the other eight vertices changes degree, so each of them
still has degree $5$. Thus the contracted graph has degree sequence $(6,5^8,4^2)$.

The contraction removes the two faces incident with the edge $[u,v]$, namely $[u,v,x]$ and $[u,v,y]$, while each of the remaining $18$ triangular faces of the icosahedron becomes a triangular face
of the contracted graph. Hence the resulting graph is again a
triangulation of the sphere. Equivalently, it has $V=11$, $E=27$, and  $F=18$. In particular, it is a 3-connected planar triangulation and hence,
by Steinitz's theorem, is the graph of a simplicial convex 3-polytope.

By Lemma~\ref{lem:Wimmer-small-reductions},
there is exactly one combinatorial type of simplicial convex
3-polytope with 11 vertices, minimum degree at least 4, and
degree sequence $(6,5^8,4^2)$. This is Type 34. This proves that Type 34 is combinatorially equivalent to the edge-contracted icosahedron.
\end{proof}

We next provide an explicit realization of Type 34 whose volume
will also serve as a convenient lower bound for the maximum volume
in $\mathcal{P}_{11}$.

\begin{lemma}\label{lem:type34-explicit-competitor}
Let $p=(1,0,0)$, and for $\varepsilon,\delta\in\{-1,1\}$, define the points
\[
\begin{alignedat}{2}
    a_\varepsilon
    &:=\left(\frac{49}{100},\varepsilon\frac{\sqrt{7599}}{100},0\right),
   \qquad& b_\varepsilon
    &:=\left(-\frac{13}{25},
    \varepsilon\frac{2\sqrt{114}}{25},0\right),\\
    c_\delta&:=
    \left(-\frac{4}{5},0,\delta\frac35
    \right), \qquad&
    d_{\varepsilon,\delta}&:=
    \left(\frac{1}{5},
    \frac{\varepsilon}{2},
    \delta\frac{\sqrt{71}}{10}
    \right).
\end{alignedat}
\]
Then the convex hull
\[
    P_{\mathrm{can}}:=
    \conv\big\{
    p,a_+,a_-,b_+,b_-,c_+,c_-,
    d_{++},d_{+-},d_{-+},d_{--}
    \big\}
\]
is a simplicial polytope of Type 34 and $o\in\operatorname{int}(P_{\mathrm{can}})$. Moreover,
\begin{equation}\label{eq:explicit-Type34-volume}
\vol(P_{\mathrm{can}})=
\frac{18}{125}
+\frac{4\sqrt{114}}{125}
+\frac{3\sqrt{71}}{50}
+\frac{43\sqrt{8094}}{6250}
+\frac{19\sqrt{539529}}{18750},
\end{equation}
and in particular $\vol(P_{\mathrm{can}})>2.3545$.
\end{lemma}

\begin{proof}
First, we show that every stated point lies on $\Sp$. Indeed, we have $\|p\|=1$ and:
\begin{align*}
    \|a_\varepsilon\|^2&=\frac{49^2+7599}{100^2}=\frac{2401+7599}{10000}=1\\
    \|b_\varepsilon\|^2&=\frac{13^2+4(114)}{25^2}=\frac{169+456}{625}=1\\
    \|c_\delta\|^2&=
    \frac{16}{25}+\frac9{25}=1\\    \|d_{\varepsilon,\delta}\|^2&=
    \frac{1}{25}+\frac14+\frac{71}{100}
    =\frac{4+25+71}{100}=1.
\end{align*}
We also have
\[
31p+50a_++50a_-+50c_++50c_-=o.
\]
The coefficients are all positive, and these five points affinely span $\R^3$. Thus $o$ lies in the interior of their convex hull, and thus in the interior of $P_{\mathrm{can}}$.

Next, let $\sigma_y(x,y,z)=(x,-y,z)$ and $\sigma_z(x,y,z)=(x,y,-z)$. The set of vertices of $P_{\mathrm{can}}$ is invariant under the group generated by $\sigma_y$ and $\sigma_z$. An elementary
calculation using supporting planes shows that the facets of $P_{\mathrm{can}}$ are precisely the images under these two reflections of the following six triangular facets:
\begin{center}
\begin{tabular}{c|c}
Representative facet & Orbit size\\
\hline
$[p,a_+,d_{++}]$ & $4$\\
$[p,d_{++},d_{-+}]$ & $2$\\
$[b_+,c_+,c_-]$ & $2$\\
$[b_+,d_{++},c_+]$ & $4$\\
$[b_+,d_{++},a_+]$ & $4$\\
$[c_+,d_{++},d_{-+}]$ & $2$
\end{tabular}
\end{center}
Thus, $P_{\mathrm{can}}$ has
\[
    4+2+2+4+4+2=18
\]
triangular facets. For completeness, the supporting plane
assertion may be checked by forming the plane through the three
vertices in each of the six representative triples, substituting
the remaining eight vertices, and checking the signs of all inner products (see the supplied program {\tt verify\_small.py} on Zenodo and Appendix~\ref{sec:appendix-verification-programs}); in each case, all remaining vertices
lie strictly in the same open halfspace. By symmetry, the same
statement holds for every member of the corresponding orbit. All eleven stated points are distinct points of $\Sp$, so each is a vertex of their convex hull. Since a convex 3-polytope with eleven vertices has at most $18$ facets, the listed supporting triangles are all its facets.

From the preceding facet list, the vertex $p$ is incident with
six facets, each of $a_+$ and $a_-$ is incident with four facets,
and each of
\[
    b_+,b_-,c_+,c_-,
    d_{++},d_{+-},d_{-+},d_{--}
\]
is incident with five facets. Since $P_{\mathrm{can}}$ is simplicial, this shows that $P_{\mathrm{can}}$ has
the vertex degree sequence $d(P_{\mathrm{can}})=(6,5^8,4^2)$.
Therefore, by Lemma~\ref{lem:Wimmer-small-reductions}, $P_{\mathrm{can}}$ is a Type 34 polytope. Finally, by Lemma~\ref{lem:type34-edge-contracted-icosahedron}, we deduce that $P_{\mathrm{can}}$ is an edge-contracted icosahedron.

It remains to compute the volume of $P_{\mathrm{can}}$. Since the origin lies in the interior of $P_{\mathrm{can}}$, the volume is the sum of the volumes of the tetrahedral cones over its facets. Thus,
\[
    \vol(P_{\mathrm{can}})
    =\frac{1}{6}\sum_{F=[x,y,z]}
    |\det(x,y,z)|
\]
where the sum is taken over all facets $F=[x,y,z]$ of $P_{\mathrm{can}}$. By the two reflection symmetries, it is enough to compute the determinant for the six representative facets above. A straightforward calculation yields the following data:
\[
\begin{array}{c|c|c}
\text{Representative facet}
&\text{Multiplicity}
&|\det(x,y,z)|\\
\hline
[p,a_+,d_{++}]
&4&
\dfrac{\sqrt{539529}}{1000}\\[2mm]
[p,d_{++},d_{-+}]
&2&
\dfrac{\sqrt{71}}{10}\\[2mm]
[b_+,c_+,c_-]
&2&
\dfrac{48\sqrt{114}}{625}\\[2mm]
[b_+,d_{++},c_+]
&4&
\dfrac{6\sqrt{114}}{625}
+\dfrac{39}{250}
+\dfrac{4\sqrt{8094}}{625}\\[2mm]
[b_+,d_{++},a_+]
&4&
\dfrac{49\sqrt{8094}}{12500}
+\dfrac{13\sqrt{539529}}{25000}\\[2mm]
[c_+,d_{++},d_{-+}]
&2&
\dfrac{3}{25}+\dfrac{2\sqrt{71}}{25}.
\end{array}
\]
Hence
\begin{align*}
6\vol(P_{\mathrm{can}})
&=4\frac{\sqrt{539529}}{1000}+2\frac{\sqrt{71}}{10}+2\frac{48\sqrt{114}}{625}
+4\left(\frac{6\sqrt{114}}{625}+\frac{39}{250}
+\frac{4\sqrt{8094}}{625}\right)\\
&+4\left(\frac{49\sqrt{8094}}{12500}
+\frac{13\sqrt{539529}}{25000}\right)+2\left(\frac{3}{25}+\frac{2\sqrt{71}}{25}\right).
\end{align*}
Collecting like terms and dividing by 6, we get
\[
\vol(P_{\mathrm{can}})
=\frac{18}{125}+\frac{4\sqrt{114}}{125}+\frac{3\sqrt{71}}{50}+\frac{43\sqrt{8094}}{6250}+\frac{19\sqrt{539529}}{18750},
\]
as claimed.

Finally, we use the following estimates to obtain a lower bound on the volume of $P_{\mathrm{can}}$:
\[
\begin{aligned}
(10.677)^2&=113.998329<114,\\
(8.4261)^2&=70.99916121<71,\\
(89.9666)^2&=8093.98911556<8094,\\
(734.5263)^2&=539528.88539169<539529.
\end{aligned}
\]
Therefore,
\[
\sqrt{114}>10.677,\qquad
\sqrt{71}>8.4261,\qquad
\sqrt{8094}>89.9666,\qquad
\sqrt{539529}>734.5263.
\]
Substituting these lower bounds into
\eqref{eq:explicit-Type34-volume}, we obtain
\begin{align*}
\vol(P_{\mathrm{can}})&>\frac{18}{125}+\frac4{125}(10.677)
+\frac{3}{50}(8.4261)
+\frac{43}{6250}(89.9666)
+\frac{19}{18750}(734.5263)\\
&=\frac{18394689}{7812500}
=2.354520192>2.3545.
\end{align*}
This completes the proof.
\end{proof}

\subsection{Type 32}
\label{subsec:Type32}

We first eliminate Type 32. Type 33 will be handled in the next subsection. Let $P\in\mathcal{P}_{11}$ have Type $32$ combinatorial type. Then its degree sequence is $d(P)=(6^2,5^6,4^3)$, and, by Table~\ref{tab:facial-degree-patterns-N11}, its 18 facets have the facial degree multiset
\begin{equation}\label{eq:Type32-facial-pattern}
    8(456)+4(455)+2(555)+4(556).
\end{equation}

Suppose that a facet has facial degree pattern $(456)$. If the
corresponding half-angles are $\theta_4$, $\theta_5$, and $\theta_6$, then we
write $W_{(456)}:=W(\theta_4,\theta_5,\theta_6)$.
Likewise, for a facet of pattern $(455)$ we write $W_{(455)}:=W(\theta_4,\theta_5,\theta_5')$. We use analogous notation for the patterns $(555)$ and $(556)$, and reserve $s_{(ijk)}$ for the fixed supporting points. To eliminate Type $32$, we use the following supporting
half-angle triples in Table~\ref{tab:type32-supporting-triples}.

\begin{table}[ht]
\centering
\renewcommand{\arraystretch}{1.15}
\begin{tabular}{c|c}
Facial degree pattern $(ijk)$
& Supporting half-angle triple $s_{(ijk)}=(s_i,s_j,s_k)$\\
\hline
$(456)$ & $(0.807473169,\;0.591467910,\;0.508075579)$\\
$(455)$ & $(0.741248151,\;0.601367680,\;0.601367680)$\\
$(555)$ & $(0.645572627,\;0.645572627,\;0.645572627)$\\
$(556)$ & $(0.679179430,\;0.679179430,\;0.554645170)$
\end{tabular}
\vspace{2mm}
\caption{Supporting half-angle triples for Type $32$.}
\label{tab:type32-supporting-triples}
\end{table}

Let us demonstrate the subsequent computations using the pattern
$(456)$ and the supporting triple
\[
    s_{(456)}
    =
    (0.807473169,0.591467910,0.508075579).
\]
At this point,
\begin{align*}
t_0&=
0.807473169+0.591467910+0.508075579-\frac{\pi}{2}=0.3362203312\ldots,\\
t_1&=\frac{\pi}{2}+0.807473169-0.591467910-0.508075579=1.2787260068\ldots,\\
t_2&=\frac{\pi}{2}
-0.807473169+0.591467910-0.508075579=0.8467154888\ldots,\\
t_3&=\frac{\pi}{2}-0.807473169-0.591467910+0.508075579=0.6799308268\ldots.
\end{align*}
In particular, each $t_i$ lies in the interval $(0,\pi)$, so the supporting point lies in the domain of $W$. Moreover,
\[
    S=\sum_{i=0}^3\cot t_i=5.2830821464\ldots,
\]
which implies
\[
    f=\frac{1}{S}=0.1892834471\ldots
\]
and
\[
    W(s_{(456)})=\frac{2}{3}f=0.1261889647\ldots.
\]
Substituting the preceding values of $t_0,t_1,t_2,t_3$ into
\eqref{eq:partials-of-W}, we obtain
\begin{equation}\label{eq:Type32-456-partials}
\begin{aligned}
W_{\theta_1}(s_{(456)})&=0.1424974702\ldots,\\
W_{\theta_2}(s_{(456)})&=0.1755363129\ldots,\\
W_{\theta_3}(s_{(456)})&=0.2112513118\ldots.
\end{aligned}
\end{equation}
Therefore,
\begin{align*}
C_{(456)}&=W(s_{(456)})
-W_{\theta_1}(s_{(456)})(0.807473169)\\
&\qquad-W_{\theta_2}(s_{(456)})(0.591467910)
-W_{\theta_3}(s_{(456)})(0.508075579)
=-0.2000296478\ldots.
\end{align*}

Continuing in the same fashion for the other three supporting
points gives the following tangent plane data in Table~\ref{tab:type32-tangent-data}.

\begin{table}[ht]
\centering
\renewcommand{\arraystretch}{1.15}
\setlength{\tabcolsep}{5pt}
\begin{tabular}{c|cccc}
Pattern $(ijk)$
& $W_{\theta_1}$
& $W_{\theta_2}$
& $W_{\theta_3}$
& $C_{(ijk)}$\\
\hline
$(456)$
& $0.1424974702$
& $0.1755363129$
& $0.2112513118$
& $-0.2000296478$\\
$(455)$
& $0.1424974715$
& $0.1755363106$
& $0.1755363106$
& $-0.1801969718$\\
$(555)$
& $0.1755363131$
& $0.1755363131$
& $0.1755363131$
& $-0.2030570125$\\
$(556)$
& $0.1755363053$
& $0.1755363053$
& $0.2112513137$
& $-0.2244696492$
\end{tabular}
\vspace{2mm}
\caption{Values of the partial derivatives and tangent plane
intercepts for the supporting half-angle triples used for Type
32.}
\label{tab:type32-tangent-data}
\end{table}

For completeness, the strict bounds needed below are verified by the outward-rounded Taylor intervals described in Appendix \ref{sec:certificates}. The displayed approximations in Table \ref{tab:type32-tangent-data} serve only to illustrate the calculation. At each supporting point we obtain
\begin{equation}\label{eq:Type32-rigorous-slope-bounds}
\begin{aligned}
W_{\theta_j}(s)&<0.142499
&&\text{whenever the $j$th entry of $s$ corresponds to degree $4$},\\
W_{\theta_j}(s)&<0.175538
&&\text{whenever the $j$th entry of $s$ corresponds to degree $5$},\\
W_{\theta_j}(s)&<0.211253
&&\text{whenever the $j$th entry of $s$ corresponds to degree $6$},
\end{aligned}
\end{equation}
for each of the four supporting points in
Table~\ref{tab:type32-supporting-triples}. The same interval
calculation yields
\begin{equation}\label{eq:Type32-rigorous-intercept-bounds}
\begin{aligned}
C_{(456)}&<-0.200028,\\
C_{(455)}&<-0.180195,\\
C_{(555)}&<-0.203055,\\
C_{(556)}&<-0.224468.
\end{aligned}
\end{equation}

Now set $\lambda_4:=0.142499$, $\lambda_5:=0.175538$, and $\lambda_6:=0.211253$. Since all spherical half-angles are positive, replacing the partial derivatives in each supporting tangent plane by the larger degree-dependent quantities $\lambda_4,\lambda_5,\lambda_6$ preserves the upper bound. Likewise, replacing the tangent plane
intercepts by the larger quantities in
\eqref{eq:Type32-rigorous-intercept-bounds} preserves the upper
bound. Consequently, we obtain:
\begin{equation}\label{eq:Type32-facial-tangent-bounds}
\begin{aligned}
W_{(456)}&\leq-0.200028+\lambda_4\theta_4+\lambda_5\theta_5+\lambda_6\theta_6,\\
W_{(455)}&\leq-0.180195+\lambda_4\theta_4+\lambda_5\theta_5+\lambda_5\theta_5',\\
W_{(555)}&\leq-0.203055+\lambda_5\theta_5+\lambda_5\theta_5'+\lambda_5\theta_5'',\\
W_{(556)}&\leq-0.224468+\lambda_5\theta_5+\lambda_5\theta_5'+\lambda_6\theta_6.
\end{aligned}
\end{equation}

A Type 32 polytope has 18 facets. Since $o\in\operatorname{int}(P)$, its volume is the sum of the volumes of its $18$ facial cones, i.e., $\vol(P)=\sum_F\vol([o,F])=\sum_F W_F$. By \eqref{eq:Type32-facial-pattern}, the sum of the constant
terms in \eqref{eq:Type32-facial-tangent-bounds} is
\begin{align*}
8(-0.200028)+4(-0.180195)+2(-0.203055)+4(-0.224468)=-3.624986.
\end{align*}

We now consider the linear terms. For each facet $F$ and each
vertex $v\in F$, let $\theta_{F,v}:=\frac{1}{2}\alpha_{F,v}$, where $\alpha_{F,v}$ is the spherical interior angle of the radial
spherical triangle corresponding to $F$ at $v$. After summing all
the facet inequalities, every corner incident with a degree-$4$
vertex has coefficient $\lambda_4$, every corner incident with a
degree-$5$ vertex has coefficient $\lambda_5$, and every corner
incident with a degree-$6$ vertex has coefficient $\lambda_6$.
Thus, the total linear contribution is
\[
    \sum_{v\in\operatorname{Vert}(P)}
    \lambda_{\deg(v)}
    \sum_{F\ni v}\theta_{F,v}.
\]
Moreover, the radial projections of the facets form a spherical
triangulation of $\Sp$. Hence, the spherical angles at the facets
incident with any fixed vertex $v$ sum to $2\pi$, and therefore
\[
    \sum_{F\ni v}\theta_{F,v}
    =\frac{1}{2}\sum_{F\ni v}\alpha_{F,v}=\pi.
\]
Consequently, the total linear contribution is $\pi\sum_{v\in\operatorname{Vert}(P)}\lambda_{\deg(v)}$.
Since Type 32 has degree sequence $(6^2,5^6,4^3)$, we get
\[
    \sum_{v\in\operatorname{Vert}(P)}\lambda_{\deg(v)}
    =3\lambda_4+6\lambda_5+2\lambda_6.
\]
Using the chosen values of the $\lambda_j$, we have
\begin{align*}
3\lambda_4+6\lambda_5+2\lambda_6&=3(0.142499)+6(0.175538)+2(0.211253)=1.903231.
\end{align*}

Putting everything together and using
$\pi<3.1415927$, we obtain
\begin{align*}
\vol(P)\leq-3.624986+\pi(1.903231)
&<-3.624986+(3.1415927)(1.903231)\\
&=2.3541906160\ldots<2.354191.
\end{align*}
Therefore,
\[
    \vol(P)<2.354191<2.3545<\vol(P_{\mathrm{can}}),
\]
which contradicts the maximality of $P$. Hence a global
volume maximizer in $\mathcal{P}_{11}$ cannot have Type 32 combinatorial type.

\subsection{Type 33}

Let $P\in\mathcal P_{11}$ have Type 33 combinatorial type. Its facial degree multiset is
\[
    2(446)+4(455)+4(456)+2(555)+6(556).
\]
Take $(\lambda_4,\lambda_5,\lambda_6)=(0.13272,0.17672,0.21065)$ and use the supporting points in Table \ref{tab:type33-short}. Each displayed decimal is exact. Substitution in \eqref{eq:partials-of-W} gives partial derivatives smaller than the corresponding $\lambda_j$, and the last column is an upper bound for the tangent intercept.

\begin{table}[ht]
\centering
\renewcommand{\arraystretch}{1.12}
\begin{tabular}{c|c|c}
Pattern & $s_{(ijk)}$ & $\overline C_{(ijk)}$\\
\hline
$(446)$ & $(0.716190,0.716190,0.501639)$ & $-0.16427$\\
$(455)$ & $(0.778791,0.580505,0.580505)$ & $-0.17413$\\
$(456)$ & $(0.861213,0.551621,0.482830)$ & $-0.19222$\\
$(555)$ & $(0.645031,0.645031,0.645031)$ & $-0.20531$\\
$(556)$ & $(0.677404,0.677404,0.558098)$ & $-0.22570$
\end{tabular}
\caption{Supporting tangent plane data for Type 33.}
\label{tab:type33-short}
\end{table}

Thus Lemma \ref{lem:tangent-plane-bound} gives
\begin{align*}
    \vol(P)&\leq-3.55876+\pi(3\lambda_4+6\lambda_5+2\lambda_6)\\
    &=-3.55876+1.87978\pi<2.346744<2.3545.
\end{align*}
Therefore, Type 33 is eliminated. Together with the Type 32 estimate and Lemma \ref{lem:Wimmer-small-reductions}, this proves that Type 34 is the combinatorial type of every global volume maximizer in $\mathcal P_{11}$. This completes the proof of Theorem \ref{mainThm-N11}. \qed

\section{Thirteen, fourteen, fifteen, and sixteen vertices}\label{sec:larger-cases}

In this section, we prove Theorem \ref{thm:larger-combinatorial-reductions}. The strategy is to use the supporting line bound \eqref{eq:degree-supporting-line-bound} and the supporting tangent plane bounds \eqref{eq:degree-tangent-plane-bound} and \eqref{eq:vertex-tangent-plane-bound}. The exact inputs and the verification programs are supplied in the accompanying Zenodo archive; see Appendix \ref{sec:certificates}.

\subsection{Lower bounds and the exclusion of trivalent vertices}

Wimmer \cite{Wimmer1998} gave the coordinate lists for the candidate configurations with $13$, $15$, and $16$ vertices in \cite[Chapter VI, pp.~67, 69--70]{Wimmer1998}. To obtain exact configurations from these printed data, we treat each coordinate as an exact terminating decimal and normalize each vector separately to length one. Using determinant intervals, we verify all listed strict supporting facets, with the origin on their interior sides. They give the following rigorous lower bounds:
\begin{equation}\label{eq:larger-candidate-lower-bounds}
\begin{aligned}
    \vol(P_{13,\mathrm{can}})>2.6128, \quad
    \vol(P_{15,\mathrm{can}})>2.8043, \quad\text{and}\quad 
    \vol(P_{16,\mathrm{can}})>2.8864.
\end{aligned}
\end{equation}
For fourteen vertices, we use Wimmer's symmetric gyroelongated hexagonal bipyramid, whose volume exceeds $2.7209$; see Remark \ref{rem:N14}.

On the other hand, by Lemma \ref{lem:Wimmer-trivalent} we get the following upper bounds if a trivalent vertex is present:
\begin{center}
\begin{tabular}{c|c|c}
$N$ & Upper bound & Supporting value $e_0$\\
\hline
13 & $2.594432$ & $1.30303856$\\
14 & $2.696426$ & $1.16892033$\\
15 & $2.786468$ & $1.05971858$\\
16 & $2.866481$ & $0.96911430$
\end{tabular}
\end{center}
Here the upper bounds are certified bounds, not rounded optimization outputs. More precisely, if we set
\[
    Q_N(e):=\Psi_3(e)+(2N-7)U\left(\frac{4\pi-e}{2N-7}\right),
\]
then $Q_N$ is concave, and
\[
    Q_N(e)\leq Q_N(e_0)+Q_N'(e_0)(e-e_0)
    \leq Q_N(e_0)+4\pi|Q_N'(e_0)|.
\]
Outward-rounded evaluation of the last expression proves the displayed bounds. Each is smaller than the corresponding lower bound. Therefore, a maximizing polytope in each of these four cases has no trivalent vertex.

\subsection{The exhaustive reduction}

We generate the nontrivalent planar triangulations using \texttt{plantri}~5.8 \cite{plantri58,BrinkmannMcKay}. There are: $525$, $2472$, $12400$, and $65619$ graphs for $N=13,14,15,16$, respectively. For each graph, we record its degree sequence and its facial degree patterns. First, we apply \eqref{eq:degree-supporting-line-bound}. For each graph not excluded by its degree sequence, we apply \eqref{eq:degree-tangent-plane-bound}. In the fifteen-vertex case, one final graph requires the vertex-dependent form \eqref{eq:vertex-tangent-plane-bound}.

\begin{table}[ht]
\centering
\renewcommand{\arraystretch}{1.15}
\begin{tabular}{l|rrrr}
$N$ & 13 & 14 & 15 & 16\\
\hline
Nontrivalent triangulations & 525 & 2472 & 12400 & 65619\\
Excluded by degree sequence & 490 & 2454 & 12324 & 65575\\
Excluded by degree-dependent planes & 34 & 17 & 74 & 42\\
Excluded by vertex-dependent planes & 0 & 0 & 1 & 0\\
Remaining combinatorial types & 1 & 1 & 1 & 2
\end{tabular}
\caption{The exhaustive combinatorial reductions.}
\label{tab:larger-reductions}
\end{table}

Every exclusion in Table \ref{tab:larger-reductions} has an interval-checked upper bound strictly below the appropriate lower bound. The verifier checks every graph in the supplied lists and verifies that the certificates match its degree sequence or its faces. The graphs that remain have the following data:

\begin{table}[ht]
\centering
\renewcommand{\arraystretch}{1.15}
\begin{tabular}{c|c|l}
$N$ & Degree sequence & Facial degree patterns\\
\hline
13 & $(6^2,5^{10},4)$ & $4(456)+10(555)+8(556)$\\
14 & $(6^2,5^{12})$ & $12(555)+12(556)$\\
15 & $(6^3,5^{12})$ & $8(555)+18(556)$\\
16 & $(6^4,5^{12})$ & $8(555)+16(556)+4(566)$\\
16 & $(6^4,5^{12})$ & $4(555)+24(556)$
\end{tabular}
\caption{Remaining graphical data for the cases $13\leq N\leq 16$.}
\label{tab:remaining-graph-data}
\end{table}

There is exactly one graph with the displayed degree sequence for each of $N=13,14,15$, and exactly two for $N=16$. This proves Theorem \ref{thm:larger-combinatorial-reductions}. In the next two subsections, we show the explicit details for the closest exclusions.

\subsection{Two thirteen-vertex competitors}

Wimmer \cite[Chapter IV, pp.~51--53]{Wimmer1998} identified the leading thirteen-vertex graph, and proved its uniqueness for its degree sequence. He also discusses two competitors with degree sequence $(6^3,5^8,4^2)$. His discussion does not supply a complete exclusion of all remaining graphs. The preceding exhaustive verification supplies the missing step in his argument. The two closest competitors, denoted here by $B$ and $C$, have facial degree multisets
\begin{align*}
    B:\quad&2(455)+4(456)+2(466)+6(555)+6(556)+2(566),\\
    C:\quad&3(455)+4(456)+(466)+3(555)+10(556)+(566).
\end{align*}
They are records $1$ and $130$ in the supplied thirteen-vertex list, with records numbered starting at zero.

For $B$, take $(\lambda_4,\lambda_5,\lambda_6)=(0.204387,0.228580,0.251712)$ and use Table \ref{tab:N13-B}. For $C$, take $(\lambda_4,\lambda_5,\lambda_6)=(0.200929,0.228211,0.252080)$ and use Table \ref{tab:N13-C}. All supporting coordinates and all intercept bounds in the tables are exact terminating decimals.

\begin{table}[ht]
\centering
\renewcommand{\arraystretch}{1.12}
\begin{tabular}{c|c|c|c}
Pattern & Multiplicity & $s_{(ijk)}$ & $\overline C_{(ijk)}$\\
\hline
$(455)$ & 2 & $(0.71579295,0.57854274,0.57854274)$ & $-0.28779704$\\
$(456)$ & 4 & $(0.77175234,0.57712352,0.50269723)$ & $-0.30030501$\\
$(466)$ & 2 & $(0.88229502,0.46791780,0.46791780)$ & $-0.31244175$\\
$(555)$ & 6 & $(0.62332912,0.62332912,0.62332912)$ & $-0.30395561$\\
$(556)$ & 6 & $(0.65470044,0.65470044,0.54742599)$ & $-0.31748244$\\
$(566)$ & 2 & $(0.71687341,0.56444047,0.56444047)$ & $-0.33156311$
\end{tabular}
\caption{Supporting tangent plane data for the thirteen-vertex competitor $B$.}
\label{tab:N13-B}
\end{table}

\begin{table}[ht]
\centering
\renewcommand{\arraystretch}{1.12}
\begin{tabular}{c|c|c|c}
Pattern & Multiplicity & $s_{(ijk)}$ & $\overline C_{(ijk)}$\\
\hline
$(455)$ & 3 & $(0.72925467,0.57165905,0.57165905)$ & $-0.28487371$\\
$(456)$ & 4 & $(0.79224900,0.56442546,0.49179523)$ & $-0.29757301$\\
$(466)$ & 1 & $(0.92642533,0.43970405,0.43970405)$ & $-0.30965115$\\
$(555)$ & 3 & $(0.62346945,0.62346945,0.62346945)$ & $-0.30326520$\\
$(556)$ & 10 & $(0.65565411,0.65565411,0.54538903)$ & $-0.31719960$\\
$(566)$ & 1 & $(0.72077782,0.56214935,0.56214935)$ & $-0.33171188$
\end{tabular}
\caption{Supporting tangent plane data for the thirteen-vertex competitor $C$.}
\label{tab:N13-C}
\end{table}

Substituting the corresponding data in \eqref{eq:partials-of-W} and using outward rounding verifies every derivative and intercept inequality needed for Lemma \ref{lem:tangent-plane-bound}. Therefore,
\[
    \vol(P_B)<2.607921<2.6128 \quad\text{and}\quad 
    \vol(P_C)<2.605773<2.6128.
\]
The other thirty-two facial exclusions are verified in the same way by the supplied certificates. The surviving graph is record $0$. Its geometric optimization remains to be done; we leave this as an open problem for further research.

\subsection{The final fifteen-vertex exclusion}

The leading fifteen-vertex graph has degree sequence $(6^3,5^{12})$. The last competing graph, record $0$ in the supplied list, has degree sequence $(6^4,5^{10},4)$ and facial degree multiset
\[
    2(456)+2(466)+6(555)+14(556)+2(566).
\]
A degree-dependent angle optimization gives approximately $2.8044377673$, which is above the explicit candidate volume. To fix this issue, for this graph we will use a separate coefficient $\lambda_v$ at every vertex. The twenty-six supporting triples, their intercept bounds, and the fifteen coefficients are listed in Appendix \ref{subsec:N15-vertex-data}. By Equation \eqref{eq:vertex-tangent-plane-bound}, we get
\begin{equation}\label{eq:N15-final-exclusion}
    \vol(P)\leq\sum_F\overline C_F+\pi\sum_v\lambda_v
    <2.803886151<2.8043.
\end{equation}
This completes the fifteen-vertex exclusion. The existence of only one graph with degree sequence $(6^3,5^{12})$ is a combinatorial fact; the volume argument above proves that every maximizer must have this degree sequence.

For sixteen vertices, both graphs with degree sequence $(6^4,5^{12})$ remain. The degree and tangent plane certificates exclude every other graph, but do not distinguish these two. Thus, these two cases remain open as well.

\begin{remark}\label{rem:N14}
For the case of fourteen vertices, Table \ref{tab:larger-reductions} recovers Wimmer's theorem that the maximizing polytope must be combinatorially equivalent to the gyroelongated hexagonal bipyramid (see \cite[Chapter IV, pp.~44--51]{Wimmer1998}). Thus, we gave an independent verification of his combinatorial conclusion. This polytope is depicted in Figure~\ref{fig:gyroelongated-hexagonal-bipyramid} below.

Wimmer also optimized the symmetric family. For antipodal poles and two regular hexagons at heights $h$ and $-h$, twisted by $\pi/6$, the volume is
\[
    V_{14}(h):=(1-h^2)[\sqrt{3}+(\sqrt{3}+2)h].
\]
Its unique positive critical height and corresponding volume are
\[
    h_* =\frac{\sqrt{6}+\sqrt{2}-1}{3+2\sqrt{3}}\quad\text{and}\quad
    V_{14}(h_*)=2.7209778993\ldots,
\]
respectively. It remains an open problem to optimize over arbitrary realizations of the combinatorial type's  graph.

There is a specific limitation of the present tangent plane constraints. Give each cap face half-angles $(u,u,\pi/6)$ and each belt face half-angles $(v,v,v)$, with $2u+3v=\pi$. Every vertex angle budget is satisfied. At the exact choice $u=0.65$, the relaxed objective satisfies
\[
    12W(u,u,\pi/6)+12W(v,v,v)>2.7311>V_{14}(h_*).
\]
Thus, these face domains and angle budgets alone do not prove that the symmetric candidate is optimal. They do not impose agreement of the spherical edge lengths calculated from neighboring faces. Wimmer encounters a related gap for hexagonal caps in \cite[Chapter III, pp.~26--29]{Wimmer1998}. A further geometric inequality is needed for the unrestricted fourteen-vertex volume theorem.
\end{remark}

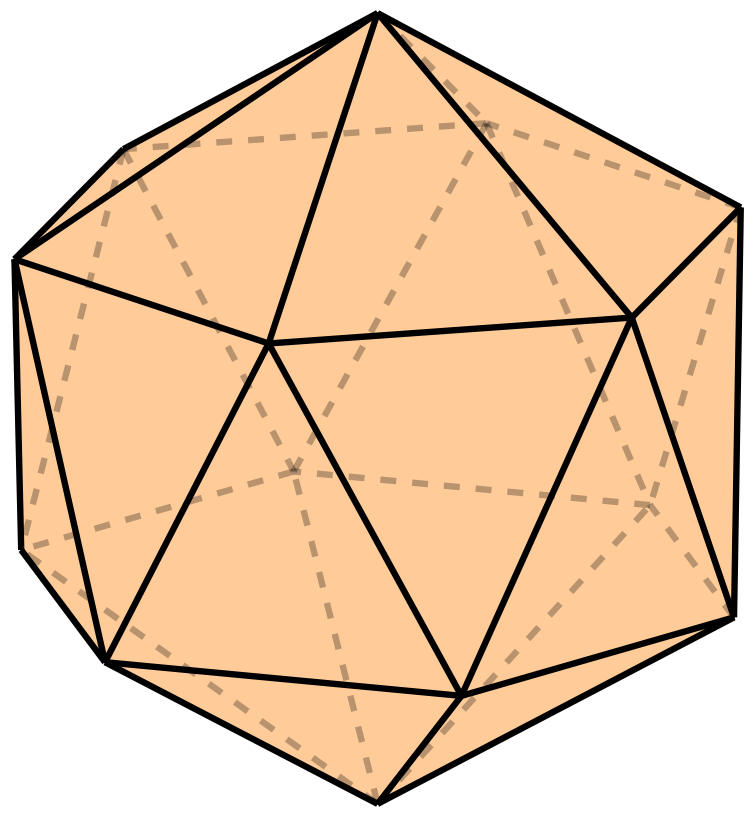
\begin{figure}[ht]
\begin{center}
\tdplotsetmaincoords{72}{122}

\begin{tikzpicture}[
    scale=2.75,
    line join=bevel,
    tdplot_main_coords
]

\def\h{0.443016443066}
\def\a{0.865965530063}
\def\b{0.232034764454}
\def\c{0.633930765610}

\coordinate (T) at (0,0,1);
\coordinate (D) at (0,0,-1);
\coordinate (A1) at (\a,\b,\h);
\coordinate (A2) at (\b,\a,\h);
\coordinate (A3) at (-\c,\c,\h);
\coordinate (A4) at (-\a,-\b,\h);
\coordinate (A5) at (-\b,-\a,\h);
\coordinate (A6) at (\c,-\c,\h);
\coordinate (B1) at (\a,-\b,-\h);
\coordinate (B2) at (\c,\c,-\h);
\coordinate (B3) at (-\b,\a,-\h);
\coordinate (B4) at (-\a,\b,-\h);
\coordinate (B5) at (-\c,-\c,-\h);
\coordinate (B6) at (\b,-\a,-\h);

\draw[fill opacity=0.4,fill=orange]
    (T)--(A5)--(A6)--cycle;

\draw[fill opacity=0.4,fill=orange]
    (A6)--(B6)--(B1)--cycle;

\draw[fill opacity=0.4,fill=orange]
    (D)--(B3)--(B2)--cycle;

\draw[fill opacity=0.4,fill=orange]
    (A2)--(B3)--(A3)--cycle;

\draw[fill opacity=0.4,fill=orange]
    (D)--(B2)--(B1)--cycle;

\draw[fill opacity=0.4,fill=orange]
    (T)--(A2)--(A3)--cycle;

\draw[fill opacity=0.4,fill=orange]
    (A2)--(B2)--(B3)--cycle;

\draw[fill opacity=0.4,fill=orange]
    (T)--(A6)--(A1)--cycle;

\draw[fill opacity=0.4,fill=orange]
    (A6)--(B1)--(A1)--cycle;

\draw[fill opacity=0.4,fill=orange]
    (T)--(A1)--(A2)--cycle;

\draw[fill opacity=0.4,fill=orange]
    (A1)--(B1)--(B2)--cycle;

\draw[fill opacity=0.4,fill=orange]
    (A1)--(B2)--(A2)--cycle;

\draw[very thick,dashed,opacity=0.28] (T)--(A4);
\draw[very thick,dashed,opacity=0.28] (D)--(B4);
\draw[very thick,dashed,opacity=0.28] (D)--(B5);
\draw[very thick,dashed,opacity=0.28] (D)--(B6);
\draw[very thick,dashed,opacity=0.28] (A3)--(A4);
\draw[very thick,dashed,opacity=0.28] (A3)--(B4);
\draw[very thick,dashed,opacity=0.28] (A4)--(A5);
\draw[very thick,dashed,opacity=0.28] (A4)--(B4);
\draw[very thick,dashed,opacity=0.28] (A4)--(B5);
\draw[very thick,dashed,opacity=0.28] (A5)--(B5);
\draw[very thick,dashed,opacity=0.28] (A5)--(B6);
\draw[very thick,dashed,opacity=0.28] (B3)--(B4);
\draw[very thick,dashed,opacity=0.28] (B4)--(B5);
\draw[very thick,dashed,opacity=0.28] (B5)--(B6);

\draw[very thick] (T)--(A1);
\draw[very thick] (T)--(A2);
\draw[very thick] (T)--(A3);
\draw[very thick] (T)--(A5);
\draw[very thick] (T)--(A6);
\draw[very thick] (D)--(B1);
\draw[very thick] (D)--(B2);
\draw[very thick] (D)--(B3);
\draw[very thick] (A1)--(A2);
\draw[very thick] (A1)--(A6);
\draw[very thick] (A1)--(B1);
\draw[very thick] (A1)--(B2);
\draw[very thick] (A2)--(A3);
\draw[very thick] (A2)--(B2);
\draw[very thick] (A2)--(B3);
\draw[very thick] (A3)--(B3);
\draw[very thick] (A5)--(A6);
\draw[very thick] (A6)--(B1);
\draw[very thick] (A6)--(B6);
\draw[very thick] (B1)--(B2);
\draw[very thick] (B1)--(B6);
\draw[very thick] (B2)--(B3);

\end{tikzpicture}
\end{center}
\caption{The combinatorial type of every global volume maximizer in $\mathcal{P}_{14}$ is that of a gyroelongated hexagonal bipyramid.}
\label{fig:gyroelongated-hexagonal-bipyramid}
\end{figure}

\section{Thirty-two vertices}\label{sec:N32}

We now prove Theorem \ref{mainThm-N32}. The proof uses a supporting plane estimate that can be summed over every triangulated graph.

\subsection{A universal facial inequality}

Set
\begin{equation}\label{eq:N32-constants}
    s:=\left(\frac{\pi}{5},\frac{\pi}{6},\frac{\pi}{6}\right),\quad
    w:=W(s)=\frac{\sqrt{10(5-\sqrt{5})}}{90},\quad \text{and} \quad A:=1.44684,
\end{equation}
and define
\begin{equation}\label{eq:N32-slopes}
\begin{aligned}
    \lambda_3&:=0.4447612,\qquad \lambda_4:=0.4500217,\\
    \lambda_5&:=\frac{5+7\sqrt{5}}{45},\qquad
    \lambda_6:=\sqrt{\frac{130+58\sqrt{5}}{1215}},\\
    \lambda_d&:=0.4614439\quad\text{for every integer }d\geq7.
\end{aligned}
\end{equation}
Substituting these values in \eqref{eq:partials-of-W} and using $\cos(\pi/5)=(1+\sqrt{5})/4$, we get
\[
    \nabla W(s)=(\lambda_5,\lambda_6,\lambda_6).
\]

\begin{lemma}\label{lem:N32-universal-face}
Let $d_1,d_2,d_3\geq3$ be integers, and let $\theta=(\theta_1,\theta_2,\theta_3)$ be in the domain of $W$. Then
\begin{equation}\label{eq:N32-universal-face}
\begin{aligned}
    W(\theta)\leq w+A\left(\sum_{j=1}^3\frac{1}{d_j}-\frac{8}{15}\right)
    +\sum_{j=1}^3\lambda_{d_j}\left(\theta_j-\frac{\pi}{d_j}\right).
\end{aligned}
\end{equation}
The inequality is strict unless the degree multiset is $(566)$. In that case, equality holds precisely when the corresponding half-angles are $(\pi/5,\pi/6,\pi/6)$.
\end{lemma}

\begin{proof}
For $p=(d_1,d_2,d_3)$, define
\[
    K_p:=w-\frac{8A}{15}+\sum_{j=1}^3\frac{A-\pi\lambda_{d_j}}{d_j}.
\]
Then the desired inequality is $W(\theta)-\sum_j\lambda_{d_j}\theta_j\leq K_p$. For $p=(566)$, it is exactly the supporting plane at $s$, since $1/5+1/6+1/6=8/15$. Strict concavity then gives the asserted equality condition.

For any other triple with entries in $\{3,4,5,6,7\}$, choose the supporting point $q$ in Table \ref{tab:N32-supporting-points}. Let $g_j=W_{\theta_j}(q)$ and $C=W(q)-\sum_jg_jq_j$. Applying concavity again, we get
\begin{equation}\label{eq:N32-residual-bound}
    W(\theta)-\sum_j\lambda_{d_j}\theta_j
    \leq C+\frac{\pi}{2}\sum_j\max\{0,g_j-\lambda_{d_j}\},
\end{equation}
because $0<\theta_j<\pi/2$. The thirty-four rows in the table satisfy
\begin{equation}\label{eq:N32-finite-slack}
    C+\frac{\pi}{2}\sum_j\max\{0,g_j-\lambda_{d_j}\}
    <K_p-\frac{1}{20000}.
\end{equation}
These are finite strict inequalities at specified points, verified by outward-rounded intervals in Appendix \ref{sec:certificates}. The smallest certified slack that arises is greater than $0.0000606593$.

For the seven rows marked $\mathcal B_d$, the exact supporting point is
\[
    \epsilon:=10^{-4},\qquad a_d:=2\sqrt{\lambda_d+\tfrac13},\qquad
    \mathcal B_d:=\left(\frac{\pi}{2}-a_d\epsilon,\epsilon,\epsilon\right).
\]
These are interior points since $a_d<2$, and the four coordinates $t_i$ are respectively $(2-a_d)\epsilon$, $\pi-(a_d+2)\epsilon$, $a_d\epsilon$, $a_d\epsilon$, all positive. Thus,  no boundary limit is required for \eqref{eq:N32-finite-slack}.

Finally, $A-\pi\lambda_7<-0.0028287<0$. For $d\geq7$ we have $\lambda_d=\lambda_7$, and
\[
    \frac{A-\pi\lambda_d}{d}\geq\frac{A-\pi\lambda_7}{7}.
\]
Replacing every degree greater than $7$ by $7$ leaves the slopes unchanged and decreases or preserves $K_p$. Therefore, the corresponding verified inequality proves \eqref{eq:N32-universal-face} for every integer degree triple. Permuting corners causes no issues because $W$ is symmetric. This proves the lemma.
\end{proof}

\subsection{Summation and the global bound}

Let $P_{32}^*\in\mathcal P_{32}$ have maximum volume. As a volume maximizer, $P_{32}^*$ has exactly $32$ vertices, is simplicial, and contains $o$ in its interior. More generally, let $P$ be any such polytope with $N$ vertices, and note that $P$ must have $F=2N-4$ facets. Apply Lemma \ref{lem:N32-universal-face} to each facet, using the degrees of its incident vertices. Every degree-$d$ vertex contributes $d$ copies of $1/d$, so
\[
    \sum_{F\in\mathcal F_2(P)}\sum_{v\in F}\frac{1}{\deg(v)}=N.
\]
Also, for each vertex $v$, the sum of its half-angles is $\pi$, which implies
\[
    \sum_{F\ni v}\lambda_{\deg(v)}\left(\theta_{F,v}-\frac{\pi}{\deg(v)}\right)=0.
\]
Summing \eqref{eq:N32-universal-face}, we therefore obtain
\begin{equation}\label{eq:N32-affine-preliminary}
    \vol(P)\leq(2N-4)\left(w-\frac{8A}{15}\right)+AN.
\end{equation}
When $N=32$, this is $\vol(P)\leq 60w=V_{32}$. In particular, it bounds the volume of the global maximizer, and hence the volume of every member of $\mathcal{P}_{32}$.

\subsection{Attainment and equality}

We begin with a regular icosahedron inscribed in $\Sp$ and add the radial projections of its twenty face centers. Taking the convex hull, the resulting polytope has $32$ vertices, $90$ edges, and $60$ triangular facets. Its twelve original vertices have degree $5$, and the twenty added vertices each have degree $6$. Every facet has pattern $(566)$, and its spherical angles are $(2\pi/5,\pi/3,\pi/3)$.

For completeness, its volume can also be calculated by a Gram determinant. Let $r=\sqrt{5}$. In one facet, the scalar product of the two degree-$6$ vertices is $b=r/3$, while each degree-$5$/degree-$6$ scalar product is positive and has square $\ell^2=(5+2r)/15$. Thus, the squared determinant of the three unit vertices is
\[
    1-2\ell^2-b^2+2\ell^2b=\frac{2(5-\sqrt{5})}{45}.
\]
Summing the sixty facial cone volumes, we get $\frac{60}{6}\sqrt{2(5-\sqrt{5})/45}=V_{32}$. To be rigorous, the exact-radical coordinates and interval checks of all strict supporting faces are included in the certificates. As noted in the introduction, this configuration and its volume are already included in Wimmer's dissertation \cite{Wimmer1998}.

Now suppose that equality holds. Then every facet must have pattern $(566)$ and half-angles $(\pi/5,\pi/6,\pi/6)$. In particular, every vertex has degree $5$ or $6$, and Euler's formula gives twelve degree-$5$ vertices. Each degree-$5$ vertex is surrounded by five congruent spherical triangles. Removing it and merging its star produces a regular spherical pentagon whose sides have equal length and its five interior angles are all $2\pi/3$.

At a degree-$6$ vertex, the neighbors alternate between degree $5$ and degree $6$, since each incident triangle contains exactly one degree-$5$ vertex. Thus, three of the pentagons meet at every remaining vertex. The resulting tiling is the regular dodecahedral spherical tiling. Indeed, after fixing one of its regular pentagons, the congruent pentagon across any one of its edges is uniquely determined on the other side. Propagation across edges determines the connected tiling. The removed degree-$5$ vertices are the spherical centers of its pentagons. Consequently, the configuration is precisely the icosahedron together with its normalized face centers, up to an orthogonal transformation. This proves Theorem \ref{mainThm-N32}. \qed

\begin{corollary}\label{cor:universal-affine-bound}
For every integer $N\geq4$ and every $P\in\mathcal P_N$,
\begin{equation}\label{eq:universal-affine-bound}
    \vol(P)\leq V_{32}+(N-32)\sigma \quad\text{where}\quad 
    \sigma:=2w-\frac{A}{15}=0.0203731360\ldots>0.
\end{equation}
Equality is possible only for $N=32$ and the configuration in Theorem \ref{mainThm-N32}.
\end{corollary}

\begin{proof}
Apply \eqref{eq:N32-affine-preliminary} to a global maximizer in $\mathcal{P}_N$, and rewrite the result using $V_{32}=60w$. Note that equality requires every face to have pattern $(566)$, which gives twelve degree-$5$ vertices and $60$ facets, hence $N=32$. The equality statement then follows from Theorem \ref{mainThm-N32}.
\end{proof}

\appendix
\section{Concavity of the facial cone function}\label{sec:appendix}

\subsection{Proof of Lemma \ref{lem:Types8-12-main-Lemma-concavity}}

Let $\mathbf{t}=(t_0,t_1,t_2,t_3)\in\mathcal{D}$, where $\mathcal{D}$ is the feasible set defined in \eqref{eq:Types8-12-feasible-set}. We shall prove that the Hessian matrix of $f$ is negative definite on $\mathcal{D}$. We begin by computing the second-order partial derivatives of $f$.

Since $t_3=\pi-t_0-t_1-t_2$, we have
\[
f(\mathbf{t})=f(t_0,t_1,t_2)=\bigl(\cot t_0+\cot t_1+\cot t_2+\cot(\pi-t_0-t_1-t_2)\bigr)^{-1}.
\]
Let $j\in\{0,1,2\}$. The first-order partial derivatives of $f$ are:
\begin{align*}
    \frac{\partial f}{\partial t_j}&=-\left(\sum_{i=0}^3\cot t_i\right)^{-2}\frac{\partial}{\partial t_j}\bigl(\cot t_0+\cot t_1+\cot t_2+\cot(\pi-t_0-t_1-t_2)\bigr)\\
    &=-\left(\sum_{i=0}^3\cot t_i\right)^{-2}\left[-\csc^2 t_j-\csc^2(\pi-t_0-t_1-t_2)\frac{\partial}{\partial t_j}(\pi-t_0-t_1-t_2)\right]\\
    &=\left(\sum_{i=0}^3\cot t_i\right)^{-2}\left(\csc^2 t_j-\csc^2 t_3\right)=f(\mathbf{t})^2\left(\csc^2 t_j-\csc^2 t_3\right).
\end{align*}

Now let $k\in\{0,1,2\}$, and denote the Kronecker delta by
\[
\delta_{jk}=\begin{cases}
    1, &\quad\text{if }j=k,\\
    0, &\quad\text{if }j\neq k.
\end{cases}
\]
The second-order partial derivatives of $f$ are:
\begin{align*}
    f_{jk}=\frac{\partial^2 f}{\partial t_k\partial t_j} &=2f(\mathbf{t})\left(\csc^2 t_j-\csc^2 t_3\right)\frac{\partial f}{\partial t_k}+f(\mathbf{t})^2\frac{\partial}{\partial t_k}\left(\csc^2 t_j-\csc^2 t_3\right)\\
    &=2f(\mathbf{t})\left(\csc^2 t_j-\csc^2 t_3\right)\cdot f(\mathbf{t})^2\left(\csc^2 t_k-\csc^2 t_3\right)\\
    &+f(\mathbf{t})^2\left[2\csc t_j(-\csc t_j\cot t_j)\frac{\partial t_j}{\partial t_k}-2\csc t_3(-\csc t_3\cot t_3)\frac{\partial}{\partial t_k}(\pi-t_0-t_1-t_2)\right]\\
    &=2f(\mathbf{t})^3\left(\csc^2 t_j-\csc^2 t_3\right)\left(\csc^2 t_k-\csc^2 t_3\right)\\
    &-2f(\mathbf{t})^2\left(\delta_{jk}\csc^2 t_j\cot t_j+\csc^2 t_3\cot t_3\right).
\end{align*}

It remains to prove that the $3\times 3$ Hessian matrix of $f$, denoted by $H:=(f_{jk})_{j,k=0}^2$, is negative definite. To simplify the notation, let $q_i=\cot t_i$, $r_i=\csc^2 t_i$, $S=\sum_{i=0}^3 q_i$, and $f=1/S$. Then, in this notation,
\[
f_{jk}=2f^3(r_j-r_3)(r_k-r_3)-2f^2(q_3 r_3+\delta_{jk}q_j r_j).
\]
Also, set
\[
v=\begin{pmatrix}
    r_0-r_3\\r_1-r_3\\r_2-r_3
\end{pmatrix},\quad A_{jk}:=q_3 r_3+\delta_{jk}q_j r_j,\quad\text{and}\quad A=(A_{jk})_{j,k=0}^2.
\]
The Hessian matrix of $f$ is then
\begin{equation}\label{eq:Hessian}
    H=2f^3 vv^\top-2f^2 A=-2f^2(A-fvv^\top).
\end{equation}

The remainder of the proof consists of two main steps. Namely, we will show that:
\begin{itemize}
    \item[I.)] $A$ is positive definite.

    \item[II.)] $\det(H)<0$.
\end{itemize}
Together, items I.) and II.) imply the negative definiteness of $H$, as we will show at the very end of the proof.

In the sequel, we will use the notation and definitions from the proof of Lemma \ref{lem:spherical-cone-formula}. Recall that
\[
S=\frac{\sin(t_0+t_1)\sin(t_0+t_2)\sin(t_0+t_3)}{\sin t_0\sin t_1\sin t_2\sin t_3},
\]
and we also showed already that each $\sin t_i>0$. Moreover, $0<t_i+t_j<\pi$ for $i\neq j$ (since the other $t_k$'s are positive and their total sum is $\pi$), so $\sin(t_i+t_j)>0$. Hence $S>0$ and $f>0$. Since $t_0+t_3=\pi-(t_1+t_2)$, we have $\sin(t_0+t_3)=\sin(t_1+t_2)$. Letting also
\[
\Pi:=\prod_{i=0}^3\sin t_i\quad\text{and}\quad K:=\sin(t_0+t_1)\sin(t_0+t_2)\sin(t_1+t_2),
\]
we have $S=K/\Pi$ and $f=\Pi/K$.

Next, we show that $A$ is the restriction of a diagonal quadratic form to a 3-dimensional subspace of $\R^4$. Let $d_i:=q_i r_i=\cot t_i\csc^2 t_i$ and $D:=\diag(d_0,d_1,d_2,d_3)$. Consider the tangent hyperplane
\[
L:=\{x=(x_0,x_1,x_2,x_3)\in\R^4: x_0+x_1+x_2+x_3=0\}.
\]
A basis of $L$ is $\{b_i\}_{i=0}^2$ where $b_i:=e_i-e_3$ and $\{e_i\}_{i=0}^3$ is the standard basis of $\R^4$. Then
\[
b_j^\top Db_k=(e_j-e_3)^\top D(e_k-e_3)=\delta_{jk}d_j+d_3=\delta_{jk}q_j r_j+q_3 r_3=A_{jk}.
\]
Thus, $A$ is the matrix representing the restriction of $D$ to $L$ in the basis $\{b_0,b_1,b_2\}$.

We claim that $A$ has at most one negative eigenvalue. Since $r_i=\csc^2 t_i>0$, we have
\[
\sgn(d_i)=\sgn(q_i r_i)=\sgn(q_i)=\sgn(\cot t_i).
\]
At most one of $t_0,t_1,t_2,t_3$ can exceed $\pi/2$ since $\sum_{i=0}^3 t_i=\pi$. Hence at most one $q_i$, and thus at most one $d_i$, can be negative. Thus, the diagonal quadratic form
\[
Q(x):=x^\top Dx=\sum_{i=0}^3 d_i x_i^2
\]
has a negative index of at most one. This means that $A$ has at most one negative eigenvalue, as claimed. Thus, to prove that $A$ is positive definite, it therefore suffices to show that $\det(A)>0$. Indeed, since $A$ is $3\times 3$, if it had exactly one negative eigenvalue and the other two were positive, then its determinant would be negative. We just showed that $A$ cannot have two or three negative eigenvalues.

In the next step, we compute $\det(A)$. Let $\mathbf{1}_3:=\begin{pmatrix}1\\1\\1\end{pmatrix}$. We have
\[
A=\begin{pmatrix}
    d_0+d_3 & d_3 & d_3\\ d_3 & d_1+d_3 & d_3\\d_3 & d_3 & d_2+d_3
\end{pmatrix}=\diag(d_0,d_1,d_2)+(d_3\mathbf{1}_3)\mathbf{1}_3^\top.
\]
We split the argument into two cases. First, assume that no $d_i$ equals zero. By Cauchy's Matrix Determinant Lemma,
\begin{align*}
    \det(A) &=\det(\diag(d_0,d_1,d_2))\left[1+\mathbf{1}_3^\top(\diag(d_0,d_1,d_2))^{-1}(d_3\mathbf{1}_3)\right]\\
    &=d_0 d_1 d_2\left[1+\mathbf{1}_3^\top\diag(d_0^{-1},d_1^{-1},d_2^{-1})(d_3\mathbf{1}_3)\right]\\
    &=d_0 d_1 d_2+d_0 d_1 d_2 d_3\mathbf{1}_3^\top\diag(d_0^{-1},d_1^{-1},d_2^{-1})\mathbf{1}_3\\
    &=d_0 d_1 d_2+d_3\mathbf{1}_3^\top\diag(d_1 d_2, d_0 d_2, d_0 d_1)\mathbf{1}_3\\
    &=d_0 d_1 d_2+d_3(1,1,1)\begin{pmatrix}
        d_1 d_2 & 0 & 0\\0 & d_0 d_2 & 0\\0 & 0 & d_0 d_1\end{pmatrix}
    \begin{pmatrix}
        1\\1\\1
    \end{pmatrix}\\
    &=d_0 d_1 d_2+d_3(d_1 d_2+d_0 d_2+d_0 d_1)
    =\left(\prod_{i=0}^3 d_i\right)\left(\sum_{i=0}^3 d_i^{-1}\right).
\end{align*}
Note that $\prod_{i=0}^3 d_i=\prod_{i=0}^3\frac{\cos t_i}{\sin^3 t_i}=\frac{\prod_{i=0}^3\cos t_i}{\Pi^3}$ and
\[
d_i^{-1}=\frac{\sin^3 t_i}{\cos t_i}=\frac{\sin t_i(1-\cos^2 t_i)}{\cos t_i}=\tan t_i-\sin t_i\cos t_i,
\]
so
\[
\sum_{i=0}^3 d_i^{-1}=\sum_{i=0}^3(\tan t_i-\sin t_i\cos t_i)=\underbrace{\sum_{i=0}^3\tan t_i}_{=:T_1}-\underbrace{\sum_{i=0}^3\sin t_i\cos t_i}_{=:T_2}.
\]
We now compute the sums $T_1$ and $T_2$. We begin with $T_1$. Assume that $t_i\neq\pi/2$ for all $i$, so that $T_1$ is defined. (We will see in the second case that this assumption is equivalent to no $d_i$ equaling zero, so this does not affect the logic of the proof. The exceptional cases $t_i=\pi/2$ will thus be handled in the second case.) Then $\sum_{i=0}^3 t_i=\pi$ and $\cos t_i\neq 0$ for $i\in\{0,1,2,3\}$. We claim that
\begin{equation}\label{eq:T1-identity}
    T_1=S\prod_{i=0}^3 \tan t_i.
\end{equation}
On one hand, note that $\sin(\sum_{i=0}^3 t_i)=\sin\pi=0$. On the other hand, using the trigonometric identities $\sin(x+y)=\sin x\cos y+\cos x\sin y$ and $\cos(x+y)=\cos x\cos y-\sin x\sin y$, we obtain
\begin{align*}
0 &=\sin(t_0+t_1+t_2+t_3)=\sin(t_0+t_1)\cos(t_2+t_3)+\cos(t_0+t_1)\sin(t_2+t_3)\\
&=(\sin t_0\cos t_1+\cos t_0\sin t_1)(\cos t_2\cos t_3-\sin t_2\sin t_3)\\
&+(\cos t_0\cos t_1-\sin t_0\sin t_1)(\sin t_2\cos t_3+\cos t_2\sin t_3)\\
&=\sin t_0\cos t_1\cos t_2\cos t_3+\cos t_0\sin t_1\cos t_2\cos t_3+\cos t_0\cos t_1\sin t_2\cos t_3\\
&+\cos t_0\cos t_1\cos t_2\sin t_3
-\sin t_0\sin t_1\sin t_2\cos t_3-\sin t_0\sin t_1\cos t_2\sin t_3\\
&-\sin t_0\cos t_1\sin t_2\sin t_3-\cos t_0\sin t_1\sin t_2\sin t_3.
\end{align*}
Thus
\begin{align*}
    &\phantom{=}\sin t_0\cos t_1\cos t_2\cos t_3+\cos t_0\sin t_1\cos t_2\cos t_3+\cos t_0\cos t_1\sin t_2\cos t_3\\&+\cos t_0\cos t_1\cos t_2\sin t_3\\
    &=\sin t_0\sin t_1\sin t_2\cos t_3+\sin t_0\sin t_1\cos t_2\sin t_3+\sin t_0\cos t_1\sin t_2\sin t_3\\
    &+\cos t_0\sin t_1\sin t_2\sin t_3.
\end{align*}
Dividing both sides by $\prod_{i=0}^3\cos t_i\neq 0$, we get
\begin{align*}
    T_1&=\sum_{i=0}^3\tan t_i\\
    &=\tan t_0\tan t_1\tan t_2+\tan t_0\tan t_1\tan t_3+\tan t_0\tan t_2\tan t_3+\tan t_1\tan t_2\tan t_3\\
    &=\left(\prod_{i=0}^3\tan t_i\right)\cot t_3+\left(\prod_{i=0}^3\tan t_i\right)\cot t_2+\left(\prod_{i=0}^3\tan t_i\right)\cot t_1+\left(\prod_{i=0}^3\tan t_i\right)\cot t_0\\
    &=\left(\prod_{i=0}^3\tan t_i\right)\left(\sum_{i=0}^3\cot t_i\right)
    =S\prod_{i=0}^3\tan t_i=\frac{S\Pi}{\prod_{i=0}^3\cos t_i}.
\end{align*}

Next, we consider $T_2=\sum_{i=0}^3 \sin t_i\cos t_i$. We will show that $T_2=2K$. Indeed,
\begin{align*}
    2T_2 &=\sum_{i=0}^3 2\sin t_i\cos t_i=\sum_{i=0}^3\sin(2t_i)=\sum_{i=0}^2\sin(2t_i)+\sin(2t_3)\\
    &=\sum_{i=0}^2\sin(2t_i)+\sin(2(\pi-t_0-t_1-t_2))=\sum_{i=0}^2\sin(2t_i)-\sin(2t_0+2t_1+2t_2).
\end{align*}
Note that $\sin(2t_0)+\sin(2t_1)=2\sin(t_0+t_1)\cos(t_0-t_1)$ and
\begin{align*}
    \sin(2t_2)-\sin(2t_0+2t_1+2t_2)&=2\cos(t_0+t_1+2t_2)\sin(-t_0-t_1)\\
    &=-2\sin(t_0+t_1)\cos(t_0+t_1+2t_2).
\end{align*}
Thus,
\[
2T_2=2\sin(t_0+t_1)\left[\cos(t_0-t_1)-\cos(t_0+t_1+2t_2)\right].
\]
Using now the identity $\cos x-\cos y=-2\sin\left(\frac{x+y}{2}\right)\sin\left(\frac{x-y}{2}\right)$, we get
\begin{align*}
    &\phantom{=}\cos(t_0-t_1)-\cos(t_0+t_1+2t_2)\\
    &=-2\sin\left(\frac{t_0-t_1+t_0+t_1+2t_2}{2}\right)\sin\left(\frac{t_0-t_1-(t_0+t_1+2t_2)}{2}\right)\\
    &=-2\sin(t_0+t_2)\sin(-t_1-t_2)\\
    &=2\sin(t_0+t_2)\sin(t_1+t_2).
\end{align*}
Consequently,
\[
T_2=2\sin(t_0+t_1)\sin(t_0+t_2)\sin(t_1+t_2)=2K.
\]
Therefore,
\begin{align*}
    \det(A)&=\left(\frac{\prod_{i=0}^3\cos t_i}{\Pi^3}\right)(T_1-T_2)
    =\left(\frac{\prod_{i=0}^3\cos t_i}{\Pi^3}\right)\left(\frac{S\Pi}{\prod_{i=0}^3\cos t_i}-2K\right)\\
    &=\frac{S}{\Pi^2}-\frac{2K}{\Pi^3}\left(\prod_{i=0}^3\cos t_i\right)=\frac{K}{\Pi^3}\left(1-2\prod_{i=0}^3\cos t_i\right).
\end{align*}
Since $K,\Pi>0$, it suffices to prove that $1-2\prod_{i=0}^3\cos t_i>0$. If one of the $t_i$ is at least $\pi/2$, then $\prod_{i=0}^3\cos t_i\leq 0$, so $1-2\prod_{i=0}^3\cos t_i\geq 1>0$. If, instead, all $t_i\in(0,\pi/2)$, then the function $x\mapsto\log(\cos x)$ is strictly concave on $(0,\pi/2)$ since $\frac{d^2}{dx^2}[\log(\cos x)]=-\sec^2 x<0$ on this interval. Hence by Jensen's inequality,
\[
\frac{1}{4}\sum_{i=0}^3\log(\cos t_i)\leq\log\left(\cos\left(\frac{1}{4}\sum_{i=0}^3 t_i\right)\right)=\log(\cos(\pi/4))=\log(1/\sqrt{2}).
\]
Thus,
\begin{align*}
    \frac{1}{\sqrt{2}}&=e^{\log(1/\sqrt{2})}\geq e^{\frac{1}{4}\sum_{i=0}^3\log(\cos t_i)}\\&=\left(e^{\sum_{i=0}^3\log(\cos t_i)}\right)^{1/4}
    =\left(\prod_{i=0}^3 e^{\log(\cos t_i)}\right)^{1/4}=\left(\prod_{i=0}^3\cos t_i\right)^{1/4},
\end{align*}
which implies $\prod_{i=0}^3\cos t_i\leq 1/4$. Therefore, in the case all $t_i\in(0,\pi/2)$, we have
\[
1-2\prod_{i=0}^3\cos t_i\geq 1-2(1/4)=1/2>0.
\]
This proves that $\det(A)>0$, and thus that $A$ is positive definite, in the case that all $d_i\neq 0$.

To complete the proof of I.), we now handle the exceptional cases $d_i=0$ and $t_i=\pi/2$. Recall that $d_i=\cot t_i\csc^2 t_i$. Since $t_i\in(0,\pi)$, we have $\csc^2 t_i>0$, so
\[
d_i=0\Longleftrightarrow \cot t_i=0\Longleftrightarrow t_i=\pi/2.
\]
Therefore, the exceptional cases are equivalent. Since each $t_i>0$ and $\sum_{i=0}^3 t_i=\pi$, at most one $t_i$ can equal $\pi/2$, say, $t_k=\pi/2$. Then $\sum_{i\neq k}t_i=\pi/2$, so $0<t_i<\pi/2$ for all $i\neq k$. Hence $\cot t_i>0$ for $i\neq k$, so $d_i>0$ for all $i\neq k$, while $d_k=0$. Recall that $L=\{x\in\R^4: \sum_{i=0}^3 x_i=0\}$ and let $x\in L\setminus\{o\}$. Then
\[
Q(x)=\sum_{i=0}^3 d_i x_i^2=\sum_{i\neq k}d_i x_i^2\geq 0.
\]
If $Q(x)=0$, then each $x_i=0$ for $i\neq k$ since $d_i>0$ for $i\neq k$, so $x_k=0$ too by the definition of $L$; this implies $x=o$, a contradiction. Thus $Q(x)>0$ for all $x\in L\setminus\{o\}$, so the matrix $A$, which represents $Q|_L$, is positive definite. This completes the proof of I.).

We now turn our attention to proving II.) and show that $\det(H)<0$. Recall that $D=\diag(d_0,d_1,d_2,d_3)$, and set $R=(r_0,r_1,r_2,r_3)^\top$. Let $E=D-fRR^\top$, and restrict $E$ to the hyperplane $L$ in the basis $b_i=e_i-e_3$, $i\in\{0,1,2\}$. Then $b_i^\top Db_j=A_{ij}$ and $b_i^\top R=r_i-r_3=v_i$. Thus, the restricted matrix is $M=A-fvv^\top$. Hence $H=-2f^2(A-fvv^\top)=-2f^2 M$, which implies
\begin{equation}\label{eq:relation-determinants-of-H-and-M}
    \det(H)=\det(-2f^2 M)=(-2f^2)^3\det(M)=-8f^6\det(M).
\end{equation}
Since $f>0$, $\sgn(\det(H))=-\sgn(\det(M))$, so our goal becomes to show that $\det(M)>0$.

We first consider the case that $t_i\neq \pi/2$ for every $i\in\{0,1,2,3\}$. Then $d_i\neq 0$ for every $i$, so $D$ is invertible. Set
\[
a:=\sum_{i=0}^3d_i^{-1},\qquad
b:=\sum_{i=0}^3\frac{r_i}{d_i},\qquad
c:=\sum_{i=0}^3\frac{r_i^2}{d_i}.
\]
Since $d_i=q_ir_i$, we have
\[
b=\sum_{i=0}^3q_i^{-1}=\sum_{i=0}^3\tan t_i
\]
and
\[
c=\sum_{i=0}^3\frac{r_i}{q_i}
=\sum_{i=0}^3\left(q_i+q_i^{-1}\right)
=S+b.
\]
Since $f=1/S$, it follows that
\[
1-fc=1-\frac{1}{S}(S+b)=-\frac{b}{S}=-bf.
\]
Moreover, using $\sum_{i=0}^3t_i=\pi$, we have
\[
b=\sum_{i=0}^3\tan t_i
=S\prod_{i=0}^3\tan t_i.
\]
Since $t_i\neq\pi/2$ for every $i$ and $t_i\in(0,\pi)$, the product on the
right is nonzero. Hence $b\neq 0$, and therefore $1-fc=-bf\neq 0$. By Cauchy's Matrix Determinant Lemma,
\begin{align*}
    \det(E)&=\det(D+(-fR)R^\top)=\det(D)\left[1+R^\top D^{-1}(-fR)\right]=\det(D)\left[1-fR^\top D^{-1}R\right].
\end{align*}
Computing the second term inside the brackets, we get
\begin{align*}
R^\top D^{-1}R&=(r_0,r_1,r_2,r_3)
    \begin{pmatrix}
        d_0^{-1} & 0 & 0 & 0\\
        0 & d_1^{-1} & 0 & 0\\
        0 & 0 & d_2^{-1} & 0\\
        0 & 0 & 0 & d_3^{-1}
    \end{pmatrix}\begin{pmatrix}
        r_0\\r_1\\r_2\\r_3
    \end{pmatrix}
    =\sum_{i=0}^3\frac{r_i^2}{d_i}=c.
\end{align*}
Hence
\[
\det(E)=\det(D)\left(1-fR^\top D^{-1}R\right)
       =\det(D)(1-fc)\neq 0,
\]
so $E$ is invertible.

Denoting the basis matrix by \[B=\begin{pmatrix}
1 & 0 & 0\\0& 1 & 0\\0 & 0 & 1\\-1 & -1 & -1\end{pmatrix}\]
we obtain
\[
B^\top EB=B^\top(D-fRR^\top)B=B^\top DB-fB^\top RR^\top B=A-fvv^\top=M.
\]
In the next step of the proof, we aim to  show that
\begin{equation}\label{eq:det-M-det-E-relation}
\det(M)=\det(B^\top EB)=\det(E)\mathbf{1}_4^\top E^{-1}\mathbf{1}_4.
\end{equation}
Let $P_0=\begin{pmatrix}
    b_0 & b_1 & b_2 & \tfrac{1}{4}\mathbf{1}_4
\end{pmatrix}$.
Since
\[
\det(P_0)=\det\begin{pmatrix}
    1 & 0 & 0 &1/4\\
    0 & 1 & 0 & 1/4\\
    0 & 0 & 1 & 1/4\\
    -1 & -1 & -1 & 1/4
\end{pmatrix}=1,
\]
we have $\det(P_0^\top EP_0)=\det(E)$. Note that \[P_0^{-1}=\begin{pmatrix}3/4 & -1/4 & -1/4 &-1/4\\-1/4 & 3/4 & -1/4 & -1/4\\-1/4 & -1/4 & 3/4 & -1/4 \\1 & 1 & 1 & 1\end{pmatrix}\] has final row $\mathbf{1}_4^\top$. Thus, the lower right entry of $(P_0^\top EP_0)^{-1}=P_0^{-1}E^{-1}P_0^{-\top}$ is $\mathbf{1}_4^\top E^{-1}\mathbf{1}_4$. We also have
\[
C=P_0^\top EP_0=\begin{pmatrix}
    B^\top\\\frac{1}{4}\mathbf{1}_4^\top
\end{pmatrix}E\begin{pmatrix}
    B & \frac{1}{4}\mathbf{1}_4\end{pmatrix}=\begin{pmatrix}
        B^\top EB & \frac{1}{4}B^\top E\mathbf{1}_4\\
        \frac{1}{4}\mathbf{1}_4^\top EB & \frac{1}{16}\mathbf{1}_4^\top E\mathbf{1}_4
    \end{pmatrix}.
\]
The upper left block of the $4\times 4$ matrix $C$ is the $3\times 3$ matrix $B^\top EB$. Using the formula $C^{-1}=\operatorname{adj}(C)/\det(C)$, we have $(C^{-1})_{ij}=\operatorname{cof}_{ji}(C)/\det(C)$ where $\operatorname{cof}_{ji}(C)$ is the $(j,i)$-cofactor of $C$. In particular, the lower right entry of $C^{-1}$ is
\[
(C^{-1})_{44}=\frac{\operatorname{cof}_{44}(C)}{\det(C)}=\frac{(-1)^{4+4}\det(C_{[4|4]})}{\det(C)}=\frac{\det(B^\top EB)}{\det(P_0^\top EP_0)}
\]
where $C_{[4|4]}$ is the minor obtained by deleting the last row and column of $C$. Therefore, $(C^{-1})_{44}=\mathbf{1}_4^\top E^{-1}\mathbf{1}_4=\det(B^\top EB)/\det(P_0^\top EP_0)$, so
\[
\det(B^\top EB)=\det(P_0^\top EP_0)\mathbf{1}_4^\top E^{-1}\mathbf{1}_4=\det(E)\mathbf{1}_4^\top E^{-1}\mathbf{1}_4.
\]
This proves \eqref{eq:det-M-det-E-relation}.

Next, we apply \eqref{eq:det-M-det-E-relation} to the matrix $E=D-fRR^\top$. By the Sherman--Morrison formula (see, e.g., \cite[Equation (2.1.5), p. 65]{Golub-VanLoan-matrix-computations}),
\[
E^{-1}=(D+(-fR)R^\top)^{-1}=D^{-1}-\frac{D^{-1}(-fR)R^\top D^{-1}}{1+R^\top D^{-1}(-fR)}=D^{-1}+\frac{fD^{-1}RR^\top D^{-1}}{1-fc}.
\]
Thus,
\begin{align*}
    \mathbf{1}_4^\top E^{-1}\mathbf{1}_4&=\mathbf{1}_4^\top\left(D^{-1}+\frac{fD^{-1}RR^\top D^{-1}}{1-fc}\right)\mathbf{1}_4=\mathbf{1}_4^\top D^{-1}\mathbf{1}_4+\mathbf{1}_4^\top\left(\frac{fD^{-1}RR^\top D^{-1}}{1-fc}\right)\mathbf{1}_4\\
    &=\underbrace{\mathbf{1}_4^\top D^{-1}\mathbf{1}_4}_{=a}+\frac{f}{1-fc}(\underbrace{\mathbf{1}_4^\top D^{-1}R}_{=b})(\underbrace{R^\top D^{-1}\mathbf{1}_4}_{=b})=a+\frac{fb^2}{1-fc}.
\end{align*}
Therefore, by \eqref{eq:det-M-det-E-relation},
\begin{equation}\label{eq:detM-penultimate}
    \det(M)=\det(D)(1-fc)\left(a+\frac{fb^2}{1-fc}\right)=\det(D)\left[a(1-fc)+fb^2\right].
\end{equation}

Using $1-fc=-bf$ in \eqref{eq:detM-penultimate}, we get
\[
\det(M)=\det(D)\left(fb^2-abf\right)=fb\det(D)(b-a).
\]
Note that
\begin{align*}
b-a&=\sum_{i=0}^3 \frac{1}{q_i}-\sum_{i=0}^3\frac{1}{q_i r_i}=\sum_{i=0}^3\frac{r_i-1}{q_i r_i}=\sum_{i=0}^3\frac{q_i^2}{q_i r_i}=\sum_{i=0}^3\frac{q_i}{1+q_i^2}=\sum_{i=0}^3\frac{\cot t_i}{1+\cot^2 t_i}\\
&=\sum_{i=0}^3\frac{\cos t_i/\sin t_i}{\csc^2 t_i}=\sum_{i=0}^3\sin t_i\cos t_i=2K,
\end{align*}
$\det(D)=\prod_{i=0}^3 q_ir_i$, and $fb=fS\prod_{i=0}^3\tan t_i=\prod_{i=0}^3 \tan t_i$. Since $q_i\tan t_i=1$, we thus obtain
\[
fb\det(D)=\prod_{i=0}^3 q_i r_i\tan t_i=\prod_{i=0}^3 r_i=\prod_{i=0}^3\csc^2 t_i=1/\Pi^2.
\]
Therefore, $\det(M)=2K/\Pi^2>0$, which implies
\begin{align*}
    \det(H)&=-8f^6\det(M)=-\frac{16Kf^6}{\Pi^2}=-\frac{16K}{\Pi^2}\left(\frac{\Pi}{K}\right)^6=-\frac{16\Pi^4}{K^5}
   \\ &=-\frac{16\sin^4 t_0\sin^4 t_1\sin^4 t_2\sin^4 t_3}{\sin^5(t_0+t_1)\sin^5(t_0+t_2)\sin^5(t_1+t_2)}.
\end{align*}
As we proved earlier, every sine term here is positive, so $\det(H)<0$
whenever $t_i\neq\pi/2$ for every $i$.

It remains to consider the exceptional case that $t_k=\pi/2$ for some
$k\in\{0,1,2,3\}$. The matrix $H$ depends continuously on
$(t_0,t_1,t_2,t_3)$ on the set
\[
\left\{(t_0,t_1,t_2,t_3)\in(0,\pi)^4:
\sum_{i=0}^3t_i=\pi\right\}.
\]
The points for which $t_i\neq\pi/2$ for every $i$ are dense in this set.
Therefore, by continuity, the identity $\det(H)=-16\Pi^4/K^5$
also holds when $t_k=\pi/2$ for some $k$. Since $\Pi>0$ and $K>0$, we
again obtain $\det(H)<0$. This proves II.).

To complete the proof, we show that I.) and II.) together imply that $H$ is negative definite. Recall that $H=-2f^2 A+2f^3 vv^\top$. We already proved that $A>0$. Consider any $x\in v^\perp$. Then $v^\top x=0$, and hence
\[
x^\top Hx=-2f^2 x^\top Ax+2f^3(v^\top x)^2=-2f^2 x^\top Ax<0
\]
for every nonzero $x\in v^\perp$. If $v=o$, then $H=-2f^2 A$ is already negative definite. So assume that $v\neq o$. Let $\lambda_1,\lambda_2,\lambda_3$ be the eigenvalues of $H$. At least two are negative, say, $\lambda_1,\lambda_2<0$. But $\lambda_1\lambda_2\lambda_3=\det(H)<0$. Since $\lambda_1\lambda_2>0$, it follows that $\lambda_3<0$. Therefore, all three eigenvalues of $H$ are negative. This proves that $H$ is negative definite. \qed

\section{Degree sequences for ten and eleven vertices}\label{sec:degree-sequence-tables}

Tables \ref{tab:N10-degree-sequences} and \ref{tab:N11-degree-sequences} list the degree sequences of all simplicial convex 3-polytopes with exactly ten and eleven vertices, respectively, and no trivalent vertices. Here $d(P)$ lists the vertex degrees in nonincreasing order, and exponents denote multiplicities. Repeated degree sequences correspond to distinct combinatorial types. The labels Types 11 and 12 for ten vertices, and Types 32, 33, and 34 for eleven vertices, agree with those used in the proofs; Types 32 and 33 are distinguished by their facial degree patterns in Table \ref{tab:facial-degree-patterns}.

The tables are obtained by enumerating the simple planar triangulations of minimum degree at least 4 using \texttt{plantri} \cite{BrinkmannMcKay,plantri58}. Steinitz's theorem identifies these graphs with the required combinatorial types. The supplied graph lists contain twelve graphs for ten vertices and thirty-four for eleven vertices. Their degree sequences give the tables below; in particular, the eleven-vertex list has sixteen distinct degree sequences.

\begin{table}[p]
\centering
\renewcommand{\arraystretch}{1.12}
\[
\begin{array}{c|c}
\text{Type} & d(P) \\
\hline
1 &(8^2,4^8)\\
2 &(7^2,5^2,4^6)\\
3 &(7^2,5^2,4^6)\\
4 &(7,6^2,5,4^6)\\
5 &(7,6,5^3,4^5)\\
6 &(6^4,4^6)\\
7 &(6^3,5^2,4^5)\\
8 &(6^2,5^4,4^4)\\
9 &(6^2,5^4,4^4)\\
10 &(6^2,5^4,4^4)\\
11 &(6,5^6,4^3)\\
12 &(5^8,4^2)
\end{array}
\]
\caption{Degree sequences of the twelve ten-vertex simplicial combinatorial types with no trivalent vertices.}
\label{tab:N10-degree-sequences}
\end{table}

\begin{table}[p]
\centering
\renewcommand{\arraystretch}{1.12}
\[
\begin{array}{c|c}
\text{\rm Type(s)} & \text{\rm Degree sequence } d(P)\\
\hline
1 &(9^2,4^9)\\
2\text{--}3 &(8^2,5^2,4^7)\\
4\text{--}5 &(8,7,5^3,4^6)\\
6 &(8,6^3,4^7)\\
7\text{--}9 &(8,6^2,5^2,4^6)\\
10 &(7^2,6^2,4^7)\\
11\text{--}14 &(7^2,6,5^2,4^6)\\
15\text{--}16 &(7^2,5^4,4^5)\\
17 &(7,6^3,5,4^6)\\
18\text{--}21 &(7,6^2,5^3,4^5)\\
22\text{--}24 &(7,6,5^5,4^4)\\
25 &(6^5,4^6)\\
26\text{--}27 &(6^4,5^2,4^5)\\
28\text{--}31 &(6^3,5^4,4^4)\\
32\text{--}33 &(6^2,5^6,4^3)\\
34 &(6,5^8,4^2)
\end{array}
\]
\caption{Degree sequences of the thirty-four eleven-vertex simplicial combinatorial types with no trivalent vertices.}
\label{tab:N11-degree-sequences}
\end{table}
\clearpage

\section{Computational certificates}\label{sec:certificates}

\subsection{Exact inputs and interval verification}

The computational parts of the proof consist of strict inequalities at explicit supporting points. Note that none of the exclusions relies on a numerical optimizer finding a global optimum. Every displayed terminating decimal used as an input is interpreted exactly. The supplied Python programs use only the standard library and perform interval arithmetic with $80$ decimal digits, rounding the endpoints outward after every operation. Sine and cosine are evaluated by $100$ Taylor terms, with a tail enclosure of $10^{-100}$ on the interval $|x|<4$. The value of $\pi$ and every radical are enclosed by intervals. Each interval used in a strict comparison is separated from the comparison threshold.

For a degree certificate, the verifier checks the supporting slope and intercept of $\Psi_m$ and applies \eqref{eq:degree-supporting-line-bound}. For a facial certificate, it checks membership of the supporting point in the domain of $W$, the three derivative bounds, the intercept bound, and the resulting sum in Lemma \ref{lem:tangent-plane-bound}. For a vertex-dependent certificate, it also checks that the listed triangles are exactly the facets of the indicated graph. The lower-bound configurations are checked by strict supporting plane determinants, incidence counts, and volume determinants after normalization.

The graph lists are generated by \texttt{plantri}~5.8, with commands of the form
\begin{verbatim}
    plantri -m4a 13 n13.txt
    plantri -m4a 14 n14.txt
    plantri -m4a 15 n15.txt
    plantri -m4a 16 n16.txt
\end{verbatim}
The verification checks every graph in these lists. Completeness of the enumeration relies on \texttt{plantri} \cite{BrinkmannMcKay,plantri58}; the unmodified program source, guide, and license are included. This is a computational proof with explicit certificates, rather than a formal verification of the graph generator or the decimal arithmetic implementation. Record identifiers refer to the supplied lists and start at zero.

\subsection{Verification programs}\label{sec:appendix-verification-programs}

The accompanying certificate archive contains the exact inputs, graph lists, coordinate data, verification programs, and the \texttt{plantri} source, guide, and license. First, extract this archive to a directory of your choice. The certificate files are not needed to compile the paper. Next, change to the extracted directory and run Python with assertions enabled:
\begin{verbatim}
    python3 verify_small.py
    python3 certify.py
    python3 verify_all.py
    python3 verify_larger.py
    python3 certify32.py
\end{verbatim}
The first program verifies the small-case facial data, the retained ten- and eleven-vertex tangent plane inequalities, and the explicit eleven-vertex lower bound. The next programs verify the two displayed thirteen-vertex tables, every exclusion in Table \ref{tab:larger-reductions}, the explicit larger candidates, and the universal thirty-two-vertex facial inequality. The thirty-two-vertex proof uses no graph enumeration.

\subsection{The vertex-dependent fifteen-vertex certificate}\label{subsec:N15-vertex-data}

The vertex labels are $0,1,\ldots,14$ in record $0$ of \texttt{n15.txt}. Each row gives a facet, its supporting half-angle triple in the same vertex order, and its intercept upper bound. All entries are exact decimals. The coefficients are:
\begin{align*}
    \lambda_{0}=0.2721560,\quad \lambda_{1}=0.2718420,\quad \lambda_{2}=0.2721560\\
    \lambda_{3}=0.2939279,\quad \lambda_{4}=0.2747870,\quad \lambda_{5}=0.2748625\\
    \lambda_{6}=0.2901189,\quad \lambda_{7}=0.2748625,\quad \lambda_{8}=0.2747870\\
    \lambda_{9}=0.2741029,\quad \lambda_{10}=0.2741029,\quad \lambda_{11}=0.2901529\\
    \lambda_{12}=0.2551400,\quad \lambda_{13}=0.2901529,\quad \lambda_{14}=0.2688174.
\end{align*}
The interval verifier checks $W_{\theta_j}(s_F)<\lambda_{v_j}$ and $C_F<\overline C_F$ at every row. Summation gives \eqref{eq:N15-final-exclusion}.

\begin{table}[!htbp]
\centering\footnotesize
\renewcommand{\arraystretch}{1.05}
\begin{tabular}{c|c|c}
Facet & $s_F$ & $\overline C_F$\\
\hline
$(0,1,2)$ & $(0.60727494,0.60930138,0.60727494)$ & $-0.38422287$\\
$(0,1,5)$ & $(0.61260660,0.61472946,0.59552117)$ & $-0.38585044$\\
$(0,2,3)$ & $(0.64617414,0.64617414,0.51894948)$ & $-0.39665980$\\
$(0,3,4)$ & $(0.65641048,0.52272398,0.63173960)$ & $-0.39834076$\\
$(0,4,5)$ & $(0.61912648,0.60162208,0.60115016)$ & $-0.38764144$\\
$(1,2,7)$ & $(0.61472946,0.61260660,0.59552117)$ & $-0.38585044$\\
$(1,5,6)$ & $(0.65141618,0.62499973,0.53725961)$ & $-0.39616414$\\
$(1,6,7)$ & $(0.65141618,0.53725961,0.62499973)$ & $-0.39616414$\\
$(2,3,8)$ & $(0.65641048,0.52272398,0.63173960)$ & $-0.39834076$\\
$(2,7,8)$ & $(0.61912648,0.60115016,0.60162208)$ & $-0.38764144$\\
$(3,4,10)$ & $(0.52608706,0.63916825,0.64551004)$ & $-0.39960805$\\
$(3,8,9)$ & $(0.52608706,0.63916825,0.64551004)$ & $-0.39960805$\\
$(3,9,10)$ & $(0.52502111,0.64291759,0.64291759)$ & $-0.39916956$\\
$(4,5,11)$ & $(0.63582186,0.63517207,0.54236131)$ & $-0.39807784$\\
$(4,10,11)$ & $(0.63324087,0.63915003,0.54107542)$ & $-0.39759390$\\
$(5,6,11)$ & $(0.68474953,0.56131912,0.56113583)$ & $-0.40724292$\\
$(6,7,13)$ & $(0.56131912,0.68474953,0.56113583)$ & $-0.40724292$\\
$(6,11,12)$ & $(0.47221760,0.47209782,0.84278978)$ & $-0.39234019$\\
$(6,12,13)$ & $(0.47221760,0.84278978,0.47209782)$ & $-0.39234019$\\
$(7,8,13)$ & $(0.63517207,0.63582186,0.54236131)$ & $-0.39807784$\\
$(8,9,13)$ & $(0.63324087,0.63915003,0.54107542)$ & $-0.39759390$\\
$(9,10,14)$ & $(0.59587814,0.59587814,0.63143039)$ & $-0.38468970$\\
$(9,13,14)$ & $(0.61813685,0.53025757,0.66545572)$ & $-0.39371883$\\
$(10,11,14)$ & $(0.61813685,0.53025757,0.66545572)$ & $-0.39371883$\\
$(11,12,14)$ & $(0.49466470,0.72800655,0.58962541)$ & $-0.38100479$\\
$(12,13,14)$ & $(0.72800655,0.49466470,0.58962541)$ & $-0.38100479$\\
\end{tabular}
\caption{The final fifteen-vertex exclusion using a coefficient at each vertex.}
\label{tab:N15-vertex-certificate}
\end{table}

\begin{table}[p]
\centering\small
\setlength{\tabcolsep}{11pt}
\renewcommand{\arraystretch}{1.08}
\begin{tabular}{@{}lrrr@{}}
\hline
Degrees&$q_1$&$q_2$&$q_3$\\
\hline
(3,3,3) & 0.56238015 & 0.56238015 & 0.56238015 \\
(3,3,4) & 0.59809178 & 0.59809178 & 0.48802708 \\
(3,3,5) & 0.64582277 & 0.64582277 & 0.38105704 \\
(3,3,6) & 0.66175407 & 0.66175407 & 0.34291752 \\
(3,3,7) & 0.65768786 & 0.65768786 & 0.35278203 \\
(3,4,4) & 0.66191094 & 0.51010571 & 0.51010571 \\
(3,4,5) & 0.78700441 & 0.50136984 & 0.37975608 \\
(3,4,6) & 0.86414906 & 0.46980112 & 0.32380945 \\
(3,4,7) & 0.84049657 & 0.48090491 & 0.33953128 \\
(3,5,5) & \multicolumn{3}{c}{$\mathcal B_3$} \\
(3,5,6) & \multicolumn{3}{c}{$\mathcal B_3$} \\
(3,5,7) & \multicolumn{3}{c}{$\mathcal B_3$} \\
(3,6,6) & \multicolumn{3}{c}{$\mathcal B_3$} \\
(3,6,7) & \multicolumn{3}{c}{$\mathcal B_3$} \\
(3,7,7) & \multicolumn{3}{c}{$\mathcal B_3$} \\
(4,4,4) & 0.56126296 & 0.56126296 & 0.56126296 \\
(4,4,5) & 0.62050090 & 0.62050090 & 0.43572993 \\
(4,4,6) & 0.63873555 & 0.63873555 & 0.39437392 \\
(4,4,7) & 0.63414423 & 0.63414423 & 0.40492923 \\
(4,5,5) & 0.79528005 & 0.43607541 & 0.43607541 \\
(4,5,6) & 0.90041429 & 0.39685350 & 0.35838050 \\
(4,5,7) & 0.86738451 & 0.41123779 & 0.38091938 \\
(4,6,6) & \multicolumn{3}{c}{$\mathcal B_4$} \\
(4,6,7) & 1.10012980 & 0.26240973 & 0.26907470 \\
(4,7,7) & 1.00295222 & 0.32041201 & 0.32041201 \\
(5,5,5) & 0.55939578 & 0.55939578 & 0.55939578 \\
(5,5,6) & 0.58659279 & 0.58659279 & 0.50336490 \\
(5,5,7) & 0.57999970 & 0.57999970 & 0.51714216 \\
(5,6,7) & 0.61783180 & 0.51919407 & 0.53897185 \\
(5,7,7) & 0.60838542 & 0.53405324 & 0.53405324 \\
(6,6,6) & 0.55869846 & 0.55869846 & 0.55869846 \\
(6,6,7) & 0.55044127 & 0.55044127 & 0.57531672 \\
(6,7,7) & 0.54316195 & 0.56661165 & 0.56661165 \\
(7,7,7) & 0.55887870 & 0.55887870 & 0.55887870 \\
\hline
\end{tabular}
\caption{Supporting points for the thirty-four non-$(566)$ degree triples in Lemma \ref{lem:N32-universal-face}. All decimals are exact.}
\label{tab:N32-supporting-points}
\end{table}
\clearpage
\section*{Acknowledgments}

The author sincerely thanks Lienhard Wimmer for bringing his dissertation to the author's attention, and for the correspondence on maximum convex hulls. The author used ChatGPT-6.1 Sol as a tool for technical and editorial assistance, including the generation of the code used in the numerical certifications. The main ideas and approaches used, namely the generalized degree sequence bound and the tangent plane method,  originated with the author. All mathematical statements and proofs were verified by the author.

\section*{Data availability statement}

The exact computational certificates, graph lists, coordinate data, and standard library verification programs used in this paper are supplied in the accompanying archive \nolinkurl{Maximum_Volume_Certificates.zip}, available at \url{https://zenodo.org/records/24013978}. Appendix \ref{sec:certificates} explains how to verify them.

\bibliographystyle{plain}
\bibliography{main}

\vspace{3mm}

\noindent {\sc Department of Mathematics \& Computer Science, Longwood University, 201 High Street, Farmville, VA 23909, U.S.A.}

\noindent {\it E-mail address:} {\tt hoehnersd@longwood.edu}

\end{document}